\documentclass[12pt]{amsart}

\usepackage{setspace}
\usepackage{amssymb}   
\usepackage{amsmath}
\usepackage{mathrsfs}
\usepackage{stmaryrd}
\usepackage{amsthm}
\usepackage{newlfont}
\usepackage{amscd}
\usepackage{mathtools}
\usepackage{bm}
\usepackage{tikz-cd}
\usepackage{graphicx}
\usepackage{hyperref}
\usepackage{enumitem}
\usepackage[all]{xy}
\usepackage{textcomp}
\usepackage{amsbsy}

\newtheorem{theorem}{Theorem}[section]
\newtheorem{thm-defn}[theorem]{Theorem/Definition}
\newtheorem{lemma}[theorem]{Lemma}
\newtheorem{proposition}[theorem]{Proposition}
\newtheorem{corollary}[theorem]{Corollary}

\theoremstyle{definition}
\newtheorem{defn}[theorem]{Definition}

\newtheorem{example}[theorem]{Example}
\newtheorem{ques}[theorem]{Question}
\theoremstyle{remark}
\newtheorem{rem}[theorem]{Remark}
\newtheorem{rems}[theorem]{Remark}

\newcommand{\wt}{\widetilde}
\newcommand{\bbA}{{\mathbb A}}

\newcommand{\bbL}{{\mathbb L}}

\newcommand{\bbQ}{{\mathbb Q}}
\newcommand{\bbR}{{\mathbb R}}

\newcommand{\bbZ}{{\mathbb Z}}

\newcommand{\cB}{{\mathcal B}}
\newcommand{\cC}{{\mathcal C}}

\newcommand{\cF}{{\mathcal F}}
\newcommand{\cG}{{\mathcal G}}
\newcommand{\cI}{{\mathcal I}}

\newcommand{\cJ}{{\mathcal J}}

\newcommand{\cO}{{\mathcal O}}

\newcommand{\cS}{{\mathcal S}}

\newcommand{\cX}{{\mathcal X}}

\newcommand{\rD}{{\mathrm D}}
\newcommand{\rE}{{\mathrm E}}
\newcommand{\rF}{{\mathrm F}}

\newcommand{\rH}{{\mathrm H}}

\newcommand{\rM}{{\mathrm M}}

\newcommand{\Spec}{{\mathrm Spec}}

\DeclareMathOperator{\im}{Im}
\DeclareMathOperator{\coker}{coker}
\DeclareMathOperator{\Tor}{Tor}

\DeclareMathOperator{\Gal}{Gal}
\DeclareMathOperator{\gr}{gr}

\newcommand{\onto}{\twoheadrightarrow}

\newcommand{\cris}{\textnormal{cris}}
\newcommand{\CRIS}{\textnormal{CRIS}}
\newcommand{\syn}{\textnormal{syn}}
\newcommand{\SYN}{\textnormal{SYN}}
\newcommand{\et}{\textnormal{\'et}}
\newcommand{\Rbar}{{\overline{R}}}

\newcommand{\Z}{{\mathbb Z}}
\newcommand{\inj}{{\hookrightarrow}}
\newcommand{\Kbar}{{\overline K}}
\newcommand{\MF}{{\mathrm{MF}}}

\newcommand{\R}{\mathrm}
\newcommand{\fm}{\mathfrak{m}}

\newcommand{\wh}{\widehat}

\newcommand{\PD}{\textnormal{PD}}
\newcommand{\Xbar}{\overline{X}}
\newcommand{\zar}{\mathrm{zar}}
\newcommand{\Mbar}{\overline{M}}

\begin{document}
\title{Relative Fontaine--Messing theory over power series rings}
\author{Tong Liu, Yong Suk Moon, and Deepam Patel}

\begin{abstract}
Let $k$ be a perfect field of characteristic $p>2$, $R \coloneqq W(k)[\![t_1, \dots, t_d]\!]$ be the power series ring over the Witt vectors, and $X$ be a smooth proper scheme over $R$. The main goal of this article is to extend classical Fontaine-Messing theory (\cite{FontaineMessing}) to the setting where the base ring is $R$. In particular, we obtain comparison theorems between torsion crystalline cohomology of $X/R$ and torsion \'etale cohomology in this setting.
\end{abstract}

\maketitle

\section{Introduction}

Let $k$ be a perfect field of characteristic $p>2$, $W(k)$ the Witt vectors over $k$, $R= W(k)[\![t_1, \dots, t_d]\!]$ the power series ring over $W(k)$ in $d$-variables, and $X$ a smooth proper scheme over $R$. The aim of this paper is to extend classical Fontaine--Messing theory as in \cite{FontaineMessing} to the setting where the base is $R$.

\subsection{Classical Fontaine--Messing theory}

We begin by briefly recalling one of the main results of classical Fontaine--Messing theory (\cite{FontaineMessing}). To simplify notation, we shall use $\varphi$ to denote various Frobenii (on $W(k)$-algebras) extending the arithmetic Frobenius on $W(k)$, i.e., $\varphi (x) = x^p \mod p$ for $x \in W(k)$. Let $K \coloneqq W(k)[\frac 1 p]$, $\Kbar$ be a fixed algebraic closure of $K$ and $G_K \coloneqq \Gal (\Kbar/ K)$. A key object in classical Fontaine--Messing theory (and Fontaine--Laffaille theory in \cite{Fonatine-Laffaille}) is a \emph{Fontaine--Laffaille module}

\begin{defn} A \emph{Fontaine--Laffaille module over $W(k)$} is a triple $(M, \rF^i (M), \varphi_i)$ where
\begin{enumerate}	
\item $M$ is a finitely generated $W_n(k)$-module for some $n \geq 1$.
\item $\rF^i (M)\subset M$ is a $W_{n}(k)$-submodule satisfying:
	\begin{enumerate}
	    \item $\rF^{i +1} (M)\subset \rF^{i} (M)$ for all $i$,
		\item $\rF^0 (M) = M$ and $\rF^{r+1}(M)= \{0\}$ for some $r$.
	\end{enumerate}
\item $\varphi_i\colon \rF^i(M) \to M$ is a $\varphi_{W(k)}$-semilinear map such that
	\begin{enumerate}
		\item $\varphi_i|_{\rF^{i+1}(M)} = p \varphi_{i+1}$,
		\item $\sum\limits_{i =0}^{r} \varphi_i (\rF^i (M))=M$ as $W(k)$-modules.
	\end{enumerate}
\end{enumerate}
\end{defn}

We denote by $\MF(W(k))$ the category whose objects are Fontaine--Laffaille modules over $W(k)$ and morphisms are $W(k)$-module morphisms compatible with $\rF^i$ and $\varphi_i$. Let $\MF^{[0, r]}(W(k))$  denote the full subcategory consisting of Fontaine--Laffaille modules such that $\rF^{r+1}(M) = \{0\}$. These categories are known to satisfy many good properties including:
\begin{enumerate}
\item[(1)] By \cite[Prop. 1.8]{Fonatine-Laffaille}, $\MF^{[0, r]}(W(k))$ is an abelian category.
\item[(2)] A sequence $0 \to M \to M'' \to M' \to 0$ in $\MF(W(k))$ is exact if and only if $0 \to \rF^i(M ) \to \rF^i (M'') \to \rF^i (M' ) \to 0$ is exact for all $i$. 
\end{enumerate}

Let $A_\cris$ denote the integral crystalline period ring over $W(k)$ equipped with the usual Frobenius $\varphi$ and filtration $\rF^i A_{\cris} \subset A_{\cris}$ (see \S \ref{subsec-period-ring} for details). For $i = 0, \ldots, p-1$, we have $\varphi (\rF^i A_\cris )\subset p^iA_\cris$. In particular, since $A_\cris$ is $\Z_p$-flat, we can define $\varphi_i\colon \rF^ i A_\cris\to A_\cris$ by setting $\varphi_i (x) = \frac{\varphi(x)}{p ^i}$. Given an object $(M, \rF^i(M), \varphi_i) \in \MF^{[0, r]}(W(k))$ with $r \leq p-1$, let $$\rF^r(A_\cris \otimes_{W(k)}M) \coloneqq \sum_{i = 0}^r \rF^i A_\cris \otimes_{W(k)} \rF^{r-i}M $$ and define $\varphi_r\colon \rF^r(A_\cris \otimes_{W(k)}M)  \to A_\cris \otimes_{W(k)} M$ via $\varphi_r \coloneqq \sum_{i = 0}^r \varphi_i|_{\rF^i A_\cris} \otimes \varphi_{r- i} |_{\rF^{r-i}M}$. Set 
$$T_\cris(M) \coloneqq ( \rF^{r}( A_\cris \otimes_{W(k)} M) )^{\varphi _r =1} (-r)=  \{ x \in \rF^{r}( A_\cris \otimes_{W(k)} M) ~|~ \varphi _r (x) = x\}(-r).$$ 
It can be shown that $T_\cris$ is independent of the choice of $r$ and gives rise to a functor from $p^n$-torsion Fontaine--Laffaille modules to $\Z/p^n\Z[G_K]$-modules. The main result of Fontaine--Laffaille theory states that if $r\leq p -2$, then $T_\cris$ is a fully faithful exact functor such that the essential image is stable under subquotients.

Given a smooth proper scheme $X$ over $W(k)$, let $X_n \coloneqq X \times_{\Z} \Z/ p ^n \Z$ and consider the crystalline cohomology group $M^i \coloneqq \rH^i _\cris (X_n/W_n(k), \cO_{X_n/W_n (k)})$ with filtration $\rF^j (M^i) \coloneqq \rH^i _\cris ( X_n / W_n(k), J ^{[j]}_{X_n / W_n(k)})$, where $ J ^{[j]}_{X_n / W_n(k)}$ denotes the usual the $j$-th  divided powers ideal on the crystalline site. In \cite{FontaineMessing}, Fontaine--Messing show that there exist $\varphi_{W(k)}$-semi-linear maps $\varphi_j\colon \rH^i _\cris ( X_n / W_n(k), J ^{[j]}_{X_n/ W_n(k)}) \to  \rH^i _\cris (X_n /W_n(k), \cO_{X_n/W_n (k)}) $ for $i \leq p-1$. The main result of
\cite{FontaineMessing} can be summarized as follows.

\begin{theorem}\label{Thm-FM} With notations as above,
	\begin{enumerate}
		\item[(1)] For $i \leq p -1 $, $\left  (\rH^i _\cris (X_n /W_n(k), \cO_{X/W_n (k)}),\  \{ \rH^i _\cris ( X_n / W_n(k), J ^{[j]}_{X/ W_n(k)})\}  , \varphi_j \right )$
		is an object of $\MF^{[0, i ]}(W(k))$.
		\item[(2)] For $i \leq p-2$, let $M^i= (M^i, \rF^j (M^i), \varphi_j)$ denote the above object in $\MF^{[0, i]}(W(k))$. Then there exists a natural  $G_K$-isomorphism
		$\iota\colon  T_\cris (M^i)\simeq \rH^i _\et (X_{\Kbar}, \Z/p^n \Z). $
	\end{enumerate}
\end{theorem}
\begin{rem} Although $\rH^{p-1} _\cris (X/W_n(k), \cO_{X/W_n (k)})$ has the structure of a Fontaine-Laffaille module, its relation with $\rH^{p-1}_\et (X_{\Kbar}, \Z/ p^n \Z)$ is unknown.
\end{rem}

\subsection{The main results}

Our main goal is to extend Theorem \ref{Thm-FM} to the setting where $X$ is a smooth proper scheme over $R = W(k)[\![t_1, \dots, t_d]\!]$. We briefly explain our main results. 

We define the Frobenius $\varphi$ on $R$ by $\varphi (t_i)= (1+t_i)^p -1$ \footnote{Note that we do not use the more standard choice $\varphi(t_i) = t_i^p$ since our choice of $\varphi$ ensures that certain embeddings of $R$ into various period rings such as $A_{\cris}(R)$ are compatible with Frobenius. See \S \ref{subsec-period-ring} for details.}. Let $R_n = R/ p^n R$ and $\widehat{\Omega}_R \coloneqq \varprojlim_n \Omega_{R_n/W_n(k)}$ be the module of $p$-adically continuous K\"ahler differentials (where $\Omega_{R_n/W_n(k)}$ denotes the module of K\"ahler differentials), which is a finite free $R$-module of rank $d$ since $R/ pR$ has a finite $p$-basis 
 given by $\{t_1, \ldots, t_d\}$ (cf. \cite[Prop. 1.3.1]{BM}). Following \cite{FaltingsCcris}, we define a \emph{Fontaine--Laffaille module over $R_n$}: 

\begin{defn} \label{def-FL}A \emph{Fontaine--Laffaille module over $R_n$} is a tuple $(M, \rF^i (M), \varphi_i, \nabla)$ where
	\begin{enumerate}	
		\item[(1)] $M$ and $\rF^i(M)$ are finitely generated $R_n$-modules such that $\rF^i(M) = M$ for all $i \leq 0$ and $\rF^{r+1}(M) = 0$ for some $r$.
		\item[(2)] for each $i$, $\rF^{i+1}(M)$ is a direct summand of $\rF^i(M)$ as $R$-modules, and $\rF^i (M)\simeq \bigoplus _j R/ p^{a_j} R$ for some integers $a_j$'s.  
		\item[(3)] $\varphi_i\colon \rF^i(M)\otimes_{R, \varphi} R  \to M$ is an $R$-linear map such that 
		\begin{enumerate}
			\item $\varphi_i|_{\rF^{i+1}(M)\otimes_{R, \varphi} R} = p \varphi_{i+1}$;
			\item $\sum\limits_{i =0}^{r} \varphi_i  ( \rF^i (M)\otimes_{\varphi, R}  R)=M$ as $R$-modules; 
				\end{enumerate}
	\item[(4)] $\nabla: M \rightarrow M\otimes_R \widehat{\Omega}_R$ is an integrable connection such that
	\begin{enumerate}
		\item Griffiths-transversality holds: $\nabla(\rF^i(M)) \subset \rF^{i-1}(M)\otimes_R \widehat{\Omega}_R$,
		\item $\varphi_i$'s are parallel: $\nabla \circ (\varphi_i|_{\rF^i (M)}) = ((\varphi_{i-1}|_{\rF^{i-1}(M)})\otimes_R d\varphi_1)\circ \nabla$ as maps from $\rF^i(M)$ to $M \otimes_R \widehat{\Omega}_R$, where $d\varphi_1 \coloneqq \frac{d \varphi}{p}\colon \widehat{\Omega}_R \rightarrow \widehat{\Omega}_R$ (see the diagram before Definition \ref{defn:FL-mod}).
	\end{enumerate}
	\end{enumerate}
\end{defn}
   
Let $\overline{R}$ denote the union of normal $R$-sub-algebras $R'$ of a fixed separable closure of ${\mathrm{Frac}}(R[\frac{1}{p}])$ such that $R'[\frac{1}{p}]$ is finite \'{e}tale over $R[\frac{1}{p}]$. Set $G_R \coloneqq \R{Gal}(\overline{R}[\frac{1}{p}] / R[\frac{1}{p}]) = \pi_1^{\text{\'{e}t}}(\R{Spec}R[\frac{1}{p}], \eta)$ where $\eta$ is a fixed geometric point. In this setting, we can construct the integral crystalline period ring $A_\cris (R)$ over $R$, which is equipped with a filtration $\rF^i A_\cris(R)$, maps $\varphi_i\colon \rF^i A_\cris (R)\to A_\cris (R)$ and a continuous $G_R$-action. Unlike the classical situation (where $R= W(k)$), there is no longer a canonical embedding of $R$ into $A_\cris(R)$. However, we can choose an explicit embedding $R \to A_\cris (R)$ which is compatible with Frobenius $\varphi$ (see \S \ref{subsec-period-ring}). 

Given a Fontaine--Laffaille module $(M , \rF^i (M), \varphi _i, \nabla)$ as above, we can define a $\Z_p[G_R]$-module 
$$T_\cris (M) \coloneqq \left(\rF^r (A_\cris (R) \otimes_R M) \right)^{\varphi_r = 1} (-r)$$ as in the classical setting (for this definition to make sense, we need to have a structure theorem (Theorem \ref{thm-Faltings1}) for Fontaine--Laffaille modules over $R_n$). We note that the connection $\nabla$ is essential to construct the $G_R$-action on $T_\cris (M)$. We refer the reader to \S \ref{sec-relativeFL} for details.

If $X$ is a smooth proper scheme over $R$ with the standard divided power structure $(R, p R, \gamma_n (-)= \frac{ (-)^n}{n !})$ and $X_n \coloneqq X \times_{\Z} \Z/ p^n \Z $, then the torsion crystalline cohomology $\rH^i _\cris (X_n / R_n)$ comes equipped with a natural filtration given by the cohomology groups $\rH^i_\cris (X_n /R_n , J^{[j]}_{X_n/R_n} )$ and  $\varphi_{R}$-semilinear maps $\varphi_j\colon \rH^i _\cris (X_n / R_n, J^{[j]}_{X_n / R_n }) \to \rH^i_{\cris} (X_n / R_n) $. Furthermore,  $\rH^i _\cris (X_n / R_n)$ comes equipped with a natural connection (Gauss--Manin connection) $\nabla\colon \rH^i _\cris (X_n / R_n) \to \rH^i _\cris (X_n / R_n) \otimes_R  \wh \Omega_R$. Our main result is 

\begin{theorem}\label{thm-main} With notations as above,
	\begin{enumerate}
		\item[(1)] For $i \leq p -2$,  
		$M^i \coloneqq \left  (\rH^i _\cris (X_n /R_n),\  \{ \rH^i _\cris ( X_n / R_n, J ^{[j]}_{X_n / R_n})\}  , \varphi_j, \nabla \right )$
		is a Fontaine--Laffaille module over $R_n$.
		\item[(2)] For $i \leq p -2 $, there exists a natural $G_R$-isomorphism
		$T_\cris (M^i) \simeq \rH^i _\et (X_{\Rbar[\frac 1 p ]}, \Z/p^n \Z). $
	\end{enumerate}
\end{theorem}

\subsection{Historical remarks and our strategy}

Classical Fontaine--Messing theory has been generalized in various directions by several authors. In \cite{Breuil}, Breuil proves an analog of Theorem \ref{Thm-FM} in the setting of semi-stable schemes over $W(k)$. In \cite{CarusoInvent}, Caruso generalized Breuil's result to the setting of possibly ramified base rings. 

In \cite{FaltingsCcris}, Faltings generalized classical Fontaine--Messing theory to the relative setting. More precisely, Faltings considers a smooth scheme $Y$ over $W(k)$, and smooth proper schemes of relative dimension $g$ over $Y$. In this setting, Faltings develops a general theory of relative FL-modules on $X$ over $Y$, and allows more general coefficients (i.e. certain types of crystals $E$ on $X$ besides $\cO_{X/Y}$ as coefficients). In \cite[Thm. 6.2]{FaltingsCcris}, Faltings establishes a Fontaine--Messing type isomorphism comparing the crystalline pushforward of $E$ to $Y$ with an etale analog, under the assumption that $i +g \leq p -2$. If $f$ is projective, then this assumption can be weakened to $i \leq p-2$ (see Remark after Thm. 6.2 \emph{loc. cit.}). Moreover, Faltings proves analogous results in the semi-stable setting. The notion of a Fontaine--Laffaille module in Falting's relative setting was revisited by  D. Xu (\cite{DXU}), who obtained a topos theoretic framework for such objects. It is likely that an analogous framework also exists in our setting where the base is a power series ring, but we do not pursue this here.

\begin{rem}
In Theorem \ref{thm-main}, we obtain a period isomorphism without any restriction coming from the relative dimension i.e. we obtain an isomorphism for all $i \leq p-2$ in the possibly non-projective but proper case (unlike in the relative setting of Faltings). 
\end{rem}

The main new idea for proving Theorem \ref{thm-main} is the observation that one can ``descend" to a setting where the results of the classical Fontaine--Messing theory are applicable in order to obtain the results in our relative setting. More precisely, let $k_g$ be the perfection of ${\mathrm{Frac}}(R/pR)$ and let $R_g \coloneqq W(k_g)$. We have an embedding $b_g\colon R \inj R_g$ compatible with $\varphi$, and set $X_g \coloneqq X \times_R R_g$. On the other hand, consider the reduction modulo $t_i$ map $\bar{b}\colon R \to W(k)$ and let $X_{s} \coloneqq X \times_R  W(k)$. A key observation is that classical Fontaine--Messing theory applied to $X_g$ and $X_s$ implies the first part of Theorem \ref{thm-main}: if $\rH ^i _\cris (X_{s, n} / W_n(k))$ and $\rH ^i _\cris (X_{R_g,n} / R_{g,n})$ are both Fontaine--Laffaille modules with the same type, then $\rH^i_\cris (X_n / R_n)$ is a Fontaine--Laffaille module over $R_n$.  We use a similar strategy to establish an isomorphism between $T_\cris (M^i) $ and $\rH^i_\et (X_{\Rbar [\frac 1 p]},  \Z/ p^n \Z)$. Namely, we show that it suffices to establish a natural map $\iota\colon T_\cris (M^i) \to \rH^i_\et (X_{\Rbar [\frac 1 p]},  \Z/ p^n \Z) $ compatible with $G_R$-actions such that the base change of $ \iota$ to $R_g$ is the classical Fontaine--Messing comparison isomorphism $\iota_g\colon T_\cris (M_g^i) \simeq \rH^i_\et (X_{\Rbar_g [\frac 1 p]} , \Z/ p^n \Z)$.

To establish the comparison map $\iota$, we extend the method of Fontaine--Messing via syntomic cohomology to our base ring $R$. Classically, the key point is to show $\rH^i_{\syn} (X_{\cO_\Kbar}, S_n ^{[i]}) \simeq \rH ^i _\et (X_{\Kbar}, \Z/p ^n \Z (i))$, where $S_n^{[i]}$ are certain sheaves on the syntomic site defined via the divided power Frobenii. This is achieved by first constructing a morphism, and a computation of $p$-adic vanishing cycles yields an isomorphism. Our relative setting, while philosophically similar, is complicated by the fact that the ring $R_n$ is not perfect. In particular, basic computations dealing with the relevant crystalline cohomology groups are more involved, and the construction of $S_n^{[i]}$ becomes complicated. Our approach is to first descend to the intermediate ring $\wt R \coloneqq W(\widetilde{R}_1)$ where $\widetilde{R}_1$ denotes the perfection of $R_1$, and then compare the latter with $R$. This allows us to construct a comparison map as in the classical case. On the other hand, instead of having an analog of the classical computations of $p$-adic vanishing cycles, we show that our comparison map is an isomorphism using the map $b_g$.

We remark that passing to the classical case via $b_g$ can also be done over a more general base ring: if $S$ is a $p$-torsion free ring equipped with a lift of Frobenius such that $(p)$ is a prime ideal, then there is a $\varphi$-equivariant map $b_g\colon S \rightarrow W(k_g)$ lifting the natural map $S/pS \rightarrow k_g$, where $k_g$ is the perfection of $\mathrm{Frac}(S/pS)$. In the context of $p$-divisible groups, such reductions via $b_g$ are used in \cite{moon-BTdeform} for example. We expect that our method can also be used to prove results analogous to Theorem \ref{thm-main} in the general relative setting (not necessarily over power-series rings).

\subsection{Applications}

We explain some applications of Theorem \ref{thm-main}. A key benefit of the relative theory (in the current setting) is that it allows us to understand the full $G_R$-action on torsion \'etale cohomology. For example, we prove that $\rH^i_\et (X_{\Rbar [\frac 1 p]}, \bbQ_p)$ is crystalline as a $G_R$-representation for $i < p-1$ (cf. Theorem \ref{Thm-cryofHet}). We remark that $\rH^i_\et (X_{\Rbar [\frac 1 p]}, \bbQ_p)$ is expected to be crystalline for any $i$, but this needs to be studied via a different method. 

Secondly, we prove a type of local invariant cycle theorem. Suppose $k = \bar k$, and let $X_k \coloneqq X \times_R k$. When $i < p-1$, we show that the \'etale cohomology of the special fiber and that of the generic fiber are associated in the following way (cf. Proposition \ref{prop-etclosedfiber}): 
\[
\rH^i _\et (X_k, \Z/ p ^ n \Z) \simeq \rH^i_\et ( X_{\Rbar [\frac 1 p]}, \Z/ p ^n \Z )^{G_R}.
\]

For the last application, consider a finite totally ramified extension $K$ over $W(k)[\frac{1}{p}]$ with ramification index $e$, and let $\cO_K$ denote the ring of integers of $K$. Let $X$ be a proper smooth scheme over $\cO_K$, and $X_k \coloneqq X \times_{\cO_K} k$. When $e(i+1) < p-1$, Caruso proved that $\rH^i_\et (X_{\overline K}, \Z/p^n \Z)\otimes_{\Z_p} W(k)$ is isomorphic to $\rH^ i _\cris(X_k/W_n(k))$ as $W(k)$-modules (\cite{CarusoInvent}), and recently Li-Liu improved the bound to $ei < p-1$ (\cite{li-liu-prism-deRham-cohom}). On the other hand, \cite{SBM} gives some examples illustrating that $\rH^i_\et (X_{\overline K}, \Z/p^n \Z)\otimes_{\Z_p} W(k)$ is not isomorphic to $\rH^ i _\cris(X_k/W_n(k))$ in general when $ei \geq p-1$. It is natural to ask the following question.

\begin{ques}
When the ramification index $e$ is large and $i < p-1$, under which conditions is $\rH^i_\et (X_{\overline K}, \Z/p^n \Z)\otimes_{\Z_p} W(k)$ isomorphic to $\rH^ i _\cris(X_k/W_n(k))$?	
\end{ques}

\noindent Faltings showed that an isomorphism holds for any $e$ under the assumptions that $i+\dim{X_k} < p-1$ and the Hodge cohomology has no torsion (\cite{Faltings}). As a consequence of Theorem \ref{thm-main}, we prove:

\begin{theorem} (cf. Proposition \ref{prop-cohom-non-defect})
Suppose $X$	admits a lifting to a smooth proper scheme over $W(k)[\![t]\!]$. Then for any $i < p-1$, we have an isomorphism of $W(k)$-modules
\[
\rH^i_\et (X_{\overline K}, \Z/p^n \Z)\otimes_{\Z_p} W(k) \cong \rH^ i _\cris(X_k/W_n(k)).
\]
\end{theorem}

\subsection{Outline}
 
In \S2, we discuss (relative) period rings, Fontaine--Laffaille modules over $R$, and the functor $T_\cris$ to Galois representations. In  \S2.1, we recall the definition and some basic results on the period rings $A_{\cris}(R)$ and $\cO A_{\cris}(R)$. The latter ring is equipped with the natural structure of a relative Fonataine--Laffaile module over $R$, and allows one to construct the functor $T_{\cris}$ naturally. In \S2.2, following Faltings' approach in the relative case, we prove some basic structural results about relative FL-modules over $R$. One key result (Proposition \ref{prop-detect}) shows that $M$ is an FL-module over $R$ if and only if $M \otimes _R W(k)$ and $M \otimes_R R_g$ are FL-modules in the classical sense with the same type. In \S2.3, we discuss the functor $T_{\cris}$ and explain how to pass to torsion Galois representations. 

In  \S3, we show that the data arising from crystalline cohomology of a smooth proper scheme over $R$ gives rise to an FL-module. More precisely, we prove  Theorem \ref{thm-main} (1) in \S \ref{sec-cohoisFL}, by assuming  the  existence of $\varphi_i$ and $\nabla$ satisfying the required properties.
 
In \S 4, we construct a morphism $\iota\colon T_\cris (M^i) \to \rH^i_\et (X_{\Rbar [\frac 1 p]},  \Z/ p^n \Z) $ via syntomic cohomology following the strategy of Fontaine--Messing. Let $\Xbar \coloneqq X \times_R \Rbar$ and consider the presheaves on $(\overline X_n)_\syn$
$$\cO _n(U) =  \rH^0_\cris (U_n /R_n, \cO_{\Xbar_n / R_n}) \text{  and } \cJ^{[j ]}_{n}(U) \coloneqq \rH^0_\cris (U_n /R_n, J^{[j  ]}_{\Xbar_n / R_n}). $$
As in the classical setting, $\cO_n$ and $\cJ ^{[i]}_n$ are in fact sheaves and their (syntomic) cohomology computes $\rH_\cris^i (\Xbar_n/ R_n)$ and $\rH^i _\cris (\Xbar _n /R_n, J^{[j]}_{\Xbar/ R_n})$ respectively. We note that, since $R$ is not perfect, the computation of these cohomology groups is not as direct as in the classical setting. We consider analogous objects over $\wt{R} = W(\widetilde{R}_1)$ and then deduce the desired results over $R_n$. These computations are carried out in $\S 4.1-4.4$. In particular, we complete the proof that our objects are FL-modules by showing the existence of $\varphi_i$ and $\nabla$ satisfying the required properties.

The cohomology groups $\rH_\cris^i (\Xbar_n/ R_n)$ have the natural Gauss--Manin connection $\nabla$ over $R$.  To relate \'etale cohomology via syntomic cohomology, we consider the `nearby cycle' sheaves $S_n^{[i]} \coloneqq (\cJ_n ^{[i]})^{\varphi _i=1, \nabla= 0}$. This is another point of departure from the classical case, which makes it difficult to connect these sheaves directly to the classically defined sheaves $S_n^{[i]}$ when the base is $R_g$, where the latter sheaves are defined as $(\cJ^{[i]}_n)^{\varphi_i =1}$. Once again, to relate these two sheaves, we pass to the intermediate ring $\wt R$ and analogous sheaves defined in that setting; the latter sheaves can then be related directly to those over $R$ and those over $R_g$. In \S4.5, using the aforementioned ideas and constructions, we construct a natural morphism $T_{\cris}(M^i) \rightarrow  \rH^i_\syn (\Xbar_n,S_n ^{[i]})$. In \S4.6, following the idea of Fontaine-Messing, we construct a natural map $ \rH^i_\syn (\Xbar_n, S_n ^{[i]}) \to \rH^i_\et (X_{\Rbar [\frac 1 p]},  \Z/ p \Z) $. This completes the construction of a required map 
$\iota $ and proves Theorem \ref{thm-main}.  

Finally, in \S 5 we discuss the applications of our relative Fontaine--Messing theory.

\subsection{Notations} We reserve $\varphi$ to denote various Frobenii compatible with the arithmetic Frobenius on $W(k)$. If necessary, we add a subscript to $\varphi$ to indicate the corresponding ring or module. For example, we sometimes write $\varphi_M\colon M \to M$ in order to avoid confusion. For any module $M$ and integer $n \geq 1$, we set $M_n \coloneqq M \otimes_{\Z} \Z/ p ^n \Z$ and denote by $\wh M$ the $p$-adic completion of $M$.  For any scheme $X$, denote $X_n \coloneqq X \times_{\mathrm{Spec}\Z}\mathrm{Spec}(\Z/p^n\Z)$. Given a ring $A$ and $x \in A $, let $[x]$ denote its Teichm\"uller lift in $W(A)$. Let ${\mathrm{Cris}}(X/S)$ be the (small) crystalline site of $X$ over $S$ and $\cO_{X/S}$ denote the usual crystalline structure sheaf. We denote by $ \rH^i _\cris (X/S) \coloneqq \rH^i _\cris ((X/S)_\cris, \cO_{X/S})$ the usual crystalline cohomology of the structure sheaf. Let $J^{[j]}_{X/S}$ denote the $j$-th PD ideal sheaf over ${\mathrm{Cris}}(X/S)$. We reserve $\gamma_m (x)$ to denote the $m$-th divided power of $x$. Since most of our base schemes are affine schemes $\Spec(A)$ with a $\Z_p$-algebra $A$, we sometimes use $A$ to denote $\Spec(A)$ to simplify notation.

\section{Relative Fontaine--Laffaille modules and associated Galois representations}

In this section, we extend the theory of Fontaine--Laffaille modules and their associated Galois representations to the base ring $R = W(k)[\![t_1,\ldots,t_d]\!]$. The exposition follows closely the relative case studied in \cite{FaltingsCcris}, but we also give several new proofs/results using a d\'evissage to the ``generic fiber."

\subsection{Relative Period Rings}\label{subsec-period-ring}

In this subsection, following \cite{brinon-relative}, we briefly recall the construction of crystalline period rings in our relative setting. We refer to \textit{loc. cit.} for details. As in the introduction, let $k$ be a perfect field of characteristic $p > 2$, $R = W(k)[\![t_1, \ldots, t_d]\!]$, and define the Frobenius endomorphism $\varphi\colon R \rightarrow R$ extending the natural Frobenius on $W(k)$ by setting $\varphi(t_i) = (1+t_i)^p-1$. We remark that this is not the usual lift obtained by setting $\varphi (t_i) = t_i^p$, but we will see later that our choice is convenient for constructing a Frobenius equivariant embedding of $R$ into the crystalline period ring. Denote by $\Omega_{R_n/W_n(k)}$ the module of K\"ahler differentials, and let $\widehat{\Omega}_R \coloneqq \varprojlim_n \Omega_{R_n/W_n(k)}$. We have $\widehat{\Omega}_R \cong \bigoplus_{i=1}^d R\cdot dt_i$ by \cite[Prop. 1.3.1]{BM}.

\begin{rem} \label{rem:base-ring}
Let $W(k)\langle T_1^{\pm 1}, \ldots, T_d^{\pm 1} \rangle$ be the $p$-adic completion of $W(k)[T_1^{\pm 1}, \ldots, T_d^{\pm 1}]$, and consider the maximal ideal $\mathfrak{m} = (p, T_1-1, \ldots, T_d-1) \subset W(k)\langle T_1^{\pm 1}, \ldots, T_d^{\pm 1} \rangle$. The $\mathfrak{m}$-adic completion of the localization $W(k)\langle T_1^{\pm 1}, \ldots, T_d^{\pm 1} \rangle_{\mathfrak{m}}$ is isomorphic to $R = W(k)[\![t_1, \ldots, t_d]\!]$, where the isomorphism is given by $T_i \mapsto t_i+1$ for $i = 1, \ldots, d$. So $R$ is indeed an example of base rings studied in \cite{brinon-relative} (see \cite[Sec. 2]{brinon-relative} for the general conditions on base rings).   
\end{rem}

Let $\overline{R}$ denote the union of normal $R$-sub-algebras $R'$ of a fixed separable closure of $\mathrm{Frac}(R[\frac{1}{p}])$ such that $R'[\frac{1}{p}]$ is finite \'{e}tale over $R[\frac{1}{p}]$. It follows that $\Spec( \overline{R}[\frac{1}{p}])$ is a universal \'etale covering of $\R{Spec}(R[\frac{1}{p}])$, and $\overline{R}$ is the integral closure of $R$ in $\overline{R}[\frac{1}{p}]$. Let $G_R \coloneqq \pi_1^{\text{\'{e}t}}(\R{Spec}(R[\frac{1}{p}]), \eta)$ where $\eta$ is a fixed geometric point. We may identify the latter with the group $\Gal(\Rbar/R)$ of $R$-linear automorphisms.  Note that there is a natural embedding $\cO_{\overline{K}} \hookrightarrow \Rbar$ compatible with $W(k) \rightarrow R$.

Given a $p$-adically separated $\Z_p$-algebra $A$, set $\displaystyle{A}^{\flat} = \varprojlim_{x \mapsto x^p} A/pA$. By \cite[Lem. 3.2]{SBM}, the canonical map $\displaystyle \varprojlim_{x \mapsto x^p} \widehat{A} \rightarrow A^{\flat}$ is an isomorphism of multiplicative monoids. By functoriality of the construction, we have a natural embedding $\cO_\Kbar^\flat \hookrightarrow \Rbar^\flat$.

We have a natural surjective map $\theta\colon W(\overline{R}^{\flat}) \rightarrow \widehat{\overline{R}}$ which lifts the projection onto the first factor (and similarly for $\cO_{\Kbar}$). For integers $n \geq 0$, we choose $p_n \in \overline{R}$ such that $p_0 = p$ and $p_{n+1}^p = p_n$. Let $\tilde{p} \coloneqq (p_n)_{n \geq 0} \in \overline{R}^{\flat}$, and $[\tilde{p}] \in W(\overline{R}^{\flat})$ be its Teichm\"{u}ller lift. Then $\R{ker}(\theta)$ is generated by $\xi \coloneqq [\tilde{p}]-p$. Similarly, choose $\widetilde{1+t}_i = (z_{i, n})_{n \geq 0} \in \overline{R}^{\flat}$ where $z_{i, 0} = 1+t_i$ and $z_{i, n+1}^p = z_{i, n}$. Let $W(\overline{R}^{\flat})^{\R{DP}}$ (resp. $W(\cO_{\Kbar}^\flat)^{\R{DP}}$) be the divided power envelope of $W(\overline{R}^{\flat})$ (resp. $W(\cO_{\Kbar}^\flat)$) with respect to $\R{ker}(\theta)$, and let $A_{\R{cris}}(R)$ (resp. $A_\cris$) be its $p$-adic completion. $\varphi$ extends naturally to $A_\cris$ and $A_{\R{cris}}(R)$ compatibly with the embedding $A_\cris \hookrightarrow A_{\R{cris}}(R)$.

Let $\lambda_1\colon R \hookrightarrow W(\Rbar^\flat)$ be the $W(k)$-linear embedding given by $t_i \mapsto [\wt{1+t_i}]-1$ (for all $i$). To see that this indeed defines a ring morphism, it suffices to check by the universal property of $p$-adic Witt vectors that $\bar{\lambda}_1\coloneqq \lambda_1 \mod p$ defines a ring morphism $R/pR \rightarrow \overline{R}^{\flat}$, since the assignment $t_i \mapsto [\wt{1+t_i}]-1$ is compatible with the Frobenii (recall $\varphi_R(t_i) = (1+ t_i)^p -1$). Let $f(t_1, \ldots, t_d) \in R/pR = k[\![t_1, \ldots, t_d]\!]$. We need to show that $\bar{\lambda}_1(f) = f(\widetilde{1+t_1}-1, \ldots, \widetilde{1+t_i}-1)$ is a well-defined element in $\displaystyle \overline{R}^{\flat} = \varprojlim_{x \mapsto x^p} \overline{R}/p\overline{R}$. Note $(z_{i, n}-1)^{p^n} = z_{i, 0}-1 = t_i$. So $f(z_{1, n}-1, \ldots, z_{d, n}-1)$ is an element in $\overline{R}/p\overline{R}$ for each $n \geq 0$. Thus, $\bar{\lambda}_1(f) = (f(z_{1, n}-1, \ldots, z_{d, n}-1))_{n \geq 0} \in \overline{R}^{\flat}$, and $\lambda_1$ indeed gives a ring map which is $\varphi$-equivariant. Using $\lambda_1$, we regard $A_\cris(R)$ as an $R$-algebra compatible with the corresponding Frobenii.   

Given $i \geq 0$, let $\rF^iA_{\R{cris}}(R)\subset A_{\R{cris}}(R)$ be the $i$-th divided power ideal. By the description of $\R{ker}(\theta)$ above, this is given by $p$-adic completion of the ideal generated by the divided powers $\gamma_j(\xi)$ for $j \geq i$. This gives a filtration on $A_{\R{cris}}(R)$.
Note that $A_\cris (R)$ is $p$-torsion free by \cite[Prop. 6.1.3]{brinon-relative}. For each $i \leq p -1$, since $\varphi (\rF^i A_\cris(R))  \subset p ^ i A_\cris (R)$, we can define $\varphi_i \coloneqq \frac{\varphi}{p^i}\colon \rF^i A_\cris (R) \to A _\cris (R)$. The constructions of $\rF^i $ and $\varphi_i$ on $A_\cris (R)$ are compatible with those on $A_\cris$.
  
In the following, several of our arguments will involve d\'evissage to the perfect case via a base change map. We fix some notations and recall relevant results to facilitate this d\'evissage. Let $R_{(p)}$ denote the localization of $R$ at $(p)$ and denote by $R_g$ the $p$-adic completion of $\varinjlim_{\varphi} R_{(p)}$. Let 
\[
b_g\colon R \rightarrow R_g
\]
denote the natural morphism, which by construction is compatible with Frobenius. By \cite[Lem. 7.1.7]{brinon-relative}, the Frobenius endomorphism $\varphi\colon R \rightarrow R$ is flat, and thus $b_g\colon R \rightarrow R_g$ is flat by \cite[Tag 0912]{stacks-project}. Let $k_g \coloneqq \varinjlim_{\varphi} \R{Frac}(R/pR)$. Then, by the universal property of $p$-adic Witt vectors, we have a Frobenius equivariant isomorphism $R_g \cong  W(k_g)$. Let $\overline{R_g}$ be the ring of integers of an algebraic closure of $R_g[\frac{1}{p}] \cong W(k_g)[\frac{1}{p}]$.

We now fix an embedding $\overline{b_g}\colon \overline{R} \rightarrow \overline{R_g}$. Then $\overline{b_g}$ induces a map $ \Rbar^\flat \rightarrow \overline{R_g}^\flat$ and hence gives a map $W(\Rbar^\flat) \rightarrow W(\overline{R_g}^\flat)$. By construction, this induces a map $A_\cris (R) \rightarrow A_\cris (R_g)$ compatible with the filtration and $\varphi$. Consider the $R$-linear extension $\theta_R\colon R\otimes_{W(k)}W(\overline{R}^{\flat}) \rightarrow \widehat{\overline{R}}$ of $\theta$, and let $\mathcal{O}A_{\mathrm{cris}}(R)$ be the $p$-adic completion of the divided power envelope of $R\otimes_{W(k)}W(\overline{R}^{\flat})$ with respect to $\ker(\theta_R)$. The ring $\mathcal{O}A_{\mathrm{cris}}(R)$ is equipped with a natural $G_R$-action, and the Frobenius $\varphi_R\otimes\varphi_{W(\overline{R}^{\flat})}$ extends to $\mathcal{O}A_{\mathrm{cris}}(R)$. We equip $\mathcal{O}A_{\mathrm{cris}}(R)$ with a filtration by setting $\rF^i(\mathcal{O}A_{\mathrm{cris}}(R)) \subset \mathcal{O}A_{\mathrm{cris}}(R)$ for $i \geq 1$ to be the $p$-adic completion of the ideal generated by elements of the form $\gamma_{j_1}(a_1)\cdots\gamma_{j_l}(a_l)$ with $a_1, \ldots, a_l \in \ker(\theta_R)$ and $j_1+\cdots j_l \geq i$. Note that $A_{\mathrm{cris}}(R)$ and $\mathcal{O}A_{\mathrm{cris}}(R)$ have natural $G_R$-actions. The natural embedding $A_{\mathrm{cris}}(R) \hookrightarrow \mathcal{O}A_{\mathrm{cris}}(R)$ is compatible with the $G_R$-action, Frobenius, and filtration.

\begin{rem} \label{rem:notation-period-ring}
    Our notation for the period rings appearing above are different from \cite{brinon-relative}: $A_{\mathrm{cris}}(R)$ and $\mathcal{O}A_{\mathrm{cris}}(R)$ defined above correspond respectively to $A_{\mathrm{cris}}^{\nabla}$ in \cite[Def. 6.1.1]{brinon-relative} and $A_{\mathrm{cris}}$ in \cite[Def. 6.1.3]{brinon-relative}.
\end{rem}

Note that $\lambda_1$ induces an $R$-algebra structure on $\mathcal{O}A_{\mathrm{cris}}(R)$. On the other hand, setting $\lambda_2\colon R \rightarrow \mathcal{O}A_{\mathrm{cris}}(R)$ to be the map given by $\lambda_2(r) = r\otimes 1$ also gives an $R$-algebra structure, which is different from the morphism induced by $\lambda_1$ above. In the following, we will make use of both $R$-algebra structures in studying various Galois actions.

We end this section by recalling two results on the structure of the rings $A_{\mathrm{cris}}(R)$ and $\mathcal{O}A_{\mathrm{cris}}(R)$ proved in \cite{brinon-relative}. Denote by $A_{\mathrm{cris}}(R)\langle X_1, \ldots, X_d\rangle$ the $p$-adic completion of the divided power polynomial ring with variables $X_i$ having coefficients in $A_{\mathrm{cris}}(R)$. Consider the $A_{\mathrm{cris}}(R)$-linear map $f\colon A_{\mathrm{cris}}(R)\langle X_1, \ldots, X_d\rangle \rightarrow \mathcal{O}A_{\mathrm{cris}}(R)$ given by $X_i \mapsto (1+t_i)\otimes 1-1\otimes[\widetilde{1+t_i}]$. The following lemmas are proved in \cite{brinon-relative} (See Remark \ref{rem:base-ring} and \ref{rem:notation-period-ring}. For Lemma \ref{lem-2}, a more general statement will be proved in Lemma \ref{lem-toOAcris}.):
 
\begin{lemma} (cf. \cite[Cor. 6.1.2]{brinon-relative}) \label{lem-1} 
	For each integer $n \geq 1$, we have an isomorphism
	\[
	W_n(\Rbar^{\flat})[\delta_0, \delta_1, \ldots]/(p\delta_0-\xi^p, p\delta_{m+1}-\delta_m^p)_{m \geq 0} \stackrel{\cong}{\longrightarrow} A_{\R{cris}}(R)/p^n A_{\R{cris}}(R)
	\]
	where $\delta_m$ maps to $\gamma^{m+1}(\xi)$ (here, $\gamma: x \mapsto \frac{x^p}{p}$). 
\end{lemma}

Note that the isomorphism of the lemma is $W_n(\Rbar^{\flat})$-linear. In particular, we may define a $G_R$-action, filtration and divided power Frobenii $\varphi_i$ on the left hand side by transporting the corresponding structures on the right hand side via the isomorphism.
 
\begin{lemma} (cf. \cite[Prop. 6.1.5]{brinon-relative}) \label{lem-2}
The map $f\colon A_{\mathrm{cris}}(R)\langle X_1, \ldots, X_d\rangle \rightarrow \mathcal{O}A_{\mathrm{cris}}(R)$ defined above is an $A_{\mathrm{cris}}(R)$-linear isomorphism. In particular, it is an $R$-algebra morphism where the algebra structures are given via $\lambda_1$.
\end{lemma}

By the above lemma, we have a natural $A_{\mathrm{cris}}(R)$-linear connection $\nabla\colon \mathcal{O}A_{\mathrm{cris}}(R) \rightarrow \mathcal{O}A_{\mathrm{cris}}(R)\otimes_{R} \widehat{\Omega}_R$ given by $\nabla(t_i\otimes 1) = dt_i$.

\begin{rem} \label{rem:comp-structure}
Once again we may define a filtration, $G_R$-action, and $\varphi_i$ on the left hand side by pulling back the corresponding structures on the right hand side. We can give an explicit description as follows:
\begin{enumerate}
\item[(1)] Given $g \in G_R$, we have $g(X_i) = X_i - g([\widetilde{1+t_i}]) + [\widetilde{1+t_i}]$.
\item[(2)] The filtration on the left hand side is given by
$$
 \rF^{r }A_{\mathrm{cris}}(R)\langle X_1, \ldots, X_d\rangle  =  \text{the } p\text{-adic completion of}\bigoplus_{\sum i_j \geq r} \rF^{i_0 } A_{\cris}(R) \gamma_{i_1} (X_1)\cdots \gamma_{i_d} (X_d).
$$
\item[(3)] $\varphi(X_i) = (X_i +[\widetilde{1+t_i}])^p - [(\widetilde{1+t_i})^p]$. 
\item[(4)] Let $\pi\colon \mathcal{O}A_{\mathrm{cris}}(R) \rightarrow A_{\mathrm{cris}}(R)$ denote the composition of $f^{-1}$ and
 $q\colon A_{\mathrm{cris}}(R)\langle X_1, \ldots, X_d\rangle \rightarrow A_{\mathrm{cris}}(R)$ given by $X_i \mapsto 0$. This is an $R$-linear map with either $R$-algebra structure (the one induced by $\lambda_1$ or $\lambda_2$) on $\mathcal{O}A_{\mathrm{cris}}(R)$.
\item[(5)] Under the isomorphism $f$, $\nabla(\frac{\partial}{\partial t_i})(X_j)$ is $1$ if $i = j$ and $0$ if $i \neq j$.

\end{enumerate}
\end{rem}

\begin{lemma} \label{lem:graded-OAcris-p-tors-free}
$\mathrm{F}^r(\mathcal{O}A_{\mathrm{cris}}(R)) / \mathrm{F}^{r+1}(\mathcal{O}A_{\mathrm{cris}}(R))$ is $p$-torsion free for each $r \geq 0$.    
\end{lemma}

\begin{proof}
By Lemma \ref{lem-2} and Remark \ref{rem:comp-structure} (2), it suffices to show that $\mathrm{F}^r A_{\mathrm{cris}}(R) / \mathrm{F}^{r+1} A_{\mathrm{cris}}(R)$ is $p$-torsion free. We have $\mathrm{F}^0 A_{\mathrm{cris}}(R) / \mathrm{F}^{1} A_{\mathrm{cris}}(R) \cong \widehat{\overline{R}}$, which is $p$-torsion free. We claim that the map 
\[
h\colon \mathrm{F}^0 A_{\mathrm{cris}}(R) / \mathrm{F}^{1} A_{\mathrm{cris}}(R) \rightarrow \mathrm{F}^r A_{\mathrm{cris}}(R) / \mathrm{F}^{r+1} A_{\mathrm{cris}}(R)
\]
given by multiplication by $\frac{\xi^r}{r!}$ is an isomorphism. Note that any element in $\mathrm{F}^r A_{\mathrm{cris}}(R)$ can be written in the form $\sum_{i=r}^{\infty} a_i\frac{\xi^i}{i!}$ for some $a_i \in W(\overline{R}^{\flat})$ with $a_i \rightarrow 0$ $p$-adically as $i \rightarrow \infty$. So $h$ is surjective. 

Let $x \in \mathrm{F}^0 A_{\mathrm{cris}}(R) = A_{\mathrm{cris}}(R)$ such that $\frac{\xi^r}{r!}\cdot x \in \mathrm{F}^{r+1} A_{\mathrm{cris}}(R)$. Then
\[
\frac{\xi^r}{r!}\cdot x = \sum_{i=r+1}^{\infty} a_i\frac{\xi^i}{i!}
\]
for some $a_i \in W(\overline{R}^{\flat})$ with $a_i \rightarrow 0$ $p$-adically. Let $B_{\mathrm{dR}}^{\nabla +}(R)$ be the de Rham period ring constructed in \cite[Def. 5.1.1]{brinon-relative}, which is $\xi$-torsion free by \cite[Prop. 5.1.4 Pf.]{brinon-relative}. We have a natural map $A_{\mathrm{cris}}(R) \rightarrow B_{\mathrm{dR}}^{\nabla +}(R)$ which is compatible with the map $\theta$ and is injective by \cite[Prop. 6.2.1]{brinon-relative} (see Remark \ref{rem:notation-period-ring}). Thus, 
\[
x = \sum_{i=r+1}^{\infty} a_i\frac{r!}{i!}\xi^{i-r}
\]
as an element in $B_{\mathrm{dR}}^{\nabla +}(R)$, and $\theta(x) = 0$. So $x \in \mathrm{F}^{1} A_{\mathrm{cris}}(R)$, and $h$ is injective. This proves that $\mathrm{F}^r A_{\mathrm{cris}}(R) / \mathrm{F}^{r+1} A_{\mathrm{cris}}(R)$ is $p$-torsion free.  
\end{proof}

Note that any element in $\mathbb{Z}_p(1)$ is given by $\tilde{\epsilon} \coloneqq (\epsilon_n)_{n \geq 0}$ with $\epsilon_0 = 1$ and $\epsilon_{n+1}^p = \epsilon_n$. Considering $\tilde{\epsilon}$ as an element in $\overline{R}^{\flat}$, we have a canonical map $\beta\colon \Z_p(1) \rightarrow \rF^1 A_{\R{cris}}(R)$ given by $\beta(\tilde{\epsilon}) = \log ([\tilde{\epsilon}])$, which is $G_R$-equivariant.

\subsection{Relative Fontaine--Laffaille Modules}\label{subsec-rel-FM}

Following \cite{FaltingsCcris}, we define the categories of Fontaine--Laffaille modules over $R$ and study their basic properties. The arguments below are similar to those of \cite{FaltingsCcris}, with some minor modifications for our setting of power series base ring. 

Denote by $\R{MF}_{\R{big}}(R)$ the category whose objects consist of a $p$-power torsion $R$-module $M$, a sequence of $p$-power torsion $R$-modules $\rF^i(M)$, and sequences of $R$-linear maps $\iota_i\colon\rF^i(M) \rightarrow \rF^{i-1}(M)$, $\pi_i\colon \rF^i(M) \rightarrow M$, and $\varphi_i\colon \rF^i(M)\otimes_{R, \varphi} R \rightarrow M$, satisfying the following conditions:

\begin{itemize}
	\item[(1)] The composition $\rF^i(M) \overset{\iota_i}{\longrightarrow} \rF^{i-1}(M) \overset{\pi_{i-1}}{\longrightarrow} M$ is the map $\pi_i\colon \rF^{i}(M) \rightarrow M$;
	\item[(2)] The map $\pi_i\colon \rF^i(M) \rightarrow M$ is an isomorphism for $i \ll 0$;
	\item[(3)] The composition of $\varphi_{i-1}$ with $\iota_i\colon \rF^i(M) \rightarrow \rF^{i-1}(M)$ is $p\varphi_i$.
\end{itemize}
Morphisms between such objects are compatible collections of $R$-linear maps between $M$'s and $\rF^i(M)$'s. Note that $\R{MF}_{\R{big}}(R)$ is an abelian category. 

\begin{rem}
An $R$-linear map $\varphi_i\colon \rF^i(M)\otimes_{R, \varphi} R \rightarrow M $ is equivalent to a $\varphi_R$-\emph{semilinear} map $\varphi_i\colon \rF^i (M)\to M$. In the following, we will use these two notions of $\varphi_i$ interchangeably.
\end{rem}

\begin{rem}
	Let $\wt{M}$ be the colimit of the following diagram:
	\[
	\cdots \xrightarrow{\iota_{i+2}} \rF^{i+1}(M) \xleftarrow{\times p} \rF^{i+1}(M) \xrightarrow{\iota_{i+1}} \rF^i(M) \xleftarrow{\times p} \rF^i(M) \xrightarrow{\iota_i} \rF^{i-1}(M) \xleftarrow{\times p} \cdots.
	\]
	Then condition (3) above is equivalent to requiring that $\varphi_i$'s induce an $R$-linear map $\varphi_M\colon \wt{M}\otimes_{R, \varphi} R \rightarrow M$.
\end{rem}

Let $\R{MF}(R)$ be the full subcategory of $\R{MF}_{\R{big}}(R)$ whose objects consist of tuples $(M, \rF^i(M), \varphi_i)$ such that $M$ and $\rF^i(M)$ are finitely generated $R$-modules, $\rF^i(M) = \{0\}$ for $i \gg 0$, and $\varphi_i$'s induce an isomorphism $\varphi_M\colon \wt{M}\otimes_{R, \varphi} R \simeq M$. We further denote by $\R{MF}^{[a,b]}(R)$ the full subcategory of $\R{MF}(R)$ consisting of objects such that $\rF^{i}(M)= \{0\}$ for $i \geq b+1$ and $\pi_i$ for $i \leq a$ is an isomorphism. We refer to objects in $\MF(R)$ as \emph{Fontaine--Laffaille modules} over $R$ without connection $\nabla$. We often drop ``without connection $\nabla$" or ``over $R$'' if no confusion arises.

Let $\R{MF}_{\R{big}}^{\mathrm{fin}}(R)$ denote the full sub-category of $\R{MF}_{\R{big}}(R)$ consisting of objects such that there exist integers $a \leq b$ so that $\rF^{i}(M)= \{0\}$ for $i \geq b+1$ and $\pi_i$ for $i \leq a$ is an isomorphism. For any such object $M$, $\wt M$ can be defined explicitly by the right exact sequence
\begin{equation}\label{eqn-tildeM}
\xymatrix{ \bigoplus_{i = a+1} ^b \rF^i (M) \ar[r]^-{\theta_M} &  \bigoplus_{i = a} ^b \rF^i (M)  \ar[r] & \wt M \ar[r] & 0,}
\end{equation}
where $\theta_M (( x_{a+1}, \dots, x_b)) = (\iota_{a+1} (x_{a+1}), -p x_{a+1} + \iota_{a+2} (x_{a+2}), \dots, - px_{b-1} + \iota_b (x_b)  , -p x_b)$. It follows that the functor $M \mapsto \wt  M$ is a right exact functor on $\R{MF}_{\R{big}}^{\mathrm{fin}}(R)$. Moreover, if $A$ is a $W(k)$-algebra equipped with a Frobenius endomorphism compatible with that on $W(k)$ and $b\colon R \to A$ is a $W(k)$-algebra morphism compatible with Frobenius, then $\wt{M\otimes_{R, b} A} \simeq \wt M \otimes_{R, b} A$ (with $\rF^i (M\otimes_{R, b} A) \coloneqq\rF^i(M)\otimes_{R, b} A$). Hence if $M$ is an object in $\MF(R)$, so is $M/ p^j M$, and we have $M_0 \coloneqq M\otimes_{R, \bar{b}} W(k) \in \MF(W(k))$ for $\bar{b}\colon R \rightarrow W(k)$ given by $t_i \mapsto 0$ and $M_g \coloneqq M\otimes_{R, b_g} R_g \in \MF(R_g)$ for $b_g\colon R \rightarrow R_g$ as in previous section.

We have the following structural result for Fontaine-Laffaille modules over $R$. As in the classical case, one consequence is that $\R{MF}(R)$ is an abelian category.

\begin{theorem}\label{thm-Faltings1}(cf. \cite[Thm. 2.1]{FaltingsCcris}) Let $(M, \rF^i(M), \varphi_i) \in \R{MF}(R)$. Then
	\begin{enumerate}
		\item[(1)]  All maps $\rF^{i+1}(M) \rightarrow \rF^i(M) \rightarrow M$ are injections onto direct summands as $R$-modules. In particular, we can consider $\rF^i(M)$'s as submodules of $M$ giving a finite filtration.
		
		\item[(2)] $\rF^i(M)$ and $\rF^i (M)/ \rF^{i+1}(M)$ are isomorphic to a direct sum of $R$-modules of the form $R/p^e R$.
		
		\item[(3)] For $M, N \in \R{MF}(R)$, any map $f\colon M \rightarrow N$ in $\R{MF}(R)$ is strict with respect to the filtration. Namely, $f (\rF^i(M)) = f(M)\cap \rF^i(N). $
	\end{enumerate}
\end{theorem}

\begin{proof} The proof of \cite[Thm. 2.1]{FaltingsCcris} applies almost verbatim for our base $R = W(k)[\![t_1, \ldots, t_d]\!]$. For the convenience of readers, we sketch the proof below. 
	
Note first that the corresponding statements for objects in $\mathrm{MF}(R_g) = \mathrm{MF}(W(k_g))$ hold by the classical Fontaine-Laffaille theory as follows. In that case, (2) and (3) follow directly from \cite[Sec. 1]{Fonatine-Laffaille}. For (1), the length of $\wt{M} \geq$ the length of $M$ for any $M \in \R{MF}(R_g)$, with equality if and only if all maps $\rF^{i+1}(M) \rightarrow \rF^i(M)$ are injective. This proves the first part of (1). The assertion about direct summands in (1) follows by applying this fact to $M/p^j M$'s (which are objects in $\R{MF}(R_g)$).
	
For our base $R$, assume first that $p\mathrm{F}^j(M) = 0$ for all $j$. For any finite $R/pR$-module $L$, denote by $\R{Fitt}_i(L)$  the $i$-th Fitting ideal of $L$ over $R/pR$. Note $M \in \R{MF}^{[a,b]}(R)$ for some integers $a \leq b$ by the definition. So $\wt{M}$ is the cokernel of the map
	\[
	\bigoplus_{a+1 \leq  i \leq b} \rF^i(M) \stackrel{\theta_M}{\longrightarrow} \bigoplus_{a \leq  i \leq b} \rF^i(M).
	\]
	 Let $r_i$ be the rank of $\rF^i(M)$ at the generic point of $R/pR$. Then $r_i$ is the smallest subscript $s$ for which $\R{Fitt}_s(\rF^i(M)) \neq\{ 0\}$. If we let $r= r_a$ be the rank of $M $ and $\wt{M}$ at the generic point, then
	\[
	\prod_{a  \leq  i \leq b } \R{Fitt}_{r_i}(\rF^i(M)) \supset \R{Fitt}_r(\wt{M})\cdot \prod_{a+1\leq  i \leq b } \R{Fitt}_{r_i}(\rF^i(M)).
	\]
	On the other hand, since $\varphi_M$ is an isomorphism, we have
	\[
	\R{Fitt}_{r_a}(\rF^a(M)) = \R{Fitt}_r(M) = \varphi_*(\R{Fitt}_r(\wt{M})) = \mbox{the ideal generated by } \varphi(\R{Fitt}_r(\wt M)).
	\]
	Since $\R{gr}_{\fm} R$ has no zero divisors where $\mathfrak{m} \subset R$ denotes the maximal ideal, this implies that $\R{Fitt}_r(M) = R/pR$, and therefore $M$ and $\wt{M}$ are free over $R/pR$. Since $p\mathrm{F}^j(M) = 0$ for all $j$, $\wt{M}$ is the direct sum of the cokernels of the maps $\rF^{i+1}(M) \rightarrow \rF^i(M)$ from the definition of $\theta_M$. So these cokernels are free over $R/pR$, and the images are direct summands. Furthermore, from the case of $\mathrm{MF}(R_g)$ above, the maps $\mathrm{F}^{i+1}(M)\otimes_R R_g \rightarrow \mathrm{F}^{i}(M)\otimes_R R_g$ are injective. Since $R_{(p)} \rightarrow R_g$ is faithfully flat, we deduce (1) and (2) by decreasing induction on $i$ in this case when $p$ annihilates all modules.
	
For the general case, we induct on the smallest integer $e$ such that $p^e$ annihilates all modules involved. Let $M \in \R{MF}(R)$. Since $\wt{M}$ defines a right exact functor of $M$, the assertions (1) and (2) hold for $M/pM$. Thus, each map $\rF^i(M)/p\rF^i(M) \rightarrow \rF^{i-1}(M)/p\rF^{i-1}(M)$ is injective. Let $N \subset M$ be the object $pM$ given by $\mathrm{F}^i(N) = p\mathrm{F}^i(M) \subset \mathrm{F}^i(M)$ with the induced maps $\iota_i, \pi_i, \varphi_i$. Since $\rF^i(M)/p\rF^i(M) \rightarrow \rF^{i-1}(M)/p\rF^{i-1}(M)$ is injective, the induced map $\wt{N} \rightarrow \wt{M}$ is injective by the snake lemma, identifying $\wt{N}$ with $p\wt{M} \subset \wt{M}$. Since $\varphi\colon R \rightarrow R$ is flat, we deduce that $\varphi_N$ is an isomorphism and $N \in \R{MF}(R)$. 

By the inductive hypothesis, (1) and (2) hold for $N$. So $\rF^{i-1}(N) \cong \rF^i (N) \oplus \rF^{i-1}(N) /\rF^i (N)$ and both $\rF^{i-1}(N)$ and $\rF^i(N)$ are direct sums of the form $R/p^a R$. That is, $\rF^{i-1}(N) \cong \bigoplus_{j=1}^n R/p^{a_j}R$ for some integers $a_j \geq 1$ and $\rF^i(N)\cong \bigoplus_{j=1}^m R/p^{a_j}R$ for $m \leq n$. Set $Q = \bigoplus_{j=1}^n R/p^{a_j+1}R \cdot e_j$ and
	$L =\bigoplus_{j=1}^m R/p^{a_j+1}R \cdot e_j$ with $e_j$ being basis elements. We have an isomorphism $h'\colon \rF^{i-1}(N) \stackrel{\cong}{\longrightarrow} pQ$ which induces an isomorphism $\rF^i(N) \cong pL$. Since $\rF^i(M)/p\rF^i(M) \rightarrow \rF^{i-1}(M)/p\rF^{i-1}(M)$ is an injection of free $R/pR$-modules, we can lift $h'$ to an $R$-module morphism $h\colon \rF^{i-1}(M) \rightarrow Q$ such that $h \circ \iota_i(\rF^i(M))\subset L$. For the same reason, ${h'}^{-1}\colon pQ \rightarrow \rF^{i-1}(N)$ and ${h'}^{-1}\colon pL \rightarrow \rF^{i}(N)$ lift to $R$-module maps $f\colon Q \rightarrow \rF^{i-1}(M)$ and $f\colon L \rightarrow \mathrm{F}^i(M)$ compatibly with $\iota_i\colon \mathrm{F}^i(M) \rightarrow \mathrm{F}^{i-1}(M)$. The composite $h \circ f\colon Q \rightarrow Q$ is a lift of $h' \circ {h'}^{-1} = \mathrm{Id}\colon pQ \rightarrow pQ$, so $(h\circ f)(e_j) = e_j+x_j$ for some $p$-torsion elements $x_j \in Q$. Since $Q \cong \bigoplus_{j=1}^n R/p^{a_j+1}R\cdot e_j$, we deduce that $h\circ f$ is an isomorphism. This implies that $Q$ injects onto a direct summand of $\rF^{i-1}(M)$, which we identify with $Q$. Similarly, $L$ is a direct summand of $\rF^i (M)$. Since $p\rF^{i-1}(M) = \rF^{i-1}(N) = pQ$, we have $\rF^{i-1}(M)/p\rF^{i-1}(M) \cong Q/pQ \oplus \rF^{i-1}(M)/Q$, and similarly $\rF^i(M)/ p \rF^i(M) \cong L/ pL \oplus \rF^i(M)/ L$. Since the map $\rF^i(M)/p\rF^i(M) \rightarrow \rF^{i-1}(M)/p\rF^{i-1}(M)$ identifies $\rF^i(M)/p\rF^i(M)$ as a direct summand of $\rF^{i-1}(M)/p\rF^{i-1}(M)$, the induced map $\mathrm{F}^i(M)/L \rightarrow \mathrm{F}^{i-1}(M)/Q$ identifies $\mathrm{F}^i(M)/L$ as a direct summand of $\mathrm{F}^{i-1}(M)/Q$. Thus, $\rF^i(M) \cong L \oplus \rF^i(M)/L$ is a direct summand of $\rF^{i-1}(M) \cong Q \oplus \rF^{i-1}(M) / Q $. We conclude that $\rF^{i-1}(M) \cong \rF^i(M) \oplus \rF^{i-1}(M)/ \rF^i(M)$, and both $\rF^i(M)$ and $\rF^{i-1}(M)/ \rF^i(M)$ are direct sums of the form $R/ p^a R$.
		
Finally, we prove assertion (3). We first show that $L \coloneqq N/f(M)$ (as an object in $\R{MF}_{\R{big}}(R)$) lies in $\MF(R)$. Note that $L \in \R{MF}_{\R{big}}^{\mathrm{fin}}(R)$. Since the functor $(-) \mapsto \widetilde{(-)}$ is right exact on $\R{MF}_{\R{big}}^{\mathrm{fin}}(R)$, $\wt M \to \wt N \to \wt L \to 0$ is exact. Since $\varphi_M\colon \wt M\otimes_{R, \varphi} R \to M$ and $\varphi_N\colon \wt N\otimes_{R, \varphi} R \to N$ are isomorphisms, the induced map $\varphi_L\colon \wt L\otimes_{R, \varphi} R \to L$ is an isomorphism. So $L = N/f(M)$ is an object in $\MF(R)$. In particular, $\rF^i(L) = \rF^i(N) / f(\rF^i(M))$ is a submodule of $L = N /f(M)$. This forces $f(M) \cap \rF^i (N) = f(\rF^i (M))$.
\end{proof}

\begin{rem}\label{rem-adpated} Any $M\in \MF(R)$ admits an \emph{adapted basis}: by above theorem (1) and (2), we can inductively choose (by decreasing induction on $i$ for $\mathrm{F}^i(M)$) elements $e_j \in M$ with positive integers $a_j$ and decreasing positive integers $\cdots \geq m _i  \geq m _{i+1} \geq \cdots$ so that if $\rF^i (M) \not = 0$,
	then $$\rF^i(M) = \bigoplus\limits_{j =1}^{m _i} R/p^{a_j} R \cdot e_j.$$
\end{rem}

\begin{corollary} \label{cor:MF(R)-abelian-category}
$\R{MF}(R)$ is an abelian subcategory of $\R{MF}_{\R{big}}(R)$.	
\end{corollary}

\begin{proof}
Let $f\colon M \rightarrow N$ be a morphism of objects in $\mathrm{MF}(R)$, and let $P = \ker(f)$ and $L = \mathrm{coker}(f)$ as objects in $\R{MF}_{\R{big}}(R)$. We showed $L \in \mathrm{MF}(R)$ in above proof of Theorem \ref{thm-Faltings1}, so it suffices to show $P \in \mathrm{MF}(R)$. We first claim that $f(M) \in \mathrm{MF}(R)$. Note that $f(M), P \in \R{MF}_{\R{big}}^{\mathrm{fin}}(R)$. The exact sequence 
\[
0 \rightarrow f(M) \rightarrow N \rightarrow L \rightarrow 0
\]
induces an exact sequence
\[
\widetilde{f(M)} \rightarrow \widetilde{N} \rightarrow \widetilde{L} \rightarrow 0.
\]
By the above theorem, the maps $\iota_i\colon \rF^{i}(N) \rightarrow \mathrm{F}^{i-1}(N)$ are injective and $\rF^i (f(M)) = \rF^i (N) \cap f(M)$. So from the presentation (\ref{eqn-tildeM}) of $\wt M$ using $\theta_M$, we deduce that $\wt{f(M)} \rightarrow \wt N$ is injective. Thus, the above sequence is also left exact, which implies $\varphi_{f(M)}\colon \wt{f(M)}\otimes_{R, \varphi} R \rightarrow f(M)$ is an isomorphism and $f(M) \in \mathrm{MF}(R)$.

Now, applying a similar argument for the exact sequence $0 \rightarrow P \rightarrow M \rightarrow f(M) \rightarrow 0$, we deduce $P \in \mathrm{MF}(R)$.
\end{proof}

In the following, we give another criterion for an object of $\MF_\R{big}(R)$ to be in $\MF(R)$. We say a $W(k)$-module $M$ and an $R_g$-module $M'$ have the \emph{same type} if $M' \simeq M \otimes_{W(k)} R_g $ as $R_g$-modules. Note that the projection $\bar{b}\colon R \rightarrow W(k)$ given by $t_i \mapsto 0$ is compatible with Frobenius.

\begin{proposition} \label{prop-detect} Suppose $M\in \MF_\R{big} (R)$ such that each $\rF^i(M)$ is finite as an $R$-module. Then $M$ is in $\MF(R)$ if and only if
	\begin{enumerate}\item[(1)] Both $ M_0 \coloneqq M \otimes_{R, \bar{b}} W(k)$ and $M_g \coloneqq M \otimes_R R_g$ are objects in $\MF(W(k))$ and $\MF(R_g)$ respectively;
		\item[(2)] $M_0$ and $M_g$ have the same type.
	\end{enumerate}
\end{proposition}

\begin{proof} 
If $M$ is an object in $\MF(R)$ then we have $M_0 \in \MF(W(k))$, $M_g \in \MF(R_g)$, and by Theorem \ref{thm-Faltings1} (2), $M_0$ and $M_g$ have the same type. For the converse statement, we first reduce to the case where $pM = \{0\}$. Consider the exact sequence $$0 \to pM \to M \to M/pM \to 0.$$ Since $M/pM \otimes_R W(k) \simeq M_0/ p M _0$ and $M/pM \otimes_R R_g \simeq  M_g/p M_g$, $M/pM$ satisfies the conditions in the proposition. Suppose $M/pM$ is in $\MF(R)$. Then all $\rF^i(M/pM)$ are finite $R/pR$-free, so the sequence
	$$ 0 \to p\rF^i(M) \otimes_R W(k) \to \rF^i (M_0) \to \rF^i (M/pM)  \otimes_R W(k) \to 0$$ is exact.
	Since $\MF(W(k))$ is an abelian category, $pM \otimes_R W(k)$ is an object in $\MF (W(k))$. Similarly, $pM \otimes_R R_g$ lies in $\MF(R_g)$. Using the above exact sequence for $i \ll 0$, we deduce that $pM \otimes_R W(k) $ and $pM \otimes_R R_g$ have the same type. By induction on the length of $M\otimes_R W(k)$ as a $W(k)$-module, we may assume $pM$ is an object in $\MF(R)$. Note that if $N$ is a finite $R$-module such that $N\otimes_R W(k) = 0$, then $N = 0$ by Nakayama's Lemma. So $M, pM \in \MF_\R{big}^{\mathrm{fin}}(R)$.
	By the commutative diagram
	\begin{equation}\label{diagram-1}
	\xymatrix{ & \widetilde {pM}\otimes_{R, \varphi} R \ar[d]^\wr \ar[r] & \wt  M\otimes_{R, \varphi} R \ar[d] \ar[r] & \widetilde {M/pM}\otimes_{R, \varphi} R \ar[r]\ar[d] ^\wr & 0 \\ 0 \ar[r] & pM \ar[r] & M \ar[r] &  {M/pM} \ar[r] & 0}
	\end{equation}
	whose rows are exact, we conclude that $\varphi_M\colon \wt M\otimes_{R, \varphi} R \to M$ is an isomorphism and $M \in \MF(R)$.
	
	Now it suffices to prove the proposition when $M$ is killed by $p$. Since $M_0 \otimes_k k_g \simeq M_g$, $M$ is finite free over $R/pR$ by Nakayama's Lemma. Pick an adapted basis $\bar e_1 , \dots , \bar e_n$ of $M_0$ so that $\bar e_1, \dots, \bar e_{m_i}$ (with $m_i \leq n$) is a basis of $\rF^i (M_0)$. Let $e_1 , \dots ,e_{m_i} \in \rF^i (M)$ be a choice of lifts of $\bar e_1 , \dots , \bar e_{m_i}$ and $e_{m_i+1}, \dots , e_n \in M$ be a choice of lifts of $\bar e_{m_i+1}, \dots, \bar e_n$. For the structure map $\iota_i\colon \rF^i(M) \to M$, we deduce from Nakayama's Lemma that $ e_1 , \dots,  e_{m_i}$ generate $\rF^i(M)$ and $\iota_i (e_1), \dots \iota_i(e_{m_i}), e_{m_i+1}, \dots, e_n$ form a basis of $M$. This forces $\iota_i$ to be injective and $e_1 , \dots ,e_{m_i}$ to be linearly independent. Thus, $\rF^i(M)$ is a finite free $R/pR$-module and $\rF^i(M)$ is a direct summand of $M$ via $\iota_i$. Moreover, $\wt M$ is a finite free $R/pR$-module with the same rank as that of $M$. Hence, $\varphi_M\colon \wt M\otimes_{R, \varphi} R \to M$ is an isomorphism since $\varphi_M \otimes_R W(k)$ is an isomorphism.
\end{proof}

\begin{rem} \label{rem-filrange} If $M \in \MF(R)$, then  $M \in \MF^{[a, b]} (R)$ if and only if $M_0 \in \MF ^{[a, b]} (W(k))$ if and only if $M_g \in \MF ^{[a, b]} (R_g)$.  
\end{rem}

Let $\mathrm{MF}_{\nabla, {\mathrm big}}(R)$ be the category whose objects are tuples $(M, \rF^i(M), \varphi_i)$ in $\mathrm{MF}_{\mathrm big}(R)$ equipped with an integrable connection $\nabla_i \colon \rF^i  M \rightarrow \rF^{i-1} M\otimes_R \widehat{\Omega}_R$ such that 
\begin{itemize}
	\item Griffiths-transversality holds: the following diagram commutes
 \[ \xymatrixcolsep{5pc}\xymatrix{ \rF ^i M \ar[r] ^{\iota_i} \ar[d]^{\nabla_i} & \rF^{i -1} M \ar[d]^{\nabla_{i -1}} \\ \rF^{i -1} M \otimes _R \wh \Omega_R  \ar[r]^{\iota_{i -1} \otimes_R \wh \Omega_R } & \rF^{i-2} M  \otimes_R \wh \Omega_R .}\]
	\item the semi-linear  maps $\varphi_i\colon \rF^i(M) \rightarrow M$ are parallel: $\nabla  \circ \varphi_i = (\varphi_{i-1}\otimes_R d\varphi_1)\circ \nabla_i$ as maps from $\rF^i(M)$ to $M \otimes_R \widehat{\Omega}_R$, where $d\varphi_1 \coloneqq \frac{d \varphi}{p}\colon \widehat{\Omega}_R\rightarrow \widehat{\Omega}_R$ is the $\varphi_R$-semilinear map given by $dt_i \mapsto \frac{d\varphi(t_i)}{p}$, $\nabla \colon M \to M \otimes_R \wh \Omega_R$ is given by $\nabla_j$ for $j \ll 0$ via $\pi_j\colon \rF^j(M) \stackrel{\cong}{\rightarrow} M$, and we use the diagram 
  \[ \xymatrixcolsep{5pc}\xymatrix{ \rF^i(M)  \ar[d] ^{\nabla_i }\ar[r]^-{\varphi_i}&  M \ar[d]^- \nabla   \\   \rF^{i-1}(M)  \otimes _R \wh \Omega_R     \ar[r]^-{\varphi_{i -1} \otimes d \varphi _1}   &  M \otimes_R \wh \Omega_R .& } \]
\end{itemize}
In the situation that each $\iota_i$ is injective (for example, when $M \in \MF(R)$), note that $\nabla_i = \nabla|_{\rF ^i (M)}$ and above Griffiths-transversality is the usual definition: $\nabla(\rF^i(M)) \subset \rF^{i-1}(M)\otimes_R \widehat{\Omega}_R$. A morphism between two objects in $\MF_{\nabla, {\mathrm big}}$ is an $R$-linear map compatible with all structures $\rF^i, \varphi_i$ and $\nabla_i$.  

Define $\R{MF}_{\nabla}(R)$ to be the full subcategory of $\MF_{\nabla, {\mathrm big}}$ whose $(M, \rF^i(M), \varphi_i)$ is in $\R{MF}(R)$. Denote by $\R{MF}^{[0, r]}_\nabla(R)$ the full subcategory of $\R{MF}_{\nabla}(R)$ with $M \in \mathrm{MF}^{[0, r]}(R)$. In most cases in this paper, $r \leq p-2$. 

\begin{defn} \label{defn:FL-mod} A \emph{Fontaine--Laffaille module over $R$} is an object  $(M , \rF^i (M), \varphi _i, \nabla)\in \MF_\nabla^{[0, r]}(R)$ for some $r \geq 0$.
\end{defn}

\begin{rem}
This definition is different from Definition \ref{def-FL}, where we use a weaker condition $\sum_i \varphi _i(\rF^i (M)) = M$ instead of $\varphi_M\colon \wt M \otimes_{R, \varphi} R \simeq M$. But the weaker condition still implies that $M_0$ and $M_g$ are objects in $\MF(W(k))$ and $\MF(R_g)$ respectively, so it follows from Proposition \ref{prop-detect} that the two definitions are equivalent.  
\end{rem}

\subsection{Torsion Representations associated to Fontaine--Laffaille Modules}\label{sec-relativeFL}

In this section, we study Galois representations associated to Fontaine--Laffaille modules $M\in \MF^{[0, r]}_\nabla (R)$ with $r \leq p-2$. Recall that we defined embeddings $\lambda_1\colon R \rightarrow A_{\cris}(R)$ and $\lambda_2\colon R \rightarrow \cO A_{\cris}(R)$ in \S \ref{subsec-period-ring}. Let $\overline M_R \coloneqq \cO A_{\cris} (R) \otimes _{\lambda_2, R} M$. We claim that the natural map $\rF^i \cO A _\cris (R) \otimes_{\lambda_2, R} \rF ^{r-i}(M) \to  \Mbar _R$ is injective so that one can define $\rF^r \Mbar_R \coloneqq \sum\limits_{i = 0}^r \rF^i \cO A _\cris (R) \otimes_{\lambda_2, R} \rF ^{r-i} (M) \subset \Mbar _R$. In order to see this, first recall that $\rF^{r-i}(M)$ is a direct summand of $M$, so that $\cO A _\cris (R) \otimes_{\lambda_2, R} \rF ^{r-i}(M)$ injects into $\Mbar_R$. Since $\rF^{r-i}(M)$ is a direct sum of modules of the form $R/p^n R$ and $\rF^i \cO A _\cris (R)/ p ^n \to \cO A _\cris (R)/ p ^n$ is injective by Lemma \ref{lem:graded-OAcris-p-tors-free}, the claim follows. 

Consider the map $$\varphi_r \coloneqq \sum_{i=0}^r \varphi_i|_{\rF^i \cO A _\cris (R)} \otimes_{\lambda_2, R} \varphi_{r-i} |_{\rF^{r-i}(M)}\colon \rF^r \Mbar_R \to \Mbar _R. $$ To see this is well-defined, we check that $\varphi _r$ agrees on the intersection of $\rF^i \cO A _\cris (R) \otimes_{\lambda_2, R} \rF ^{r-i}(M)$  and $\rF^{j} \cO A _\cris (R) \otimes_{\lambda_2, R} \rF ^{r-j}(M)$ for $0 \leq i < j \leq r$. Since $M$ has an adapted basis by Remark \ref{rem-adpated}, we may assume an element in such an intersection has the form $a \otimes x$ with $x \in \rF^{r-i}(M)$ and $a \in \rF ^{j} \cO A _\cris (R)$. We have $\varphi_i (a) \otimes \varphi_{r-i} (x) =  p^{j-i}\varphi_{j} (a) \otimes \varphi_{r-i} (x)  = \varphi_{j} (a) \otimes \varphi _{r-j} (x)$, and so $\varphi_r$ is well-defined.

$\Mbar_R$ has a natural $G_R$-action defined by the diagonal action where $G_R$ acts on $\cO A_\cris (R)$ as before and trivially on $M$. We have an induced connection defined by $\nabla (a\otimes m) = \nabla_{\cO A_\cris (R)} (a) \otimes m+ a\otimes \nabla_M(m)$. 
Let
$$T_\cris (M) \coloneqq ( \rF ^r \Mbar_R )^{\varphi_r =1, \nabla = 0}(-r) = \{ x \in \rF^r \Mbar_R ~|~ \varphi _r(x) = x , \nabla (x) = 0\}(-r).$$ 
Note that $T_\cris(M)$ is a $\Z_p[G_R]$-module. 

\begin{rem}
The Tate twist $(-r)$ above assures that $T_\cris$ is independent of the choice of $r\leq p-2$ (see Remark \ref{rem:indep_r}).
 \end{rem}

We now give an alternate construction of $T_\cris (M)$, which will be more convenient for comparison under base change via the map $b_g$. Let $\Mbar \coloneqq  A_{\cris} (R) \otimes_{\lambda_1, R} M$ and  $\rF^r \Mbar \coloneqq \sum_{i = 0}^r \rF^i  A _\cris (R) \otimes_{\lambda_1, R} \rF ^{r-i} M \subset \Mbar$. Define $\varphi_r\colon \rF ^r \Mbar \to \Mbar$ similarly as above. 

\begin{rem} \label{rem:tors-vanish}
For any $M \in \mathrm{MF}_{\nabla}^{[0, r]}(R)$, we have $\mathrm{Tor}_R^1(\rF^i A_{\mathrm{cris}}(R), \rF^j(M)) = 0$ for each $i$ and $j$ since $A_{\mathrm{cris}}(R)$ is $p$-torsion free by \cite[Prop. 6.1.3]{brinon-relative} and $\rF^j(M)$ is a direct sum of modules of the form $R/p^a R$. 
\end{rem}

\begin{rem} \label{rem:Fr-graded}
For $M \in \mathrm{MF}_{\nabla}^{[0, r]}(R)$ and each $i \geq 0$, we define $\rF^i (A_{\mathrm{cris}}(R)\otimes_R M) = \mathrm{F}^i(\overline{M})$ in a similar manner. Choose an adapted basis as in Remark \ref{rem-adpated}. Then for each $i \geq 0$, 
\[ 
\rF^i (\overline{M}) = \bigoplus _{j = 1}^{m_0} \rF ^{r_j} A_\cris (R) / p ^{a_j} \rF ^{r_j} A_\cris (R) \cdot e_j
\]
where $r_j = \max\{0, i-h_j\}$ for $h_j \coloneqq \max\{ \ell ~|~ e_j \in \rF^\ell (M)\}$. So $\gr^i (\Mbar) = \bigoplus\limits_{j + \ell = i} \gr^j A _\cris (R) \otimes_R \gr^{\ell} M$.
\end{rem}

\begin{lemma} \label{lem-exact-fil} 
Given an exact sequence $0 \rightarrow M_1 \rightarrow M_2 \rightarrow M_3 \rightarrow 0$ in $\mathrm{MF}_{\nabla}^{[0, r]}(R)$,  the induced sequence 
\[
0 \rightarrow \rF^r(A_{\mathrm{cris}}(R)\otimes_R M_1) \rightarrow \rF^r(A_{\mathrm{cris}}(R)\otimes_R M_2) \rightarrow \rF^r(A_{\mathrm{cris}}(R)\otimes_R M_3) \rightarrow 0
\] 
is exact.     
\end{lemma}

\begin{proof}
First, note that the sequence $0 \rightarrow A_{\mathrm{cris}}(R)\otimes_R M_1 \rightarrow A_{\mathrm{cris}}(R)\otimes_R M_2 \rightarrow A_{\mathrm{cris}}(R)\otimes_R M_3 \rightarrow 0$ is exact by Remark \ref{rem:tors-vanish}. It follows that the sequence $0 \rightarrow \mathrm{F}^r(\overline{M_1}) \rightarrow \mathrm{F}^r(\overline{M_2}) \rightarrow \mathrm{F}^r(\overline{M_3}) \rightarrow 0$ is exact on the left and right by construction. It remains to show the exactness in the middle, i.e. $\overline{M_1} \cap \mathrm{F}^r(\overline{M_2}) = \mathrm{F}^r(\overline{M_1})$ as submodules of $\overline{M}_2$. 

By Theorem \ref{thm-Faltings1} (3), the sequence $0 \rightarrow \mathrm{gr}^i(M_1) \rightarrow \mathrm{gr}^i(M_2) \rightarrow \mathrm{gr}^i(M_3) \rightarrow 0$ is exact for each $i \geq 0$. Since $\mathrm{gr}^i(M_{\ell})$ for $\ell = 1, 2, 3$ is a direct sum of modules of the form $R/p^a R$ and $\mathrm{gr}^{j} A_\cris (R)$ for each $j \geq 0$ is $p$-torsion free by the proof of Lemma \ref{lem:graded-OAcris-p-tors-free}, we deduce from Remark \ref{rem:Fr-graded} that the sequence
\[
0 \rightarrow \mathrm{gr}^s(\overline{M_1}) \rightarrow \mathrm{gr}^s(\overline{M_2}) \rightarrow \mathrm{gr}^s(\overline{M_3}) \rightarrow 0
\]
for each $s \geq 0$ is exact. In particular, the map $\mathrm{gr}^s(\overline{M_1}) \rightarrow \mathrm{gr}^s(\overline{M_2})$ is injective for $s = 0, \ldots, r-1$. Thus, $\overline{M_1} \cap \mathrm{F}^r(\overline{M_2}) = \mathrm{F}^r(\overline{M_1})$.
\end{proof}

The projection $\pi\colon \cO A_\cris (R)\onto A_\cris (R)$ given by $X_i \mapsto 0$ (see Remark \ref{rem:comp-structure}) induces a map $\pi_M\colon \Mbar _R \onto \Mbar$, which is compatible with $\varphi_r$ and maps $\rF^r \Mbar_R$ to $\rF^r \Mbar$. 

\begin{lemma}\label{lem-Tcris-compatible} The morphism $\pi _M$ induces an isomorphism $T_\cris (M)\simeq (\rF^r \Mbar) ^{\varphi_r=1}$ of $\mathbb{Z}_p$-modules.
\end{lemma}

\begin{proof} Consider $ \Mbar_R^{\nabla= 0 } \coloneqq \{x \in \Mbar_R ~|~ \nabla (x) = 0\}$, which is an $A_\cris (R)$-module. Define $\rF^r \Mbar _R^{\nabla= 0} \coloneqq \rF ^r \Mbar _R \bigcap \Mbar _R^{\nabla= 0}$. Note that $\varphi_r$ maps $\rF^r \Mbar _R^{\nabla= 0}$ to $\Mbar _R^{\nabla= 0}$ since $\varphi_i$'s are parallel. It suffices to prove $\Mbar_R^{\nabla= 0 } \subset \Mbar_R$ is a section of $\pi_M$ so that $\pi_M\colon \Mbar_R^{\nabla = 0} \to \Mbar$ is an $A_\cris (R)$-linear isomorphism which induces an isomorphism on $\rF^r(-)$ and compatible with $\varphi _r$. To show this, for any $x \in M \subset \Mbar_R$, set 
\begin{equation}\label{eqn-xhat}
\hat x \coloneqq \sum _I \gamma _I (-X) \nabla (\partial)^I(x).  
\end{equation}
Here $I = (i_1, \dots , i_d)$ with $i_k \geq 0$ is a multi-index, and $\displaystyle \gamma_I (-X) = (-1)^{i_1 + \dots + i_d}\frac{X_1^{i _1}}{i _1!} \cdots\frac{X_d^{i _d}}{i _d!}$ with $X_i$ as in Lemma \ref{lem-2}, and 
 $\nabla (\partial)^I = \Theta_1^{i_1} \cdots \Theta_d^{i_d}$ with $\Theta_i = \nabla(\frac{\partial}{\partial t_i})\colon M \to M$ (note that $\Theta_i$ and $\Theta_j$ commute as $\nabla$ is integrable). Using $\Theta_j (\gamma_I (-X))= -\gamma_{I'} (-X)$ where $I' = (i _1, \dots, i_j -1, \dots , i_d)$, we see $\hat x \in \Mbar _R^{\nabla = 0}$.  Let $M ^{\nabla} \coloneqq \{ \hat x ~|~ x \in M\} \subset \Mbar_R ^{\nabla = 0}$. 
Clearly $\pi_M( M^\nabla ) = M \subset \Mbar$. For each $a\in R\subset \cO A_\cris (R) $ (here we use the embedding via $\lambda _2$), we have $ \hat a  = \lambda_1 (a) \in A_{\cris}(R) $ by considering the corresponding Taylor expansion. So for any $a \in R $ and $x\in M\subset \Mbar _R$, we have $\wh{ax}= \hat  a  \hat x= \pi(a) \hat x \in  \Mbar_R$ where $\pi(a) = \lambda_1(a) \in A_\cris(R)$. Therefore, $M^\nabla$ is an $R$-module via $\lambda_1\colon R \to A_\cris  (R)$, and $\pi_M\colon M ^\nabla \to M$ is an isomorphism of $R$-modules. So the composite of the maps $A_\cris (R) \otimes_{\lambda_1, R} M^{\nabla} \to \Mbar_R^{\nabla = 0} \subset \Mbar_R\overset{\pi_M}{\longrightarrow} \Mbar = A_\cris (R) \otimes_{R} M $ is an isomorphism of $A_\cris (R)$-modules. 

Now we claim that the natural map $ \iota\colon \cO A_{\cris}(R) \otimes_{\lambda_1, R} M^\nabla \to \Mbar_R$ is an isomorphism. Since $\cO A _\cris (R)$ is $p$-torsion free, we can reduce to the case where $M$ is killed by $p$ by d\'{e}vissage using the exact sequence $ 0 \to p M \to M \to M / p M \to 0 $ in $\MF_\nabla^{[0, r]} (R) $. So we can assume $M$ is finite $R/ p R$-free. Let $e_1 , \dots, e_s$ be a $R/ pR$-basis of $M$. Note that $(\hat e_1, \dots, \hat e_s) = (e_1, \dots , e_s) (I_s + Y)$ where $I_s$ is the identity matrix and $Y$ is a nilpotent matrix with entries in $\cO A _\cris (R)/ (p)$. So $\{\hat{e}_j\}$ is also a $\cO A _\cris (R)/ (p)$-basis of $\Mbar_R$, and $\cO A_{\cris}(R) \otimes_{\lambda_1, R} M^\nabla \simeq \Mbar_R $ via $\iota$. Taking $\nabla= 0$ on both sides, we deduce that $A_\cris (R) \otimes_{R} M^\nabla = \Mbar_R ^{\nabla = 0}$. Thus, $\Mbar_R^{\nabla = 0}$ is an $A_\cris (R)$-linear section of $\pi_M$.  

Since $\pi_M (\rF^r \Mbar_R) \subset \rF^r \Mbar$, we have $\pi_M (\rF^r \Mbar_R^{\nabla = 0}) \subset \rF^r \Mbar$. To show $\pi_M (\rF^r \Mbar_R^{\nabla = 0}) =\rF^r \Mbar$, it suffices to check that if $x \in \rF^r M$, then $\hat x \in \rF^r \Mbar _R$. This follows from Griffiths transversality and $\gamma_I (-X) \in \rF^{\sum i_j} \cO A_\cris (R)$. Finally, $\pi_M: \Mbar_R^{\nabla= 0} \to \Mbar$ is compatible with $\varphi_r$ since $\pi_M: \Mbar_R \to \Mbar$ is. This completes the proof. 
\end{proof}

We give an alternate description of $T_\cris (M)$ using previous Lemma. First define a $G_R$-action on $\Mbar = A_\cris (R) \otimes_{R} M$ via the isomorphism $\pi_M\colon \Mbar _R^{\nabla = 0} \to \Mbar$. This $G_R$-action on $\Mbar$ can be described explicitly as follows. For any $g \in G_R$, let $\beta_i (g) \coloneqq g([\wt {1+t_i}])-[\wt{1+t_i}]\in \rF^1 A_\cris (R)$ and set
   \begin{equation}\label{eqn-Gaction}
   	g (a \otimes x) = \sum _I g (a) \gamma_I (\beta(g)) \otimes \nabla (\partial)^I(x), \ \  \forall a \in A_\cris (R), \ \forall x \in M 
   	\end{equation}
where $\displaystyle \gamma _I (\beta(g)) \coloneqq \frac{\prod_l \beta_l (g)^{i_l}}{\prod_l i_l !} $. Since $g (X_i)= -\beta_i (g) + X_i$, we see from formula \eqref{eqn-xhat} that $\pi_M(g (\hat x)) = \sum _I \gamma_I (\beta(g)) \nabla (\partial)^I(x)$ for any $x \in M$ as required.  

Recall from \S \ref{subsec-period-ring} that $z_{i, n }$ is a fixed choice of $p^n$-th root of $1+ t_i$, which is used to define $\wt {1+t_i}$. Let $R_{\infty} = \bigcup_{n \geq 1} R[z_{1, n}, \ldots, z_{d, n}]$, which is a subextension in $\overline R$, and  $G_\infty \coloneqq \R{Gal}(\overline{R}[\frac{1}{p}]/R_\infty[\frac{1}{p}]]) \subset G_R$. By the proof of Lemma \ref{lem-Tcris-compatible}, the $G_R$-action defined on $\Mbar= A_\cris (R)\otimes_{R} M$ via $\Mbar_R^{\nabla = 0}$ satisfies the following properties:
 
\begin{itemize}
   \item The equation (\ref{eqn-Gaction}) gives a well-defined $A_\cris(R)$-semilinear $G_R$-action: for $m \in A_\cris(R) \otimes_{R} M$ and $a \in A_\cris(R)$ we have $g (am) = g (a) g(m)$. Moreover, $M \subset (A_\cris(R) \otimes_{R} M)^{G_\infty}$ as $\beta_i (g)= 0$ for $g \in G_\infty$;

  \item The $G_R$-action preserves filtration, i.e., $g (\rF^i(A_\cris (R)\otimes_R M))\subset \rF^i(A_\cris (R)\otimes_R M)$, and commutes with $\varphi_i$.
\end{itemize}

\noindent Hence, we have an alternate description of $T_\cris (M)$:
  	$$T_\cris (M) = (\rF^{r}(A_\cris(R) \otimes_R M))^{\varphi_r =1}(-r) = \{x\in \rF^{r}(A_\cris(R) \otimes_R M) ~|~ \varphi_r (x) = x \}(-r). $$

We say an object $M \in \R{MF}(R)$ and a $\Z_p$-module $L$ have the \textit{same type} if $M \cong L\otimes_{\Z_p} R$ as $R$-modules.

\begin{theorem} \label{thm-Faltings2}
Let $M \in \R{MF}^{[0,r]}_\nabla(R)$. Then $T_\cris (M)$ is a finite $\Z_p$-module of the same type as $M$, and the map $\varphi_r -1\colon \rF^{r}(A_\cris (R)\otimes_R M) \to A_\cris (R)\otimes_R M$ is surjective. Furthermore, the functor $T_\cris(-)$ is exact. 
\end{theorem}

\begin{proof}
We follow the proof of \cite[Thm. 2.4]{FaltingsCcris}, with some minor modifications as our $T_\cris(-)$ is covariant. Let $M[p] \coloneq \ker(M \stackrel{\times p}{\longrightarrow} M)$. By Corollary \ref{cor:MF(R)-abelian-category}, $M[p]$ with the induced structure is a Fontaine--Laffaille module, and we have an exact sequence
\[
0 \rightarrow M[p] \rightarrow M \stackrel{\times p}{\longrightarrow} pM \rightarrow 0
\]
of objects in $\R{MF}^{[0,r]}_\nabla(R)$. Consider the induced diagram
\[
\xymatrix{ 0 \ar[r] & T_\cris (M[p])(r)\ar[d] \ar[r] & T_\cris (M)(r) \ar[d] \ar[r]^{\times p} & T_\cris (pM)(r)\ar[d]\ar[r] & 0\\ 0 \ar[r] & \rF^r  (A_\cris (R)\otimes_R M[p]) \ar[d]^{\varphi _r-1}\ar[r] & \rF^r (A_\cris (R) \otimes_R M)\ar[d] ^{\varphi_r -1} \ar[r]^{\times p} & \rF^r (A_\cris (R) \otimes_R pM)\ar[r]\ar[d]^{\varphi_r -1} & 0 \\ 0 \ar[r] & A_\cris (R) \otimes_R M[p] \ar[r] & A_\cris (R)\otimes_R M \ar[r]^{\times p} & A_\cris (R)\otimes_R pM \ar[r] & 0. }
\]
Note that the middle and bottom rows are exact by Lemma \ref{lem-exact-fil} and Remark \ref{rem:tors-vanish}. By d\'{e}vissage, to prove that $\varphi_r-1$ is surjective, it suffices to consider the case $pM = 0$. Furthermore, assuming $\varphi_r-1$ is surjective, we deduce by the snake lemma that the top row is exact. So by Theorem \ref{thm-Faltings1}, if the $\mathbb{Z}_p$-modules $T_\cris (M[p])$ and $T_\cris (pM)$ are finitely generated and have the same type as $M[p]$ and $pM$ respectively, then $T_\cris (M)$ is a finite $\mathbb{Z}_p$-module having the same type as $M$ by the structure theorem of finite $\mathbb{Z}_p$-modules. Thus, it suffices to prove the statements for the case $pM = 0$.  

Suppose now $pM = 0$. To compute $T_{\mathrm{cris}}(M) = (\rF^r (M \otimes_R A_\cris(R)/pA_\cris(R)))^{\varphi _r =1}$, we first replace $A_\cris(R)/pA_\cris(R)$ by a simpler object as follows. By Lemma \ref{lem-1}, $$A_\cris (R)/(p) \simeq \Rbar^\flat / \tilde p ^p  \Rbar^\flat [\overline{\delta}_0, \overline{\delta}_1, \dots ] / (\overline{\delta}_0^p, \overline{\delta}_1^p, \ldots) ,$$ where $\overline{\delta}_i \coloneqq \delta_i \mod p$. Note that for any $x \in M$ and $i \geq 0$, we have $x\otimes \overline{\delta}_i \in \rF^r (M \otimes_R A_\cris(R)/(p))$ since $\delta_i \in \mathrm{F}^p A_\cris(R)$ and $r \leq p-2$. Furthermore, $\varphi_r(m\otimes \overline{\delta}_i) = 0$ as $\overline{\delta}_i^p = 0$. Since $M$ is finite free over $R_1 = R/pR$, we deduce $(\rF^r (M \otimes_R A_\cris(R)/(p)))^{\varphi _r =1} = (\rF^r (M \otimes _R \Rbar^\flat / \tilde p ^p \Rbar^\flat)) ^{\varphi _r =1}$. Equip $\Rbar_1 \coloneqq \overline{R}/p\overline{R}$ with an $R_1$-module structure via
\[
R/pR \rightarrow R_{\infty}/pR_{\infty} \stackrel{\varphi^{-1}}{\longrightarrow} R_{\infty}/pR_{\infty} \rightarrow \Rbar_1,
\]
and a filtration by setting $\rF^i(\Rbar_1) \coloneqq p^{i}_1\Rbar_1$ (recall $p_1^p = p$). Define $\varphi_i$ on $\rF^i(\Rbar_1)$ by $\varphi_i(x) = (-1)^iy^p$ for $x = p_1^iy$. This gives $\Rbar_1$ a structure of an object in $\R{MF}_{\R{big}}(R)$. It is easy to check that $\theta\circ \varphi ^{-1}\colon \Rbar^\flat / \tilde p ^p  \Rbar^\flat \to \Rbar_1$ is an isomorphism compatible with all structures. Hence, it suffices to show that $\dim_{\mathbb F_p} (\rF^r (M \otimes _R \Rbar_1)^{\varphi _r =1})= {\mathrm rank}_{R_1} M$ and that $\varphi_r -1\colon \rF^r (M \otimes_R \Rbar_1 ) \to M \otimes_R \Rbar_1$ is surjective. For this, we choose an adapted basis of $M$ and compute the Jacobians of associated equations as follows.

Fix an adapted $R_1$-basis $e_1, \ldots , e_m$ of $M$ and integers $n_1 , \ldots , n_m$ such that $n _i \leq r$ and $\{p_1^{n _i} e_i, ~i =1, \ldots, m\}$ form a basis for $\rF^r (M \otimes_R  \Rbar _1)$. We have  $\varphi _r (p_1^{n _1} e_1 , \dots , p _1 ^{n _m} e_m) = (e_1, \dots , e_m ) A$ for some $A \in {\mathrm GL}_m (R_1)$ (note that $\varphi _i (p_1^i )= (-1)^i$). In particular, $x = \sum _i x_i p_1^{n _i} e_i \in \rF^r (M \otimes _R \Rbar_1)^{\varphi _r =1}$ if and only if $X = (x_1, \ldots, x_m)^T$ satisfies the equation $A X^{p} = \Lambda X$, where $X^p \coloneqq (x_1^p, \ldots x_m^p)^T$ and $\Lambda = \langle p_1 ^{n_1}, \ldots , p_1 ^{n_m} \rangle$ denotes the diagonal matrix with $p_1^{n_1}, \ldots , p_1^{n_m}$ on the main diagonal. 

Choose a lift of $A$ to an element in ${\mathrm GL}_m(R)$, which we also denote by $A$ abusing the notation. Consider the equation $A X^p = \Lambda X$ defined over $R[p_1]$. As explained in \cite[Thm. 2.4 Pf.]{FaltingsCcris}, we note that the set of solutions for $AX^p = \Lambda X$ in $\overline{R}$ corresponds bijectively to that of $A X^p \equiv  \Lambda X \mod p$ in $\overline{R}_1$. Indeed, if $X_0$ is a solution of $AX^p \equiv \Lambda X \mod p$, then we can find $Y$ such that $X_0 + p_1^2 Y$ satisfies the equation $AX^p \equiv \Lambda X \mod p_1^{p+2}$ since $n_i \leq p-2$. We can construct inductively a compatible system of solutions to $AX^p \equiv \Lambda X \mod p_1^{p + 2\ell}$ in a sufficiently big finite $R$-subalgebra of $\overline{R}$, and therefore obtain a solution $X$ of $AX^p = \Lambda X$. Furthermore, such a lift $X$ of $X_0$ is unique since $\overline{R}$ is $p$-adically separated by \cite[Prop. 2.0.3]{brinon-relative}.  

Let $B = R[p_1][x_1, \ldots, x_m]/ (AX^p - \Lambda X)$. Note that $B$ is finite as an $R$-module. The Jacobian of $AX^p-\Lambda X$ is $pA \langle x_1^{p-1}, \dots , x_m ^{p-1} \rangle - \Lambda$, which can be written in the form $-\Lambda (I_m - p\Lambda^{-1}A')$ where $p\Lambda^{-1}A'$ has entries in $p_1^2 R[p_1][x_1, \ldots, x_m]$. The determinant of the Jacobian is a unit in $B[\frac{1}{p}]$. Hence, $B[\frac{1}{p}]$ is \'etale over $R[\frac 1 p]$, and all the solutions $X$ for $AX^p = \Lambda X$ lie inside $\Rbar$. Since the set of solutions of $AX^p = \Lambda X $ has size $p^m$, we have $\dim_{\mathbb F_p} (\rF^r (M \otimes _R \Rbar_1)^{\varphi _r =1})= {\mathrm rank}_{R_1} M$. 
  
To show that $\varphi_r -1$ is surjective, we apply a similar argument as above to the equation $AX^p = \Lambda X + Y$ with $Y$ having fixed entries in $\Rbar$. Choose a finite $R$-algebra $S$ containing elements in $Y$ such that $S[\frac 1 p]$ is \'etale. Then $ S[p_1][\frac 1 p][x_1 , \dots , x_m ]/ (A X^p - \Lambda X-Y)$ is \'etale over $S[p_1][\frac 1 p]$ by a similar calculation of the Jacobian. 

It remains to prove the exactness of $T_{\mathrm{cris}}(-)$. By Remark \ref{rem:tors-vanish}, given an exact sequence $0 \rightarrow M_1 \rightarrow M_2 \rightarrow M_3 \rightarrow 0$ in $\mathrm{MF}_{\nabla}^{[0, r]}(R)$, the induced sequence $0 \rightarrow A_{\mathrm{cris}}(R)\otimes_R M_1 \rightarrow A_{\mathrm{cris}}(R)\otimes_R M_2 \rightarrow A_{\mathrm{cris}}(R)\otimes_R M_3 \rightarrow 0$ is exact. Furthermore, the sequence $0 \rightarrow \rF^r(A_{\mathrm{cris}}(R)\otimes_R M_1) \rightarrow \rF^r(A_{\mathrm{cris}}(R)\otimes_R M_2) \rightarrow \rF^r(A_{\mathrm{cris}}(R)\otimes_R M_3) \rightarrow 0$ is exact by Lemma \ref{lem-exact-fil}. Since $\varphi_r-1$ is surjective, we deduce that 
\[
0 \rightarrow T_{\mathrm{cris}}(M_1) \rightarrow T_{\mathrm{cris}}(M_2) \rightarrow T_{\mathrm{cris}}(M_3) \rightarrow 0.
\] 
is exact. 
\end{proof}

\begin{rem} 
One can also show that $T_\cris(-)$ is fully faithful and stable under subquotients.
\end{rem}	
	
Let $M \in \MF^{[0, r]}(R)$ and recall that $M_g \coloneqq M\otimes_{R, b_g} R_g$ is an object in $\MF^{[0, r]}(R_g)$. We have a natural $G_{R_g}$-equivariant morphism $T_\cris (M)\to  T_\cris (M_g)$.

\begin{corollary}\label{cor-abstract-comp}
The natural map $T_\cris (M) \to T_\cris (M_g)$ is an isomorphism of $G_{R_g}$-representations. 

\end{corollary}
\begin{proof} By Theorem \ref{thm-Faltings2}, both $T_\cris (M)$ and $ T_\cris (M_g)$ have the same type as $M$. Therefore, it suffices to show that the map $T_\cris (M)\to  T_\cris (M_g)$ is injective. By d\'{e}vissage using the exact sequence $0 \rightarrow M[p] \rightarrow M \rightarrow pM \rightarrow 0$ as in the proof of Theorem \ref{thm-Faltings2}, we may assume $pM = 0$. In this case, by the proof of Theorem \ref{thm-Faltings2}, $T_\cris (M)(r)$ can be identified with the set 
	$\{X\in \Rbar ^m ~|~ AX^p = \Lambda X\}$. Similarly, $T_\cris (M_g)(r)$ can be identified with the set $\{X\in \overline{R_g} ^m ~|~ AX^p = \Lambda X\}$. Since the map $\overline{b_g}\colon \Rbar \to \overline{R_g}$ is injective, we deduce that $T_\cris (M)\to  T_\cris (M_g)$ is injective. 
\end{proof}

Let $\tilde{\epsilon} \coloneqq (\epsilon_n)_{n \geq 0} \in \mathbb{Z}_p(1)$ with $\epsilon_0 = 1$, $\epsilon_{n+1}^p = \epsilon_n$, and $\epsilon_1 \neq 1$. Recall that $\beta(\tilde \epsilon) = \log ([\tilde \epsilon])$ for the map $\beta\colon \Z_p(1) \rightarrow \rF^1(A_{\R{cris}}(R))$. Note that $\varphi_1 (\beta(\tilde \epsilon)) = \beta(\tilde \epsilon)$. 

\begin{rem}\label{rem:indep_r}
 The construction $T_\cris (M)$ does not depend on the choice of $r$. That is, if $r \leq s \leq p-2$, we have 
$$(\rF^{r}(A_\cris(R) \otimes_R M))^{\varphi_r =1}(-r) \simeq (\rF^{s}(A_\cris(R) \otimes_R M))^{\varphi_s =1}(-s).$$
To see this, note that $\beta(\tilde{\epsilon})^{s-r}\cdot (\rF^{r}(A_\cris(R) \otimes_R M) ) ^{\varphi_r =1} \subset (\rF^{s}(A_\cris(R) \otimes_R M)) ^{\varphi _s =1}$. Since both sides have the same type as $M$ by Theorem \ref{thm-Faltings2}, $\beta(\tilde{\epsilon})^{s-r}\cdot (\rF^{r}(A_\cris(R) \otimes_R M) ) ^{\varphi_r =1} = (\rF^{s}(A_\cris(R) \otimes_R M)) ^{\varphi _s =1}$ as $\mathbb{Z}_p$-modules. As $\mathbb{Z}_p[G_R]$-modules, we have $\beta(\tilde{\epsilon})^{s-r}\cdot \mathbb{Z}_p \simeq \Z_p (s-r)$.  
\end{rem}

\section{Fontaine--Laffaille data and torsion crystalline cohomology}\label{sec-cohoisFL}

Let $X$ be a smooth proper scheme over $R = W(k)[\![t_{1},\ldots,t_{d}]\!]$, and for a fixed integer $n \geq 1$, let $X_{n}$ denote its base change over $R_n = R/ p^n R$. The ideal $(p) \subset R$ has a canonical divided power structure, and we consider $R_n$ with the corresponding quotient divided power structure. Below, we let $M^i \coloneqq \rH^{i}_\cris (X_n/ R_n, \cO_{X_n /R_n })$ and $\rF^ j (M^i) \coloneqq \rH^i_\cris ( X_n/ R_n , J ^{[j]}_{X_n/R_n})$ denote the corresponding crystalline cohomology groups. One of our main results is to show that $(M^i,\rF^j(M^i))$ for $i \leq p-2$ gives rise to an object of $\MF_{\nabla}^{[0,p-2]}(R)$. More precisely, we have the following theorem.
 
\begin{theorem}\label{thm-cons-phi}
Let $X$ be a smooth proper scheme over $R$.  
\begin{enumerate}
\item[(1)] For $i \leq p-1$, there exist $\varphi_R$-semilinear maps $\varphi_j\colon \rF^ j (M^i) \rightarrow M^i$ and a connection $\nabla_j\colon \mathrm{F}^j(M^i) \to \mathrm{F}^{j-1}(M^i) \otimes_R \wh \Omega_R$ such that $(M^i,  \rF^j(M^i), \varphi_j, \nabla_j)$ is an object of $\R{MF}_{\nabla, \R{big}}(R)$.
\item[(2)] { Let $k \subset k'$ be a field extension where $k'$ is perfect, $R' \coloneqq W(k')$ equipped with the Witt vector Frobenius, and suppose we have a morphism $R \rightarrow R'$ of divided power rings compatible with Frobenius. Let $X' \coloneqq X \times_R R'$. Then the natural morphism (induced by functoriality) $$\rH^{i}_\cris (X_n/ R_n, \cO_{X_n /R_n }) \to \rH^{i}_\cris (X'_n/ R'_n, \cO_{X'_n /R'_n })$$ is compatible with filtration and $\varphi_j$.}
\item[(3)] For $i \leq p-2$, the data $(M^i, \rF^j(M^i), \varphi_j)$ is an object of $\MF^{[0, i]}(R)$. Hence, $(M^i,  \rF^j(M^i), \varphi_j, \nabla_j)$ is an object of $\MF_{\nabla}^{[0, i]}(R)$.
\end{enumerate}
\end{theorem}

We shall prove Theorem \ref{thm-cons-phi} (1) (2) in section \S \ref{subsec-divpowersheaf} (see e.g. Corollary \ref{cor:firstpartmainthm}). In this section, we prove part (3) of the theorem assuming parts (1) and (2). We first study the case $d=1$ in the following subsection, and then consider the general case via induction.

\subsection{The case $d=1$} In this subsection, we study the special case $R = W(k)[\![t_1]\!]$. We begin with some preliminary remarks. Let $\bar b\colon R \to W(k)$ denote the natural quotient morphism and $X_{s} \coloneqq X \times_{R, \bar b} W(k)$. By the base change theorem for crystalline cohomology (cf. Theorem \ref{thm-cry-BC} and Example \ref{Ex-BC}), $\bar{b}$ induces a natural injective map:

\begin{equation}\label{eq-modt2}
\bar \alpha\colon \rF^j (M^i)/ t_1 \rF^j (M^i) \hookrightarrow \rF^j ( M_s^i)\coloneqq \rH^i_\cris (X_{s, n}/W_n (k), J^{[j]}_{X_{s, n}/W_n(k)}).
\end{equation}

\begin{proposition} \label{prop-alphabar-suj} With notations as above, $\bar \alpha$ is an isomorphism for all $0 \leq i \leq p -2$.
\end{proposition}

Before proving the proposition, we explain how to deduce Theorem \ref{thm-cons-phi} (3) assuming the proposition. Let $b_g\colon R \to R_{g}$ be as before, and let $X_g \coloneqq X \times_{R}  R_{g}$. Since $b_g$ is flat, Theorem \ref{thm-cry-BC} implies that the induced map
 
\begin{equation}\label{eq-tensor-g} \alpha_g\colon \rF^j (M^i)\otimes_{R_n} R_{g, n} \rightarrow \rF^j (M^i_{g})\coloneqq \rH_\cris^i (X_{g, n} / R_{g, n}, J^{[j]}_{X_{g, n} /R_{g, n}})
\end{equation}
is an isomorphism.

\begin{rem}
The maps $\alpha_g$ and $\bar \alpha$ are compatible with $\varphi_j$ by Theorem \ref{thm-cons-phi} (2).  
\end{rem}

\begin{proof}[P\textbf{roof of Theorem} \ref{thm-cons-phi} (3) \textbf{when} $d=1$]
By Proposition \ref{prop-detect} and Remark \ref{rem-filrange}, it suffices to prove that $M^{i}\otimes_R R_g$ and $M^i\otimes_R W(k)$ with the induced $\varphi_j$ and filtrations are objects in $\MF(R_g)$ and $\MF(W(k))$ respectively, and that they have the same type. Since $\alpha_g$ is an isomorphism, it follows from Theorem \ref{thm-cons-phi} (1), (2) that $M^{i}\otimes_R R_g \simeq \rH^{i}_\cris (X_{g, n}/ R_{g, n})$ as objects of $\MF_{\text{big}} (R_g)$. On the other hand, the latter object is in $\MF(R_g)$ by classical Fontaine--Messing theory (over $R_g$) summarized in Theorem \ref{Thm-FM}. Similarly, by Proposition \ref{prop-alphabar-suj}, $M^i\otimes_R W(k)$ is in $\MF(W(k))$. On the other hand, $M^i\otimes_R W(k)$ (resp. $M^{i}\otimes_R R_g$) has the same type with  $\rH^i _\et (X_{s, \overline K}, \Z/p^n \Z )$ (resp.  $\rH^i _\et (X_{g, \overline{R_g}[\frac{1}{p}]}, \Z/p^n \Z )$) again by Theorem \ref{Thm-FM} and Theorem \ref{thm-Faltings2}. Since $\rH^i_\et (X_{s, \overline K}, \Z/p ^n \Z )\simeq \rH^i_\et (X_{g, \overline{R_g}[\frac{1}{p}]}, \Z/p ^n \Z)$ as $\mathbb{Z}_p$-modules by smooth and proper base change theorems, we conclude that $M^{i}\otimes_R R_g$ and $M^i\otimes_R W(k)$ have the same type.
\end{proof}

To prove Proposition \ref{prop-alphabar-suj}, we need the following lemma. For a ring $A$ and a finite length $A$-module $N$, we denote by $\ell_A(N)$ the $A$-length of $N$. 

\begin{lemma}\label{lem-length-s-g} For any finite $R_n$-module $M$, we have $$\ell_{W(k)} ( M \otimes_R W(k)) \geq \ell _{W(k_g)} (M \otimes_R W(k_g)) .$$ 
\end{lemma}

\begin{proof} We shall proceed via induction on the minimal integer $h$ such that $p^h M[\frac{1}{t_1}] = 0$. If $h = 0$, then $M\otimes_R W(k_g) = 0$ and the desired inequality is trivial. If $h = 1$, then $M[\frac{1}{t_1}] \cong (M/pM)[\frac{1}{t_1}]$. Since $M/pM$ is a finite module over $R_1 = k[\![t_1]\!]$, we have $M/pM \simeq N_\R{tor} \oplus N_{\R{free}}$ where $N_{\R{tor}}$ is killed by a power of $t_1$ and $N_{\R{free}}$ is a free $k[\![t_1]\!]$-module. Thus, 
\[
\ell_{W(k)} (M \otimes_R W(k)) \geq \ell_{W(k)} ((M/pM) \otimes_R W(k)) \geq \ell_{W(k)}(N_{\R{free}} \otimes_R W(k)) = \ell_{W(k_g)} (M\otimes_R W(k_g)).
\]
Here the last equality follows from $M\otimes_R W(k_g) \cong (M/pM)\otimes_R W(k_g) \cong N_{\R{free}}\otimes_R W(k_g)$.

For $h \geq 2$, consider the map $f\colon M \to (M/pM)[\frac{1}{t_1}]$ and the associated exact sequence $0 \to \ker(f) \to M \to f(M) \to 0$. Note that $\ker(f)[\frac{1}{t_1}]$ is killed by $p^{h-1}$, and $f(M)$ has no $t_1$-torsion. In particular, the following sequence
\[
0 \to \ker(f)\otimes_R W(k) \to M\otimes_R W(k) \to f(M)\otimes_R W(k) \to 0
\]
is exact. By the inductive hypothesis, we have
\begin{eqnarray*}
\ell_{W(k)} (M\otimes_R W(k)) &=& \ell_{W(k)} (\ker(f)\otimes_R W(k))+ \ell_{W(k)} (f(M)\otimes_R W(k)) \\ &\geq&    \ell_{W(k_g)} (\ker(f) \otimes_R  W(k_g))+ \ell_{W(k_g)} (f(M) \otimes_R W(k_g)) \\ &=& \ell _{W(k_g)} (M\otimes_R W(k_g)).
\end{eqnarray*}
\end{proof}

\begin{proof}[P\textbf{roof of Proposition} \ref{prop-alphabar-suj}]
By crystalline base change theorem (cf. Example \ref{Ex-BC}), $\bar{\alpha}$ is an isomorphism if and only if $\rF^j (M^{i+1})$ is $t_1$-torsion-free. We first treat the case $j = 0$. Since $\bar \alpha$ is injective, Lemma \ref{lem-length-s-g} implies that $$\ell_{W(k)} (M_s^i) \geq \ell _{W(k)} (M^i/ t_1 M^i) \geq \ell_{W(k_g)} (M^i_g).$$ 
On the other hand, by Theorem \ref{Thm-FM} for $X$ base changed to $W(k)$ and $R_g$, we have for $0 \leq i \leq p-2$ that 
	$$\ell _{W(k)} ( M_s^i) =  \ell_{\Z_p} (\rH^i_\et(X_{s, \overline K}, \Z/ p^n\Z))= \ell _{\Z_p} (\rH^i_\et(X_{g, \overline{R_g}[\frac 1 p]}, \Z/ p^n\Z))=  \ell_{W(k_g)} (M_g ^i). $$
Thus, $\ell_{W(k)} (M_s^i) = \ell _{W(k)} (M^i/ t_1 M^i)$ and $\bar \alpha$ is bijective. 	

Now consider the case $j>0$. We claim that $\rF^j(M^{i+1})$ is $t_1$-torsion-free for $0 \leq i \leq p-2$. Suppose otherwise. In that case, since $M^{i+1}$ is $t_1$-torsion-free, the map $\rF^j (M^{i+1}) \to M^{i+1}$ is not injective. Then the induced map $\rF^j (M^{i+1})/t_1 \rF^j (M^{i+1}) \to M^{i+1} / t_1 M^{i+1}$ is also not injective, since $M^{i+1}$ is $t_1$-torsion-free and $\rF^j (M^{i+1})$ is $t_1$-adically separated. Consider the commutative diagram: 
$$\xymatrix{  \rF^j (M^{i+1})/ t_1 \rF^j (M^{i+1})   \ar@{^{(}->}[r] \ar[d] & \rF^j (M_s^{i+1}) \ar@{_{(}->}[d] \\  M^{i+1}/ t_1 M^{i+1} \ar@{^{(}->}[r] &  M^{i+1}_s}  $$
Since $(M_s^l , \rF^j (M_s^l), \varphi_j)$, for $0 \leq l \leq p-1$, is a Fontaine--Laffaille module over $W(k)$ by \cite[Cor. II.2.7]{FontaineMessing}, the right column is injective. But this contradicts that the left column is not injective, which proves the claim.
\end{proof}

\begin{rems} The previous proof uses that $(M_s^{p-1}, \rF^j (M_s^{p-1}), \varphi_j)$ is a Fontaine--Laffaille module. It is not known whether $T_\cris (M_s^{p-1}) $ is isomorphic to $\rH^{p-1}_\et (X_{s, \overline K}, \Z/ p^n\Z)$, which we do not need here. 
\end{rems}

\begin{rems} The proof of Proposition \ref{prop-alphabar-suj} \emph{does not} use the $\varphi_j$-structure on $\rF^j (M^i)$. In particular, the proof does not use the $\varphi$-structure on $R$. This is important for the proof of Proposition \ref{prop-multi-key} which proceeds via induction on $d$ and base changing along maps $R \to W(k _{\mathfrak q})[\![s_1 , \dots , s_m ]\!]$ with $m \leq d-2$ to reduce to the case $d=1$. In the constructions below, those base change maps may not be compatible with Frobenius. But we can still do the induction since Proposition \ref{prop-alphabar-suj} does not use $\varphi$ on $R$. 
\end{rems}

\subsection{The general case} \label{subsection:crys-cohom-base-change-general-case}

We consider the general case $d \geq 1$. Consider the base change map $R \to R_d \coloneqq R/t_d R = W(k) [\![t_1, \dots t_{d-1}]\!]$, and let $X_{d} = X\times_{R} R_{d}$. By the crystalline base change theorem (cf. Theorem \ref{thm-cry-BC} and Example \ref{Ex-BC}), for each $j \geq 0$, we have a long exact sequence
\begin{equation}\label{eq-long-seqence}
{\scriptstyle \cdots \longrightarrow \rH_\cris ^i (X_n/ R_n, ~J^{[j]} _{X_n/ R_{n}}) \overset{t_d}{\longrightarrow} \rH_\cris ^i (X_n/ R_{n}, ~J^{[j]} _{X_n/ R_{n}}) \longrightarrow \rH_\cris ^{i} (X_{d, n}/ R_{d, n}, ~J^{[j]} _{X_{d, n}/ R_{d, n}}) \longrightarrow \rH_\cris ^{i+1} (X_n/ R_n, ~J^{[j]} _{X_n/ R_{n}}) \longrightarrow \cdots}
\end{equation}
Let $\rF^j (M^i)$ and $\rF^j (M_d^i)$ denote  $\rH_\cris ^i (X_n/ R_{n}, J^{[j]} _{X/ R_{n}})$  and  $\rH_\cris ^{i} (X_{d, n}/ R_{d, n}, J^{[j]} _{X_{d, n}/ R_{d,n}})$ respectively. We show the following analog of Proposition \ref{prop-alphabar-suj}.

\begin{proposition} \label{prop-multi-key} For $0 \leq i \leq p-2$,  $\rF^j(M^i)/ t_d \rF^j (M^i) \cong \rF^j (M_d^i)$; or equivalently, $\rF^j (M^{i+1})$ is $t_d$-torsion-free.	
\end{proposition} 

Given the previous proposition, the proof of Theorem \ref{thm-cons-phi} (3) (for any $d \geq 1$) proceeds exactly as in the case $d = 1$, which will not be repeated here. Note that if Proposition \ref{prop-multi-key} holds for $t_d$, then the analogous statement with $t_d$ replaced by any of the $t_i$'s will also hold. We shall proceed via induction on $d$, with the base case $d = 1$ proved in the previous subsection. Assume $d \geq 2$. 

The long exact sequence \eqref{eq-long-seqence} gives an injective map
$\bar \alpha\colon \mathrm{F}^j(M^i)/ t_d \mathrm{F}^j(M^i) \inj \mathrm{F}^j(M_d^i)$, and we have
$$\coker (\bar \alpha) \cong N \coloneqq \{x \in \mathrm{F}^j(M^{i+1}) ~|~ t_d x = 0\}.  $$
Note that $N$ is an $R_d$-module. Let $U$ be the $p$-adic completion of $R_d[\frac{1}{t_1}] $. The following lemma is a key input.
 
\begin{lemma}\label{lem:equals-zero}
With notations as above, $U \otimes _{R_d} N = 0$.
\end{lemma}

Before proving the lemma, we explain how to show the proposition assuming the lemma. Let $N' \coloneqq \{x \in \mathrm{F}^j(M^{i+1}) ~|~ t_d^m  x = 0 \text{ for some } m\}$. Then $N \subset N ' \subset \mathrm{F}^j(M^{i+1})$. Furthermore, the map $N'/t_d N' \rightarrow \mathrm{F}^j(M^{i+1})/t_d \mathrm{F}^j(M^{i+1}) \subset \mathrm{F}^j(M_d^{i+1})$ is injective.

\begin{proof}[P\textbf{roof of Proposition} \ref{prop-multi-key}] We need to show $N = 0$. Suppose $N \neq 0$. Note that  $N$ is $p$-power torsion and $U\otimes_{R_d} N = 0$, we have $t_1^s\cdot N = 0$ for some integer $s \geq 1$. We claim that $N'$ is also killed by a power of $t_1$, which holds by the following fact (suggested by anonymous referee): Let $A$ be a commutative ring, $x, y \in A$, and $M$ an $A$-module killed by a power of $x$. Suppose that every element in $M[x] \coloneqq \{m \in M ~|~ xm = 0\}$ is killed by a power of $y$. Then each $m\in M$ is killed by a power of $y$. We induct on $s \geq 1$ such that $M$ is killed by $x^s$. If $s=1$, then $M = M[x]$ and the assertion holds. For general $s$, note that $xM$ is killed by $x^{s-1}$. By inductive hypothesis, for each $m \in M$, $xm$ is killed by $y^l$ for some $l \geq 1$. Then $y^l m \in M[x]$, so it is killed by $y ^{l'}$ for some $l' \geq 1$. Thus, $y^{l + l'} m = 0$ as required. 

Since $N'/t_dN' \rightarrow \mathrm{F}^j(M^{i+1}_d)$ is injective, we have contradiction to the inductive hypothesis that $\mathrm{F}^j(M_d^{i+1}) = \rH^{i+1}_\cris (X_{d, n}/R_{d, n} , J^{[j]} _{X_{d, n}/ R_{d,n}})$ is $t_1$-torsion-free. Hence, $N= 0$.
\end{proof}

To prove the lemma, we apply base change maps $b\colon R \to A$ to the exact sequence \eqref{eq-long-seqence} for various rings $A$ as follows. Abusing notation, we denote $(X\times_R A) / A$ by $X/ A$. For any maximal ideal $\mathfrak{q}$ of $U$, let $\hat{U}_{\mathfrak{q}}$ be the $\mathfrak{q}$-completion of $U$. We will use the following lemma to prove Lemma \ref{lem:equals-zero}.

\begin{lemma} \label{lem-structrue}
\begin{enumerate}
 \item[(1)] $\hat{U}_{\mathfrak{q}} \simeq \cO_{\mathfrak q} [\![s_1, \dots , s_m]\!]$ where $\cO_{\mathfrak{q}}$ is a Cohen ring whose maximal ideal is $(p)$ and $m \leq d-2$. There exists a Frobenius endomorphism $\varphi$ on $\mathcal{O}_{\mathfrak{q}}$ lifting the $p$-power Frobenius on $\mathcal{O}_{\mathfrak{q}}/(p)$.

 \item[(2)] Let $k_{\mathfrak{q}} \coloneqq \varinjlim_{\varphi} \mathcal{O}_{\mathfrak{q}}/(p)$. Then there exists a unique map $f\colon \cO_{\mathfrak{q}} \rightarrow W(k_{\mathfrak{q}})$ which lifts the map $\mathcal{O}_{\mathfrak{q}} \rightarrow k_{\mathfrak{q}}$ and is compatible with (given choices of) Frobenius. Furthermore, $f$ is faithfully flat.  
\end{enumerate}
\end{lemma}

\begin{proof}
Note that $\hat{U}_{\mathfrak{q}}$ is a complete regular local ring such that $(p)$ is a prime ideal. Hence, by the Cohen structure theorem, we have $\hat{U}_{\mathfrak{q}} \simeq \cO_{\mathfrak{q}}[\![s_1, \dots , s_m]\!]$ for some Cohen ring $\mathcal{O}_{\mathfrak{q}}$. The Krull dimension of $\hat{U}_{\mathfrak{q}}$ is less than $d-1$, so $m \leq d-2$. Moreover, $R/pR$ has a $p$-basis given by $\{t_1, \ldots, t_d\}$, which is also a $p$-basis of $\hat{U}_{\mathfrak{q}}$. In particular, $\mathcal{O}_{\mathfrak{q}}$ has a finite $p$-basis, and thus there exists a lift of Frobenius endomorphism $\varphi\colon \mathcal{O}_{\mathfrak{q}} \rightarrow \mathcal{O}_{\mathfrak{q}}$. This proves (1).

For (2), since $\mathcal{O}_{\mathfrak q}$ is $p$-torsion free, we obtain a unique lift $f\colon \cO_{\mathfrak q} \inj W(k_{\mathfrak{q}})$ commuting with given choices of Frobenius by Cartier's Dieudonn\'e--Dwork lemma. Since $W(k_{\mathfrak{q}})$ is $p$-torsion-free and $f$ is a local homomorphism, $f$ is faithfully flat. 	
\end{proof}

\begin{proof}[P\textbf{roof of Lemma} \ref{lem:equals-zero}]
Consider the base change of \eqref{eq-long-seqence} along the morphism $b\colon R \to U[\![t_d]\!]$ given by sending $t_i$ to $t_i$ for each $i$. Since the map $b$ is flat, we obtain 
 $$\bar \alpha_{b} = U[\![t_d]\!] \otimes_{b , R} \bar \alpha\colon  \rH ^i_\cris  (X/U_n[\![t_d]\!], J^{[j]})/ t_d \inj \rH ^i_\cris (X/U_n, J^{[j]}). $$
It suffices to prove that $\bar \alpha_{b}$ is an isomorphism, by the functoriality of the base change map applied to the cartesian square:
$$
\xymatrix{
X_U \ar[r] \ar[d] & X_{U[\![t_d]\!]} \ar[d] \\
X_{R_d} \ar[r] & X }
$$

For each maximal ideal $\mathfrak{q}$ of $U$, we apply the natural flat base change $U[\![t_d]\!] \to \hat{U}_{\mathfrak{q}}[\![t_d]\!]$ to the exact sequence \eqref{eq-long-seqence} and obtain $\bar \alpha_{b, \mathfrak{q}} = \hat{U}_{\mathfrak{q}}[\![t_d]\!] \otimes_{U[\![t_d]\!] } \bar \alpha_{b}$. We need to show that $\bar \alpha_{b, \mathfrak{q}}$ is an isomorphism. Consider further a base change map $\hat{U}_{\mathfrak{q}} \rightarrow V_{\mathfrak{q}} \coloneqq  W(k_{\mathfrak{q}})[\![s_1, \dots, s_m]\!]$ obtained from $f\colon \cO_{\mathfrak{q}} \rightarrow W(k_{\mathfrak{q}})$ as in Lemma \ref{lem-structrue} applied to \eqref{eq-long-seqence}. Since $\hat{U}_{\mathfrak{q}} \rightarrow V_{\mathfrak{q}}$ is faithfully flat, $\bar \alpha _{b, \mathfrak{q}}$ is an isomorphism if and only if $$ V_{\mathfrak{q}}[\![t_d]\!] \otimes_{\hat{U}_{\mathfrak{q}}[\![t_d]\!]} \bar \alpha _{b, \mathfrak{q}}\colon \rH ^i_\cris(X/V_{\mathfrak{q}, n}[\![t_d]\!], J^{[j]})/ t_d \inj \rH ^i_\cris (X/V_{\mathfrak q, n}, J^{[j]})$$
is an isomorphism. The above map is an isomorphism by inductive hypothesis, which proves Lemma \ref{lem:equals-zero}. 
\end{proof}

\section{Construction of Relative Fontaine--Laffaille data and the comparison map}

As before, let $\Xbar \coloneqq X \times_R \Rbar$. The aim of this section is to construct a morphism $ T_\cris (\rH^i_\cris (X_n / R_n) ) \to \rH^i_\et (X_{\Rbar [\frac 1 p]}, \Z/ p^n \Z)$ which is compatible with $G_R$-actions. Along the way, we also complete the proof of Theorem \ref{thm-cons-phi} (1) (2); we show in particular the existence of a Gauss--Manin connection $\nabla$ and $\varphi_i$ satisfying the conditions for a relative Fontaine--Laffaille module. Following Fontaine--Messing, we use syntomic cohomology $\rH^i _\syn (\Xbar)$ to relate $A_\cris (R)  \otimes_R \rH^i_\cris (X_n / R_n )$ and $\rH^i_\et (X_{\Rbar [\frac 1 p]}, \Z/ p^n \Z)$. The computation of syntomic (-crystalline) cohomology in the relative setting is complicated due to that our base ring $R$ is not perfect (as opposed to the classical case over $W(k)$). We circumvent this issue by systematically making base change from $R$ to an intermediate ring $\wt R$ (to be defined below) and descending back. We sometimes abuse notation and denote by $A$ the affine scheme $\Spec(A)$ if no confusion arises. For example, we sometimes denote $\rH ^i_\cris (\Spec (A)/ \Spec (R_n))$ by $\rH ^i_\cris ( A/ R_n)$. 

\subsection{The ring $\wt R_n$}\label{subsec-tildeRn} 

Recall $R_\infty \coloneqq \bigcup_{k \geq 0} R[(1+t_1)^{(k)}, \dots, (1+ t_d)^{(k)}] \subset \Rbar$ where $(1+ t_i)^{(k)}$ is a fixed choice of $p^k$-th root of $(1+t_i)$ satisfying $((1+t_i)^{(k+1)}) ^p = (1+ t_i)^{(k)}$. Its $p$-adic completion $\wh R_\infty$ is a subring inside $\wh \Rbar$, and $\wh R_\infty^\flat \subset \Rbar^\flat$. We embed $R$ into $W(\wh R^\flat_\infty)$ via $t_i\mapsto [\wt{1+ t_i}]-1$, which is compatible with the earlier choice of embedding $\lambda_1\colon R \hookrightarrow W(\Rbar^\flat)$.  Recall that for any $p$-adically complete and separated ring $A$, we have a functorial morphism $\theta\colon W(A^\flat) \to A$ as in \cite[Sec. 3.1]{SBM}. In particular, we have a commutative diagram:

 $$\xymatrix{ R\ar[d]^{\mathrm{id}} \ar[r] & W(\wh R_\infty^\flat)\ar[r] \ar[d]^{\theta}_{\wr}  &   W(\Rbar^\flat) \ar[d]^{\theta} \\ R\ar[r]  & \wh R_\infty \ar[r] & \wh \Rbar } $$

To see that the middle column is an isomorphism, note that the inclusion  $R \subset \wh R_\infty$ induces an isomorphism $\wt R_1 \simeq \wh R_\infty^\flat$ where $\wt R_1$ denotes the perfection of $R_1$.\footnote{Given a ring $A$ of characteristic $p$, its perfection means the ring $\displaystyle \tilde{A} \coloneqq \varinjlim_{\varphi\colon A \rightarrow A} A$. This is a perfect ring which is universal for maps into perfect rings.} Since $\wh R_\infty/(p) \simeq \wt R_1$, the middle column mod $p$ is an isomorphism, and hence the middle column is an isomorphism. We set $\wt R \coloneqq W(\wt R_1) \simeq W(\wh R_\infty^\flat)$ and $\wt R_n:=  W_n (\wt R_1)  = \wt R/ p ^n \wt R$. 

\begin{lemma}\label{lem-faithflat} With notations as above, the inclusion $R \inj \wt R$ is faithfully flat. 
\end{lemma}

\begin{proof}
We first show $\wt R $ is flat over $R$. Note that $\wt R/(p) \cong \wt R_1$ is flat over $R_1$. By \cite[\href{https://stacks.math.columbia.edu/tag/0AGW}{Tag 0AGW}]{stacks-project}, it suffices to show $\Tor_1^R (\wt R, R/pR) =0$, which follows from that $\wt R = W(\wt R_1)$ has no $p$-torsion. Moreover, for the maximal ideal $\fm$ of $R$, we have $\fm \wt R \not = \wt R$. Thus, $\wt R$ is faithfully flat over $R$. 
\end{proof}

\begin{rem}
From the commutative diagram above, we obtain embeddings $\wt R \hookrightarrow \wh\Rbar$ and $\wt R \hookrightarrow W(\Rbar^\flat)$.
\end{rem}

\subsection{Comparison of crystalline and Syntomic cohomology}

In this subsection, we interpret $\rH^i_\cris (X_n / R_n)$ and $\rH^i _\cris (\Xbar_n / R_n)$ with their filtrations as certain syntomic cohomology groups with filtrations. Given a scheme $Y$ over $S$, we denote by $\CRIS(Y /S )_\SYN$ (resp. $\SYN(Y)$) the big crystalline-syntomic site (resp. big syntomic site). Let $(Y/S )_{\CRIS-\SYN}$ and $(Y)_{\SYN}$ denote the corresponding topoi. Abusing notation, we will denote by $\cO_{Y/S}$ (resp. $J^{[r]}_{Y/S}$) the corresponding structure sheaf (resp. divided power ideals) in the crystalline-syntomic site. We refer to Appendix \ref{appendix:syn-cris} for a recollection of definitions and general facts on the syntomic and crystalline-syntomic sites and topoi.

We consider the sites $\CRIS(Y_n /S_n )_\SYN$ in the following two settings:
\begin{enumerate}
\item[(1)] We set $Y =X$ and $S_n = \Spec(R_n)$ with the standard divided power structure $(R_n, p R_n, \gamma_j (-))$,
\item[(2)] We set $Y = \Xbar$. In this case, the base $S_n$ is either $ \Spec(W_n(k)), \Spec(\wt R_n)$, or $\Spec(R_n)$.
\end{enumerate}
We have a natural morphism of topoi $u_{Y_n / S_n }\colon (Y_n /S_n )_{\CRIS-\SYN}\to (Y_n /S_n )_{\SYN}$; set $u \coloneqq u_{Y_n / S_n }$ to simplify the notation. We begin by showing that the higher direct images under $u$ of the sheaves $J^{[r]}_{Y_n / S_n}$ vanish in any of the settings above. In the case $Y=X$, $S=\Spec(W(k))$ and $r = 0$, the result is proved in \cite[Prop. II.1.3]{FontaineMessing}. The proof below follows closely that of \cite[Prop. 1.17]{Bauer}, which generalizes the result of Fontaine--Messing to coefficients in quasi-coherent crystals. We note that, while the sheaves $J^{[r]}_{Y_n / S_n}$ are not crystals, the vanishing of higher direct images still holds in this case with essentially the same proof.

\begin{theorem}\label{thm-u-vanish}
With notations as above, $\bbR^{i}u_{*}(J^{[r]}_{Y_n / S_n}) = 0$ for all $i > 0$ and $r \geq  0$.
\end{theorem}

\begin{proof} 
It suffices to show that the sheaf associated to the pre-sheaf
$$U \rightarrow F(U) \coloneqq \rH^{i}((Y_n/S_n)_{\CRIS-\SYN}/\tilde{U}, J^{[r]}_{Y_n/ S_n}), \  i > 0$$
	vanishes in the syntomic topology on $Y_n$ for all $ i >0$. Here $\tilde{U}$ is the pull back of the sheaf represented by $U$. On the other hand, $F(U)$ can be identified with $\rH^{i}_{\cris}(U/S_n, J^{[r]}_{U/ S_n})$ (see Remark \ref{rem:mor-crissyn-syn} and  Corollary \ref{cor-alpha-vanish}). Since this is a local statement, we may assume $U = \Spec(A)$ where $\Spec(A)$ is flat over $Y_n$. As in the proof of \cite[Prop. 1.17]{Bauer}, it suffices to show there exists a family of faithfully flat syntomic morphisms $\{V_i \rightarrow U\}_{i \in I}$ such that $\varinjlim_{i \in I} F(V_i) = 0$. We will show below that the same construction as in \textit{loc. cit.} also works in our case. Namely, let $A_0 =A$ and $A_{i+1} = A_i[(T_a)_{a \in A_i}]/((T_a^p -a)_{a \in A_i})$. Each $A_i$ is a faithfully flat syntomic cover of $A$, and we are reduced to showing that $ \varinjlim_{i \in I} F(A_i) = F(\hat A) = 0$, where $\hat A \coloneqq \varinjlim A_i$. The Frobenius on $\hat A/ p \hat A$ is surjective by construction. Now the theorem follows from the explicit calculation of $\rH^i_\cris (\hat A/S_n, J^{[r]}_{\hat A/S_n})$ below in Proposition \ref{prop-calculateH0-1} and \ref{prop-calculateH0}.
\end{proof}

To finish the proof of Theorem \ref{thm-u-vanish}, it suffices to show $\rH^i_\cris (A/S_n, J^{[r]}_{A/S_n}) = 0$ for all $i>0$ and $S_n$-algebras $A$ such that the Frobenius on $A/pA$ is surjective. We first prove this vanishing in the cases when $S_n = \Spec(W_n(k))$ or $\Spec(\wt R_n)$. In this setting, the base ring is perfect and the construction of a final object in the corresponding crystalline site is simple. In the case of $S_n = \Spec(W_n(k))$, this is due to Tamme (cf. \cite[Prop. 1.16]{Bauer}). On the other hand, the case of $S_n = \Spec(R_n)$ is slightly more involved, and is studied separately in Proposition \ref{prop-calculateH0}.

\begin{proposition} \label{prop-calculateH0-1}
Let $S_n = \Spec(W_n(k))$ or $\Spec(\wt R_n)$. Let $\mathrm{Spec}(A)$ be an $S_n$-affine scheme such that the Frobenius $\varphi$ on $A/pA = A_1$ is surjective. Then $\rH^i_{\cris}(A/S_n, J^{[r]}_{A/S_n}) = 0$ for all $i > 0$. 
\end{proposition}

\begin{proof} 
Abusing notation, we also denote by $S_n$ the underlying ring of the affine scheme $S_n$. To compute the relevant crystalline cohomology groups, we will construct a final object for $\CRIS(A/ S_n)$. More precisely, we shall work with the subcategory of affine objects $(\Spec B, \Spec C, i, \delta )$, where $B$ is an $A$-algebra, $i\colon \Spec B \inj \Spec C $ is a closed immersion, and $\delta $ is a divided power structure on the ideal of the definition of $i$ which is compatible with that on $(S_n,  p S_n)$ (see \S \ref{appendix-crystalline-site} for general facts on crystalline site).   
Note that this is sufficient for computing the relevant crystalline cohomology. Let $B$ be an $A$-algebra and $C$ be an $S_n$-linear PD-thickening of $B$. We first construct a morphism  
	$\theta_{C, n}\colon W_n (A_1^\flat) \to C$ of $ S_n$-algebras as follows. Given $x = (x_0, x_1, \dots, x_{n -1}) \in W_n (A_1^\flat)$ with $x_i = (x_i ^{(m)})_{m \geq 0} \in A_1^\flat$ and $x^{(m)}_i \in A_1$, let $y _i\in C$ be a lift of $x_i ^{(n)} \in A_1$. Set $$ \theta_{C, n}(x) \coloneqq \sum_{i = 0}^{n -1} p ^i y_i^{p ^{n-i}}. $$
Since $C$ is a PD-thickening of $B$, it is standard to check that $\theta_{C, n}$ does not depend on the choice of the lift $y_i$ of $x_i^{(n)}$. We claim that $\theta_{C, n}$ is a ring homomorphism. Consider the diagram 
	$$\xymatrix{W_n (A_1^\flat) \ar[r]^\alpha &  W_n (A_1) \ar[r] ^{\alpha'} &  W_n (B_1) &  W_n (C)\ar[r]^-\beta \ar@{->>}[l]_-{\gamma} &  C}  $$
where $\beta$ is defined by the ghost component of $W_n(C)$ and $\alpha$ by sending $x$ to $(x^{(n)}_0, \dots , x_{n-1}^{(n)})\in W_n (A_1)$. It is clear that  $\alpha$ is a ring homomorphism. $\theta_{C, n}$ is the ``composite" of the above ring homomorphisms (for $x \in W_n (A_1^\flat)$, if we choose any $ y \in W_n (C)$ such that $\gamma (y) = (\alpha'\circ \alpha) (x)$, then $\theta_{C,n} = \beta (y)$), and thus a ring homomorphism.  
	
From the construction, we obtain the following commutative diagram: 
{	\begin{equation}\label{diag-theta}
	\xymatrix{W_n (A_1^\flat) \ar[r]^-{\theta_{C, n}}\ar[d]^{\theta_{A, n}} & C \ar@{->>}[d]^ q \\ A \ar[r]^f & B }
	\end{equation}}
We claim that $\theta_{C, n}$ is the unique map making the above diagram commute. More precisely, if there exists another ring homomorphism $g\colon W_n (A_1^\flat) \to C$ so that $q \circ g = f \circ \theta_{A, n}$, then $\theta_{C, n} = g$. To prove this, note that $x= \sum p^i [\varphi^{-i}(x_i)]$ for any $x = (x_0, \dots, x_{n-1}) \in W_n(A_1^\flat)$. So it suffices to show $g ([x_0])= \theta_{C, n} ([x_0])= y_0^{p ^n}$ where $y_0\in C$ is a lift of $f( x^{(n)}_0) \mod p \in B_1$. Let $z= \varphi ^{-n} (x_0)$. Then $g([z])$ is a lift of $f (\theta_{A, n}([z])) \mod p \in B_1$. Since $f(\theta_{A, n} ([z]))= f(x_0^{(n)}) \mod p$, both $g([z])$ and $y_0$ are lifts of $f(x_0^{(n)}) \mod p$. Hence, $g([x_0])= (g([z]))^{p^n}= y _0 ^{p ^n}= \theta_{C, n}([x_0])$. 
	
Next, we claim that the above diagram is compatible with $S_n$-algebra structures. Since $A_1$ is an $S_1$-algebra, we have a natural map $S^\flat_1 \to A^\flat_1$. The map $\theta_0\colon S_1^\flat \to S_1$ given by $(x^{(n)})_{n \geq 0} \mapsto x^{(0)}$ is an isomorphism since $S_1$ is perfect. We have $S_n = W_n(S_1)$, so the morphism $\theta_0$ gives $W_n (A_1^\flat)$ an $S_n$-algebra structure. To see that $\theta_{C, n}$ is an $S_n$-algebra map, let $x = (x_0 , \dots , x_{n -1})\in S_n = W_n(S_1)$ with $x_i \in S_1$. Then $x = \sum^{n -1} _{i=0} p ^i[\varphi^{-i} (x_i)]$. To show $\theta_{C, n} (x) = x$, it thus suffices to consider the case when $x= {[x_0]}$, which is clear from the construction. 	
	
We now show $\theta_{A, n}: W_n (A_1 ^\flat) \to A$ is surjective by induction on $n$. Since $A_1^\flat$ is perfect, we have $W_n (A_1^{\flat})/ p^{n-1} W_n (A_1^{\flat})\simeq W_{n-1}(A_1^\flat).$ Moreover, from the construction, $\theta_{A, n}\mod p ^{n-1}= \theta_{A, n-1}$. Note that $\theta_{A, n}$ when $n=1$ is surjective since $\varphi\colon A_1 \to A_1$ is surjective. Suppose $\theta_{A, n-1}$ is surjective. Then for any $x\in A$, there exists $y \in W_n (A_1^\flat)$ such that $\theta_{A, n} (y)- x = p ^{n -1} z$ for some $z \in A$. Pick $w \in W_n (A_1^\flat) $ so that $\theta_{A, n} (w)= z \mod p$. Then $\theta_{A, n}(y - p^{n-1} w)= x$, so $\theta_{A, n}$ is surjective.  
	
Finally, let $W^\PD_n(A_1^\flat)$ be the divided power envelope of $W_n (A_1 ^\flat)$ with respect to the ideal $\ker \theta_{A, n}$. From the above discussion and the universal property of divided power envelopes, we see that $W_n ^\PD(A_1 ^\flat)$ is the desired final object. Thus, 
	$$\rH^i_\cris (A/ S_n, J^{[r]}_{A/ S_n })= \rH^i_{\mathrm zar} (\Spec(W^\PD_n (A_1^\flat)), J ^{[r]}_{W^\PD_n (A_1^\flat)})$$ 
 where $J^{[r]}_{W^\PD_n(A_1 ^\flat)}$ denotes the $r$-th PD ideal of ${W^\PD_n(A_1 ^\flat)}$. In particular, $\rH^i_\cris (A/ S_n, J^{[r]}_{A/ S_n })= 0 $ for $i > 0$. 
\end{proof}	

\begin{rem}\label{rem-calculateH0-1}
From the above proof, we also obtain $$\rH^0_\cris(A/ S_n, J^{[r]}_{A/ S_n}) \cong J^{[r]}_{W^\PD_n(A_1 ^\flat)}.$$
\end{rem}
 
\begin{rem}
The proof above uses the assumption that $S_n$ is perfect to define a natural $S_n$-algebra structure on $W^\PD_n(A_1^\flat)$ and to show that $\theta_{C, n}$ is an $S_n$-algebra map.
\end{rem}

We now prove an analog of Proposition \ref{prop-calculateH0-1} in the setting where $S_n = \Spec(R_n)$. We begin by constructing a final object in the full subcategory of affine objects of $\CRIS(A/R_n)$. Let $W_n (A_1^\flat)_R \coloneqq R \otimes_{W(k)} W_n (A_1^\flat)$. Suppose $(\Spec(B), \Spec(C))$ is an object in  $\CRIS (A/R_n)$. We can consider $(\Spec(B), \Spec(C))$ as an object in $\CRIS(A/W_n(k))$. Then by the proof of Proposition \ref{prop-calculateH0-1}, we have a unique map $\theta_{C, n}\colon W_n (A_1 ^\flat)\to C$ such that diagram \eqref{diag-theta} commutes, and both $\theta_{A, n}$ and $\theta_{C, n}$ are compatible with $W_n(k)$-algebra structures. Thus, $\theta_{A, n}$ and $\theta_{C, n}$ naturally extend to maps $\theta_{A, n}\colon W_n (A_1^\flat)_R \to A$ and $\theta_{C, n}\colon W_n (A_1^\flat)_R \to C$. The analog of diagram \eqref{diag-theta} obtained from replacing $W_n (A_1^\flat)$ by $W_n (A_1^\flat)_R$ still commutes, and all such morphisms are compatible with $R_n$-algebra structures. It remains to show $\theta_{C, n}$ is unique in the sense that if there exists another $R_n$-algebra morphism $g\colon W_n (A_1^\flat)_R \to C$  satisfying $ q \circ g = f \circ \theta_{A, n}$ then $g = \theta_{C, n}$. If we restrict the diagram to $W_n(A_1^{\flat})$, then we get back diagram \eqref{diag-theta}. Thus, $g$ restricted to $W_n(A_1^\flat)$ is $\theta_{C, n}$ by the uniqueness. Since $g$ is $R$-linear, we have $g = \theta_{C, n}$. Let $\mathcal{O}A_{\mathrm{cris}, n}(A)$ be the divided power envelope with respect to the ideal $\ker (\theta_{A, n})$ in $W_n(A_1^\flat)_R$. Then $(\Spec A, \Spec(\mathcal{O}A_{\mathrm{cris}, n}(A)))$ is the desired final object, and we conclude the following:

\begin{proposition} \label{prop-calculateH0}  Suppose $A$ is an $ R_n$-algebra such that 
Frobenius $\varphi$ on $A/pA = A_1$ is surjective. Then  for each $i >0$ and $r \geq 0$, we have $\rH ^i_\cris (A / R_n , J_{A/ R_n }^{[r]}) = 0$. Furthermore, $$\rH ^0_\cris (A / R_n , J_{A/ R_n }^{[r]}) \cong \rF^r\mathcal{O}A_{\mathrm{cris}, n}(A)$$
where $\rF^r\mathcal{O}A_{\mathrm{cris}, n}(A)$ denotes the $r$-th PD-ideal of $\mathcal{O}A_{\mathrm{cris}, n}(A)$. 
\end{proposition}

 \begin{rems} A similar proof as above shows the following general statement: Let $(S, I, \gamma)$ be a PD $\Z_p$-algebra compatible with the standard PD-structure of $(\Z_p, p \Z_p, \gamma)$, and let $A$ be an $S$-algebra such that $\varphi$ on $A_1$ is surjective. Then the crystalline site $(A_n/S_n)_{\cris}$ has an affine final object $\Spec (B)$, where $B$ is the divided power envelope of $W_n (A_1^\flat)\otimes_{\Z_p}S$ with respect to the kernel of $W_n (A_1^\flat)\otimes_{\Z_p}S \to A$. We thank an anonymous referee for pointing this out. 
\end{rems}

By specializing to the case $A = \Rbar$ and comparing the constructions of $\cO A_{\cris} (R)$ and $\rF^r\cO A_\cris (R)$ in \S \ref{subsec-period-ring}, we obtain:
\begin{corollary}\label{cor-acris}
	\begin{enumerate}
		\item[(1)] For $r \geq 0$, we have canonical isomorphisms $$\rH^0_\cris(\Rbar_n / R_n, J ^{[r]}_{\Rbar_n / R_n}) \cong \rF^r\cO A_\cris (R)/ (p^n)$$
		where $\rF^r\cO A_\cris (R)$ denotes the $r$-th PD-ideal of $\cO A_\cris (R)$.
		\item[(2)] For all $ i > 0$ and $r \geq 0$,  we have $\rH^i_\cris(\Rbar_n / R_n, J ^{[r]}_{\Rbar_n / R_n}) = 0$.
	\end{enumerate}
\end{corollary}

\begin{rem} \label{rem:frobenius-divided-powers}
In the situations of Proposition \ref{prop-calculateH0-1} and \ref{prop-calculateH0}, $W^\PD_n(A_1^{\flat})$ and $\mathcal{O}A_{\mathrm{cris}, n}(A)$ have natural Frobenius endomorphisms induced by the Witt-vector Frobenius and $\varphi_R$. To see this, note that $W^\PD_n(A_1^{\flat})$ remains the same if we take the PD-envelope of $W_n (A_1^\flat)$ by the ideal $J'\coloneqq (p) + \ker (\theta_{A, n})$. Since $\varphi(J') \subset J'$ under $\varphi\colon W_n (A_1^\flat) \to W_n (A_1^\flat)$, $\varphi$ naturally extends to $W_n^\PD (A_1^\flat)$. A similar argument holds for $\mathcal{O}A_{\mathrm{cris}, n}(A)$.
\end{rem}

Consider an $\wt R_n$-algebra $A$ such that $\varphi$ on $A_1$ is surjective. We make some observations on the structure of $\mathcal{O}A_{\mathrm{cris}, n}(A)$ which will be useful later. We first relate it to the ring $A_{\mathrm{cris},n}(A) \coloneqq W^\PD_n(A_1^{\flat})$ and prove an analog of Lemma \ref{lem-2}. For each $i = 1, \ldots, d$, recall that we chose $\widetilde{1+t_i} \in \wt R_1 \cong \wt R_1^{\flat}$ satisfying $\widetilde{1+t_i}^{(0)} = 1+t_i$. Denote also by $\widetilde{1+t_i}$ the corresponding image under the induced map $\wt R_1^{\flat} \rightarrow A_1^{\flat}$. Let $A_{\mathrm{cris, n}}(A)\langle X_1, \ldots, X_d\rangle$ be the divided power polynomial ring in the variables $X_i$ with coefficients in $A_{\mathrm{cris}, n}(A)$, and let $\theta_{A, n}\colon A_{\mathrm{cris, n}}(A)\langle X_1, \ldots, X_d\rangle \rightarrow A$ be the map extending that on $A_{\mathrm{cris, n}}(A)$ given by $X_i \mapsto 0$. Consider the $A_{\mathrm{cris}, n}(A)$-linear map $$f\colon A_{\mathrm{cris, n}}(A)\langle X_1, \ldots, X_d\rangle \rightarrow \mathcal{O}A_{\mathrm{cris}, n}(A)$$ given by $X_i \mapsto (1+t_i)\otimes 1 - 1\otimes [\widetilde{1+t_i}]$. Note that $f$ is compatible with $\theta_{A, n}$.

\begin{lemma}\label{lem-toOAcris}
The map $f$ as above is an isomorphism.	
\end{lemma}
	
\begin{proof}
Consider the map $\bar{g}\colon R_1 = k[\![t_1, \ldots, t_d]\!] \rightarrow A_1^{\flat}$ given by $t_i \mapsto \widetilde{1+t_i}-1$. This indeed gives a ring homomorphism since for any $\displaystyle \sum_I a_I \underline{t}^I \in R_1$ with $a_I \in k$ and $I = (i_1, \ldots, i_d)$ a multi-index, we have
\[
\bar{g}(\sum_I a_I \underline{t}^I) = \sum_I a_I \prod_{j=1}^d(\widetilde{1+t_j}-1)^{i_j} \in A_1^{\flat}
\]	
with $\displaystyle (\bar{g}(\sum_I a_I \underline{t}^I))^{(0)} =  \sum_I a_I \underline{t}^I \in A_1$. By Cartier's Dieudonn\'e--Dwork lemma, there exists a unique lifting $g\colon R \rightarrow W(A_1^{\flat})$ of $\bar{g}$ compatible with Frobenius. Note that $\varphi$ on $W(A_1^{\flat})$ is an isomorphism since $A_1^{\flat}$ is perfect. Since $g$ is compatible with Frobenius, $g$ is given by $t_i \mapsto [\widetilde{1+t_i}]-1$. Indeed, suppose $g(t_i) = [\widetilde{1+t_i}]-1 + p^m z$ with $z \in W(A_1^{\flat})$ for some $m \geq 1$. Since $\varphi (g(t_i)) = g (\varphi (t_i))$, we deduce $\varphi(z) \in p W( A_1^\flat)$. So $z \in p W(A_1^\flat)$, which implies $z = 0$. Denote by $g_n\colon R \rightarrow W_n(A_1^{\flat})$ the induced map mod $p^n$.

Now, let $h\colon R \rightarrow A_{\mathrm{cris, n}}(A)\langle X_1, \ldots, X_d\rangle$ be the $W(k)$-linear map given by $t_i \mapsto X_i+[\widetilde{1+t_i}]-1$. This gives a ring homomorphism since $g_n$ as above is a ring homomorphism and $\underline{X}^I = 0$ for $|I| \gg 0$. Denote also by $h\colon R\otimes_{W(k)} W_n(A_1^{\flat}) \rightarrow A_{\mathrm{cris, n}}(A)\langle X_1, \ldots, X_d\rangle$ the $W_n(A_1^{\flat})$-linear extension. Since $h$ commutes with $\theta_{A, n}$, it extends uniquely to a map $h\colon \mathcal{O}A_{\mathrm{cris}, n}(A) \rightarrow A_{\mathrm{cris, n}}(A)\langle X_1, \ldots, X_d\rangle$. Note that $(h \circ f)(X_i) = X_i$ and $(f \circ h)(t_i) = t_i$. Since $f$ and $h$ are compatible with $\theta_{A, n}$, this implies that $h \circ f$ is the identity map on the PD-envelope $A_{\mathrm{cris, n}}(A)\langle X_1, \ldots, X_d\rangle$, and $f \circ h$ is the identity map on the PD-envelope $\mathcal{O}A_{\mathrm{cris}, n}(A)$.
\end{proof}

\begin{rem}\label{rem-GroupAction}
Suppose $A$ and $B$ are $\Rbar_n$-algebras such that $\varphi$ on $A_1$ and $B_1$ are surjective. Let $g \in G_R$, and suppose we have a $g$-semi-linear morphism $A \rightarrow B$ of $\Rbar_n$-algebras. By functoriality, it induces a $W_n(k)$-algebra map $A_{\mathrm{cris},n}(A) \rightarrow A_{\mathrm{cris},n}(B)$ and an $R_n$-algebra map $\cO A_{\mathrm{cris},n}(A) \rightarrow \cO A_{\mathrm{cris},n}(B)$, both of which are compatible with Frobenius. Under the isomorphisms $\mathcal{O}A_{\mathrm{cris}, n}(A) \cong A_{\mathrm{cris, n}}(A)\langle X_1, \ldots, X_d\rangle$ and $\mathcal{O}A_{\mathrm{cris}, n}(B) \cong A_{\mathrm{cris, n}}(B)\langle X_1, \ldots, X_d\rangle$ of the previous lemma, the induced map $\cO A_{\mathrm{cris},n}(A) \rightarrow \cO A_{\mathrm{cris},n}(B)$ sends $\gamma_j(X_i)$ to $g(\gamma_j(X_i)) = \gamma_j(X_i+[\widetilde{1+t_i}]-g([\widetilde{1+t_i}]))$, and is therefore compatible with the $A_{\mathrm{cris, n}}(A)$-linear (resp. $A_{\mathrm{cris, n}}(B)$-linear) connection given by $\nabla(X_i) = dt_i$.     	
\end{rem}

\subsection{Fontaine--Laffialle data associated to $(\mathcal{O}A_{\mathrm{cris}, n}(A),\rF^r\mathcal{O}A_{\mathrm{cris}, n}(A))$}

Let $f_1, \ldots, f_c$ be elements in $\Rbar[x_1, \dots, x_m]$ whose images in $\Rbar_n[x_1, \dots, x_m]$ (which, by abusing notation, we also denote by $f_1, \ldots, f_c${)} generate a Koszul-regular ideal, and suppose that $\Rbar_n[x_1, \ldots, x_m]/(f_1, \ldots, f_c)$ is flat over $\Rbar_n$. Let $A= \Rbar_n[x_1^{\frac{1}{p ^\infty}}, \ldots, x_m^{\frac{1}{p ^\infty}}]/ (f_1, \ldots, f_c)$ and note that Frobenius on $A_1$ is surjective. In this section, we show that $(\mathcal{O}A_{\mathrm{cris}, n}(A),\rF^r\mathcal{O}A_{\mathrm{cris}, n}(A))$ can be equipped with the structure of a relative Fontaine--Laffaille module functorially in $A$. In particular, we construct divided power Frobenii $\varphi_r$ and a connection $\nabla$. We begin by making some standard observations on the filtration $\rF^r\mathcal{O}A_{\mathrm{cris}, n}(A)$. 

\begin{lemma}\label{lem-info-fil}
With notations as above,
\begin{enumerate} 
\item[(1)] For any $r \geq 0$ and $1\leq m < n$, we have an exact sequence 
	$$0 \longrightarrow  J^{[r]}_{W^{\PD}_{n-m} (A_1^\flat)} \overset{\times p^{m}}{\longrightarrow}  J^{[r]}_{W^{\PD}_n (A_1^\flat)} \overset{\mod p^{m}}{\longrightarrow}   J^{[r]}_{W^{\PD}_{m} (A_1^\flat)} \longrightarrow 0 .$$

\item[(2)] We have an analogous exact sequence in the setting where $\Rbar$ is replaced by $\wt R$, and with $A$ constructed in an analogous manner (i.e. with $f_1, \ldots, f_c \in \wt R[x_1, \ldots, x_m]$ such that $(f_1, \ldots, f_c) \subset \wt R_n[x_1, \ldots, x_m]$ is Koszul-regular and $\wt R_n[x_1, \ldots, x_m] / (f_1, \ldots, f_c)$ is flat over $\wt R_n$, and $A = \wt R_n[x_1^{\frac{1}{p ^\infty}}, \ldots, x_m^{\frac{1}{p ^\infty}}] / (f_1, \ldots, f_c)$).

\end{enumerate}	
\end{lemma}

\begin{proof}
Consider $C = \Rbar [x_1^{\frac{1}{p ^\infty}} , \dots, x_m^{\frac{1}{p ^\infty}} ]$ with the ideal $I = (f_1, \dots, f_c)$. By \cite[Lem. 3.10]{SBM}, the $p$-adic completion of $C$ is a perfectoid algebra, and the kernel of $\theta_{C, n}\colon W_n (C_1^\flat) \rightarrow C_n$ is generated by $\xi$. Let $\theta'$ be the composite $\theta'\colon W_n(C_1^\flat) \stackrel{\theta_{C, n}}{\longrightarrow} C_n \rightarrow A$, which is surjective. Denote $\xi^{\flat} \coloneqq \xi \mod p \in C_1^{\flat}$. Let $\wt f_i \in W_n (C_1^\flat)$ be a lift of $f_i$ along $\theta_{C, n}$, and denote $f_i^\flat \coloneqq \wt f_i \mod p \in C_1^\flat$ and $\bar{f}_i \coloneqq f_i \mod p$. Then $\theta_{C, 1}(f_i^\flat) = \bar f_i$. Since $\xi$ is a nonzero divisor in $W_n(C_1^\flat )$, by \cite[\href{https://stacks.math.columbia.edu/tag/0669}{Tag 0669}]{stacks-project}, the sequence $\xi , \wt f_1 , \dots , \wt f_c$ is Koszul-regular in $W_n(C_1^\flat)$. Set $I^\flat= (\xi^\flat, f_j^\flat) C_1^\flat$ and $\wt I = (\xi, \wt f_j) W_n (C_1^\flat)$. Note that the surjection $C_n \onto A$ induces a map $\iota\colon C_1^\flat \to A_1^\flat$. 

We first show that $A_1^\flat$ is the $I^\flat$-completion $\widehat C_1^\flat$ of $C_1^\flat$. We construct a map $ \bar \iota\colon \widehat C_1^\flat \to A_1^\flat$ compatible with  $\iota$ as follows. Suppose for each $m \geq 1$, we are given $y_m = (y^{(l)}_m)_{l \geq 0} \in C_1^\flat$ such that $y_{m+1} - y_m \in (I^\flat)^{(c+1)p^m}$. Then $y^{(l)}_{m +1} - y_m^{(l)} \in \bar I \coloneqq (\bar f_1 , \dots , \bar f_c) C_1$ for $l\leq m$. So as $m \to \infty $, for all $l \geq 0$, $y^{(l)}_{m}$ determines an element $\bar y^{(l)}\in A_1 = C_1/ \bar I$ satisfying $(\bar y^{(l+1)})^p = \bar y^{(l)}$. This defines a ring map $\bar \iota\colon \widehat C_1^\flat\to A_1^\flat$ via $\{y_m\} \mapsto (\bar y^{(l)})_{l \geq 0}$. By construction, $\bar \iota$ is compatible with $\iota\colon C_1^\flat \to A_1^\flat$. To check $\bar \iota$ is injective, note that if $\bar y^{(l)}= 0 \in A_1$ then $y_m^{(l)} \in \bar I C_1$ for $m$ sufficiently large. Since $C_1^\flat$ is perfect, there exists $x \in (f^\flat_1, \dots, f_c^\flat)^l C_1^\flat$ such that $x^{(l)} = y _m ^{(l)}$. Since $\varphi\colon C_1/ (p_1) \rightarrow C_1$ is an isomorphism, we have $x-y_m \in (p^\flat)^l C_1^\flat$, and thus $y_m \in (I^\flat) ^l$. This proves $\bar \iota$ is injective. To show that $\bar \iota$ is surjective, let $y = (y^{(l)}) _{l \geq 0} \in A_1^\flat$. For each $l$, let $x_l \in C_1^\flat$ such that the image of $x^{(l)}_l$ in $A_1$ is $y^{(l)}$. It suffices to check that $x_{l+1}-x_l \in (I ^\flat)^l$, since then $\{x_l\} \in \widehat C_1^\flat$ satisfies $\bar \iota (\{x_l\} ) = y$. Note that $ x^{(l)}_{l+1} -x^{(l)}_l\in  \bar I$, so there exists $z \in  (f_1^\flat , \dots, f^\flat_c)^lC_1^\flat$ such that $x ^{(l)}_{l +1} = x^{(l)}_l + z^{(l)}$. Thus, $x_{l+1}-x_l -z \in (p ^\flat) ^l C_1^\flat $, i.e. $x_{l+1} -x_l \in (I^\flat)^l$. 

We now claim that for $m = l p^n $ with $l \geq n $, we have $\wt I ^m\subset W_n ((I^\flat) ^{m}) \subset \wt I^{l -n}$. Note that $\tilde f_j \equiv [f_j ^\flat] \mod p$ and $\xi \equiv [p ^\flat]\mod p$. Any $x \in \tilde I ^m$ can be written as $x= \sum_{i = 0}^{n -1} x_i p ^{i}$ where $x_i = \sum_j [y_{ij}]$ with $y_{ij} \in (I^\flat) ^{m-i}$. Since $p = \varphi V$ where $V$ is Verschiebung, we have $p^i [y _{ij}]\in W_n ((I^\flat) ^{p ^i(m -i)} )$. So $x \in W_n ((I^\flat)^m)$. On the other hand, any $x \in W_n((I^\flat)^m)$ can be written as $x = \sum_i p ^i [x_i]$ with $x_i \in (I^\flat)^ {m / p ^i}$. It follows by induction on $n$ that if $y \in (I^\flat) ^l$, then $[y] \in \wt I^{l-n}$. Thus, $x \in \wt I^{l-n}$.

Let $W_n(C_1^\flat)^\wedge$ denote the $\wt I$-completion of $W_n(C_1^\flat)$. From the above discussion, we have
\[
W_n(C_1^\flat)^\wedge \cong W_n(\widehat{C}_1^\flat) \cong W_n(A_1^\flat).
\] 
Hence, we obtain the commutative diagram:
$$\xymatrix{ W_n(C_1^\flat )^\wedge\ar[dr]_{\simeq} & \ar[l]W_n (C_1^\flat)\ar[dr]^{\theta'} \ar[r]^-{\theta_{C, n}}\ar[d]_{\iota} &  C_n \ar[d] ^{\mod (f_j)} \\  & W_n (A_1^\flat) \ar[r]^-{\theta_{A, n}} & A} $$
In particular, $\ker(\theta_{A, n}) = \wt I W_n(A_1^\flat)$. The natural map from the PD-envelope $D_{W_n(C_1^\flat)}(\wt I)$ of $\wt I$ in $W_n(C_1^{\flat})$ (over the PD-ring $(\Z/(p^n), (p))$) to $W_n^{\PD}(A_1^\flat)$ is an isomorphism of PD-algebras by \cite[Rem. 3.20 (7)]{crystal}. Thus, $J^{[r]}_{W^{\PD}_{n} (A_1^\flat)}$ for any $r \geq 0$ is $\Z/p^n \Z$-flat by Lemma \ref{lem-Jflat} below, and statement (1) follows.  

For (2), the above arguments go through after following modifications. Let $C = \wt R [x_1^{\frac{1}{p ^\infty}}, \dots, x_m^{\frac{1}{p ^\infty}} ]$ with the ideal $I = (f_1, \dots f_c)$. Note in this case that $\theta_{C, n}\colon  W_n(C_1^\flat) \rightarrow C_n$ is an isomorphism since $C_1$ is perfect. Let $\wt I = (\wt f_1, \ldots, \wt f_c)W_n(C_1^\flat)$ and $I^\flat = (f_1^\flat, \ldots, f_c^\flat)C_1^\flat$. It can be checked easily that there is an isomorphism between $A_1^\flat$ and the $I^\flat$-completion of $C_1^\flat$ compatible with the natural map $C_1^\flat \rightarrow A_1^\flat$, and that this induces an isomorphism between $W_n(A_1^\flat)$ and the $\wt I$-completion of $W_n(C_1^\flat)$. Hence, $\ker(\theta_{A, n}) = \wt I W_n(A_1^\flat)$, and (2) follows similarly as above.
\end{proof}	

\begin{lemma}\label{lem-Jflat}
Let $B$ be a flat $\Z/p^n\Z$-algebra, $I \subset B$ a finitely generated Koszul-regular ideal, and suppose that $B/I$ is flat over $\Z/p^n\Z$. Then the PD-envelope $D_{B}(I)$ and its $r$-th divided power ideal $J^{[r]}$ for any $r \geq 1$ are $\Z/p^n\Z$-flat.
\end{lemma}

\begin{proof}
By \cite[\href{https://stacks.math.columbia.edu/tag/068Q}{Tag 068Q}]{stacks-project}, there exists a faithfully flat smooth morphism $\Spec(B') \rightarrow \Spec(B)$ such that 
$IB'$ is generated by a regular sequence. Since $D_B(I) \otimes_B B' \cong D_{B'}(IB')$ by \cite[Prop. 3.21]{crystal}, we can reduce to the case where $I \subset B$ is generated by a regular sequence $(f_1, \ldots, f_m)$.

Write $\Spec(B) = \varprojlim \Spec(B_{\alpha})$ as an inverse limit where $B_\alpha$'s are finitely generated flat $\Z/p^n\Z$-algebras, and similarly write $\Spec(B/I) = \varprojlim \Spec(B_{\alpha}/I_{\alpha})$ as a limit of finitely generated flat $\Z/p^n\Z$-algebras where $I_\alpha$ denotes the pull-back of $I$ to $B_{\alpha}$. By \cite[19.8.2]{EGAIV.4}, the map $\Spec(B/I) \hookrightarrow \Spec(B)$ is a regular immersion if and only if $\Spec(B_\alpha/I_\alpha) \hookrightarrow \Spec(B_\alpha)$ is a regular immersion for sufficiently large $\alpha$. Since divided power envelopes are compatible with direct limits, we may assume that $B$ is noetherian.

Consider the natural cartesian diagram
$$
\xymatrix{
\Spec(B/I) \ar[r] \ar[d] & \Spec(B) \ar[d]^{g}\\
\Spec(\Z/p^n\Z) \ar[r] & \bbA_{\Z/p^n\Z}^{m}}
$$
where the bottom horizontal map is given by $T_i \mapsto0$ and right vertical by $T_i \mapsto f_i$. Since the fiber over the origin is flat and the top horizontal map is a regular immersion, by \cite[19.2.4]{EGAIV.4}, we deduce that for every point $x \in \Spec(B/I)$, the map $g$ is flat in a neighborhood of $x$. Thus, there is an open neighborhood $U \subset \Spec(B)$ of $\Spec(B/I)$ such that the resulting diagram
$$
\xymatrix{
\Spec(B/I) \ar[r] \ar[d] & U \ar[d]\\
\Spec(\Z/p^n\Z) \ar[r] & \bbA_{\Z/p^n\Z}^{m}}
$$
is cartesian and the right vertical map is flat. We have isomorphisms of PD-algebras
\[
D_B(I) \cong D_U(I) \cong D_{\bbA_{\Z/p^n\Z}^{m}}((T_1, \ldots, T_d))\otimes_{\bbA_{\Z/p^n\Z}^{m}} U.
\] 
Since $D_{\bbA_{\Z/p^n\Z}^{m}}((T_1, \ldots, T_d))$ and its $r$-th divided power ideals are flat over $\Z/p^n\Z$, the lemma holds. 
\end{proof}

\begin{corollary}\label{cor-exact-fil} For any $r \geq 0$ and $1\leq m < n$, we have an exact sequence 
$$0 \longrightarrow  \rF^r \cO A_{\cris,n-m} (A)  \overset{\times p ^{m}}{\longrightarrow}  \rF^r \cO A_{\cris,n} (A) \overset{\mod p ^{m}}{\longrightarrow}   \rF^r \cO A_{\cris,m} (A) \longrightarrow 0 .$$
\end{corollary}

\begin{proof}
The isomorphism $f$ of Lemma \ref{lem-toOAcris} is compatible with filtrations given by the divided power ideals, since $f$ is compatible with $\theta_{A, n}$. In particular, we have 
$$\rF^r \cO A_{\cris, n} (A) = \bigoplus_{\sum i_j \geq r} \rF^{0} A_{\cris, n} (A) \gamma_{i_1} (X_1)\cdots \gamma_{i_d} (X_d) \bigoplus_{\ell = 1}^r\bigoplus_{\sum i_j = r-\ell} \mathrm{F}^{\ell} A_{\mathrm{cris}, n}(A) \gamma_{i_1} (X_1)\cdots \gamma_{i_d} (X_d)$$
where $\gamma_{i_j}$'s denote the corresponding divided power operations. Since $\mathrm{F}^{\ell} A_{\mathrm{cris}, n}(A) = J^{[\ell]}_{W_n^{\mathrm{PD}}(A_1^{\flat})}$ (see Remark \ref{rem-calculateH0-1}), the assertion follows from Lemma \ref{lem-info-fil}. 
\end{proof}

We now explain the construction of $\varphi_r$ and connection on $\cO A_{\cris,n}(A)$. Note that both $A_{\cris, n}(A)$ and $\cO A_{\cris, n}(A)$ have natural Frobenii induced by those on $W_n (A_1^\flat)$ and $R_n$ (see Remark \ref{rem:frobenius-divided-powers}). We have $\varphi (J^{[r]}_{W_n^\PD(A_1^\flat)}) \subset p ^r W_n ^\PD (A_1^\flat)$ for any $0 \leq r\leq  p-1$. Define $\varphi_r\colon \rF^r A_{\cris , n} (A) \to A_{\cris , n}(A)$ for $0 \leq r \leq p-1$ as follows. First, given $x \in \rF^r A_{\cris , n} (A)$, choose a lift $\hat x \in \rF^r A_{\cris , n+r} (A)$. We have $\varphi (\hat x) = p^r \hat y $ for some  $\hat y \in A_{\cris, n+r} (A)$. Define $\varphi_r (x) \coloneqq \hat y \mod p^r \in A_{\cris, n}(A)$. It follows from Lemma \ref{lem-info-fil} that $\varphi_r$ is well-defined. Similarly, we define $\varphi_r\colon \rF^r \cO A_{\cris, n}(A) \to \cO A_{\cris,n}(A)$ using Corollary \ref{cor-exact-fil}. 

Consider the natural connection $$\nabla\colon A_{\cris,n}(A)\langle X_1,\cdots,X_d \rangle \rightarrow A_{\cris,n}(A)\langle X_1,\cdots,X_d \rangle \otimes_{R} \wh \Omega_{R}$$
given by setting $\nabla(A_{\cris,n}(A)) = 0$ and $\nabla(X_i) = dt_i$. By Lemma \ref{lem-toOAcris}, this defines a natural connection $\nabla\colon \cO A_{\cris, n} ( A) \to \cO A_{\cris , n}(A) \otimes_R \wh \Omega_{R}$. The following corollary is immediate.

\begin{corollary}\label{cor-varphiandconenction}
Let $0 \leq r \leq p-1$.
\begin{enumerate} 
\item[(1)] The morphism $\varphi_r\colon \rF^r\cO A_{\cris , n} (A) \to \cO A_{\cris , n}(A) $ is $\varphi$-semilinear, and for $r \leq p-2$, $\varphi_r|_{\rF^{r+1}\cO A_{\cris , n} (A)} = p \varphi _{r+1} $. We have a similar statement for $A_{\cris, n} (A)$. 
\item[(2)] The connection $\nabla$ on $\cO A_{\cris,n}(A)$ is a flat connection such that:
\begin{enumerate}
 \item $\nabla ( \rF^r\cO A_{\cris , n} (A) ) \subset  \rF^{r-1}\cO A_{\cris, n}(A)\otimes_R \wh \Omega_{R}$.
 \item $\nabla\circ \varphi_r = (\varphi_{r-1}\otimes  d \varphi_1) \circ \nabla$ (when $r = 0$, $\varphi_{-1}$ is understood to be $p\varphi_0$). 
\end{enumerate}
\end{enumerate}
The above structures $(\cO A_{\cris , n}(A), \rF^r\cO A_{\cris , n} (A), \varphi_r, \nabla)$ are functorial in $A$.
\end{corollary}

\begin{rem}
Let $(f_1,\ldots,f_c) \subset \Rbar_n[x_1,\ldots,x_m]$ be a Koszul-regular sequence such that $\Rbar_n[x_1, \ldots, x_m]/(f_1, \ldots, f_c)$ is flat over $\Rbar_n$. We may choose lifts of $f_i$ in $\Rbar[x_1,\ldots,x_m]$, and Lemma \ref{lem-info-fil} and previous corollary hold for algebras $A$ as above. This observation will be useful below as we work with syntomic sites modulo $p^n$.
\end{rem}

\subsection{The divided power sheaves $\cJ^{[r]}_{Y_n/S_n}$ on the syntomic site} \label{subsec-divpowersheaf}

In this section, we make some observations on the divided power sheaves $\cJ^{[r]}_{Y_n/S_n}$ on the (small) syntomic site. Recall that $X$ is a smooth proper scheme over $R$. We work with $Y/S$ where $Y = \Xbar$ and $S$ is either $\Spec(R), \Spec(\wt R)$ or $\Spec(W(k))$, or where $Y = X$ and $S = \Spec(R)$.

Given integers $n > 0$ and $r \geq 0$, consider the following presheaf on the small syntomic site $(Y)_{\syn}$:
$$\cJ^{[r]}_{ Y_n/S_n}\colon  U \mapsto \rH^0_\cris (U/ S_n, J^{[r]}_{U_n/ S_n}), \ \  \forall  U \in (Y_n)_\syn. $$

\begin{lemma}\label{lem:divpowersheafcomp}
With notations as above, 
\begin{enumerate}
\item[(1)] \label{lem-itissheaf} $\cJ^{[r]}_{Y_n/S_n}$ is a sheaf on $(Y_n)_\syn$. 
\item[(2)] \label{lem-syn=cry} We have isomorphisms $$\rH ^i _\cris (Y_n /S_n , J ^{[r]}_{Y_n / S_n}) \simeq \rH ^i _{\syn-\cris}(Y_n/S_n , J ^{[r]}_{Y_n / S_n} ) \simeq \rH^i_\syn (Y_n, \cJ^{[r]}_{Y_n/S_n}). $$
\end{enumerate}
\end{lemma}

\begin{proof}
The second part is a consequence of the first part, Corollary \ref{cor-alpha-vanish}, and Theorem \ref{thm-u-vanish}. The first part follows from the fact that $\cJ^{[r]}_{ Y_n/S_n}$ is the push-forward of the analogous divided power sheaf on the corresponding syntomic-crystalline site.
\end{proof}

We now explain how to define $\varphi_r$ (for $0 \leq r \leq p-1$) and $\nabla$ on $\cJ ^{[r]}_{Y_n /S_n}$ when $S = \Spec(R)$. We first consider the case when $Y= \Xbar$. Let $\cC_{Y_n} \subset (Y_n)_{\syn}$ denote the full subcategory of affine schemes $\Spec(A)$ over $Y_n$ such that $A \simeq \Rbar_n[x_1,\ldots,x_m]/I$  where $I$ is a Koszul-regular ideal. Since $\Xbar$ is syntomic (in fact smooth) over $\Rbar$, to define a morphism of sheaves on $(Y_n)_{\syn}$, it suffices to construct such a morphism when restricted to this subcategory. 

Consider $\Spec(A_{(0)}) \in \cC_{Y_n}$ where $A_{(0)} \simeq \Rbar_n[x_1, \dots, x_m]/ (f_1, \dots f_c)$ such that $(f_1,\ldots,f_c)$ is Koszul-regular, and let $A_{(i)} \coloneqq \Rbar_n[x_1^{\frac{1}{p^i}}, \dots, x_m^{\frac{1}{p^i}}]/ (f_1, \ldots f_c)$. Let $A = \varinjlim_{i} A_{(i)} \cong \Rbar_n[x_1^{\frac{1}{p ^\infty}}, \dots, x_m^{\frac{1}{p ^\infty}}]/ (f_1 , \dots f_c)$. By Lemma \ref{lem:divpowersheafcomp}, we have an exact sequence
 $$0 \rightarrow \xymatrix {\cJ^{[r]}_{Y_n/S_n}(A_{(0)})  \ar[r] &  \cJ^{[r]}_{Y_n/S_n}(A)  \ar@<1ex>[r]\ar[r] & \cJ^{[r]}_{Y_n/S_n}(A \otimes _{A_{(0)}} A)}. $$
Note that $\cJ^{[r]}_{Y_n/S_n}(A) = \varinjlim_i \cJ^{[r]}_{Y_n/S_n}(A_{(i)}) = \rH^0_\cris (A/S_n, J^{[r]}_{A/S_n})$. By Corollary \ref{cor-varphiandconenction}, $\varphi_r$ and $\nabla$ are well-defined on both  $\cJ^{[r]}_{Y_n/S_n}(A) $ and $\cJ^{[r]}_{Y_n/S_n}(A \otimes _{A_{(0)}} A)$, since $A \otimes _{A_{(0)}} A \cong \Rbar_n[x_1^{\frac{1}{p^\infty}}, \ldots, x_m^{\frac {1} {/p^\infty}}, y_1^{\frac{1}{p^\infty}}, \ldots , y_m^{\frac {1}{p^\infty}}]/ (f_1(x), \ldots f_c(x), x_1-y_1, \ldots, x_m- y_m)$. By functoriality, $\varphi_r$ and $\nabla$ are compatible with the maps $\xymatrix {\cJ^{[r]}_{Y_n/S_n}(A)  \ar@<1ex>[r]\ar[r] & \cJ^{[r]}_{Y_n/S_n}(A \otimes _{A_{(0)}} A) }$. Thus, we have an induced Frobenius $\varphi_r\colon \cJ^{[r]}_{Y_n/S_n}(A_{(0)}) \to \cO_{Y_n /S_n} (A_{(0)})$ and a flat connection $\nabla \colon \cO_{Y_n /S_n} (A_{(0)}) \rightarrow \cO_{Y_n /S_n} (A_{(0)}) \otimes _R \wh \Omega_{R} $ resulting from the exact sequence above (as an application of Corollary \ref{cor-varphiandconenction}). Moreover, we may view $R$ and $\wh \Omega_R$ as  constant sheaves. The differential $d\colon R \rightarrow \wh \Omega_R$ can then be viewed as a morphism of constant sheaves. It follows that we have a morphism of sheaves
$\nabla \colon \cO_{Y_n /S_n} \to  \cO_{Y_n /S_n}\otimes_R \wh \Omega_R $ which is a connection in the following sense. For local sections $f$ of $R$ and $s$ of $\cO_{Y_n /S_n}$, $\nabla(sf) = sdf + \nabla(s)$. Moreover, $\nabla^2 = 0 $ (where $\nabla\colon \cO_{Y_n /S_n}\otimes_R \wh \Omega_R \rightarrow \cO_{Y_n /S_n}\otimes_R \wh \Omega^2_R$ is defined in the usual way).

We claim that this construction is functorial for maps $h\colon A_{(0)} \to B_{(0)}$ in $\cC_{Y_n}$, and so does not depend on the presentation $A_{(0)} = \Rbar_n[x_1, \dots, x_m]/ (f_1 , \dots f_c)$. First note that for any $a \in \Rbar_n[x_1, \dots, x_m]$, we have $\Rbar_n[x_1, \dots, x_m, x_{m+1}] / (x_{m+1}-a) \cong \Rbar_n[x_1, \dots, x_m]$, and $(f_1, \ldots, f_c, x_{m+1}-a)$ is Koszul-regular in $\Rbar_n[x_1, \dots, x_m, x_{m+1}]$. The natural map $A_{(0)} = \Rbar_n[x_1, \dots, x_m]/ (f_1, \dots f_c) \rightarrow \Rbar_n[x_1, \dots, x_m, x_{m+1}] / (f_1, \ldots, f_c, x_{m+1}-a)$ lifts to a map $A \rightarrow \Rbar_n[x^{\frac{1}{p ^\infty}}_1, \dots, x^{\frac{1}{p ^\infty}}_m, x^{\frac{1}{p ^\infty}}_{m+1}] / (f_1 , \dots, f_c, x_{m+1}-a)$. Via a diagram chase, we deduce that using the presentation $\Rbar_n[x_1, \dots, x_m, x_{m+1}] / (f_1, \ldots, f_c, x_{m+1}-a)$ for $A_{(0)}$ gives the same $\varphi_r$ and $\nabla$. Now, if $B_{(0)} \cong \Rbar_n [y_1, \dots, y_{\ell}]/ (g_1 , \dots , g_s)$,  then there exists a map $\hat{h}\colon \Rbar_n [x_1, \dots, x_m ] \to \Rbar_n [y_1, \dots, y_{\ell}]$ which lifts $h$. Furthermore, we may use the presentation $B_{(0)} \cong \Rbar_n [y_1, \dots, y_{\ell}, z_1, \ldots, z_m]/ (g_1 , \dots , g_s, z_1-\hat{h}(x_1), \ldots, z_m-\hat{h}(x_m))$. Then $\hat{h}$ extends to a morphism $$\hat{h}\colon  A \to B \coloneqq \Rbar_n [y^{\frac{1}{p^\infty}}_1, \dots, y^{\frac{1}{p ^\infty}}_{\ell}, z^{\frac{1}{p^\infty}}_1, \dots, z^{\frac{1}{p ^\infty}}_m]/(g_1 , \dots, g_s, z_1-\hat{h}(x_1), \ldots, z_m-\hat{h}(x_m))$$ given by $\hat h (x_i^{\frac{1}{p ^r}}) = z_i^{\frac{1}{p ^r}}$. Thus, $\hat{h}$ induces a map $\cJ ^{[r]}_{Y_n/S_n}(A) \to \cJ ^{[r]}_{Y_n/S_n}(B)$ which is compatible with the map $\cJ ^{[r]}_{Y_n/S_n}(A_{(0)})  \to \cJ ^{[r]}_{Y_n/S_n}(B_{(0)})$ induced by $h$, and a diagram chase gives the desired functoriality. This completes the construction of $\varphi_r$ and $\nabla$ when $Y = \Xbar$, and by construction they satisfy the properties analogous to those stated in Corollary \ref{cor-varphiandconenction}. We also remark that these constructions are natural in $Y_n$.

When $Y = X$, we construct Frobenius and connection on $\cJ ^{[r]}_{Y_n /S_n}$ in an analogous manner, starting with $A_{(0)} = R_n[x_1,\ldots,x_m]/(f_1, \ldots, f_c)$ (where $f_i$'s form a regular sequence), and setting $A = \wt R_n [x_1^{\frac{1}{p^\infty}}, \ldots, x_m^{\frac{1}{p^\infty}}]/(f_1, \ldots, f_c)$. Note that $\wt R$ is a direct limit of syntomic covers of $R$, and therefore $A$ is a direct limit of syntomic covers of $A_{(0)}$. 

The following corollary summarizes the above discussion.

\begin{corollary}\label{cor:firstpartmainthm}
Suppose $Y_n = \Xbar_n$ or $X_n$. Let $0 \leq r \leq p-1$ and $0 \leq i \leq p-1$.
\begin{enumerate} 
\item[(1)] There are natural (in $Y_n$) morphisms $\varphi_r\colon \cJ^{[r]}_{Y_n /R_n} \rightarrow \cO_{Y_n/R_n} $ in $(Y_n)_{\syn}$ such that $\varphi_r$ is $\varphi$-semilinear, and for $r \leq p-2$, $\varphi_r|_{\cJ^{[r+1]}_{Y_n /R_n}} = p \varphi_{r+1}$. 
\item[(2)] There is a natural flat connection $\nabla\colon \cO_{Y_n/R_n} \rightarrow \cO_{Y_n/R_n} \otimes_R \wh \Omega_{R}$ such that:
\begin{enumerate}
 \item $\nabla ( \cJ^{[r]}_{Y_n /R_n} ) \subset \cJ^{[r-1]}_{Y_n /R_n}\otimes_R \wh \Omega_{R}$.
 \item $\nabla\circ \varphi _r = (\varphi_{r-1} \otimes d \varphi _1) \circ \nabla$ (when $r = 0$, $\varphi_{-1}$ is understood to be $p\varphi_0$).
\end{enumerate}
 \item[(3)] There are natural morphisms $\varphi_r\colon \rH^i_\cris (Y_n/R_n, J_{Y_n/R_n}^{[r]}) \rightarrow \rH^i_\cris (Y_n/R_n,\cO_{Y_n/R_n}) $ and connection $\nabla_r\colon \rH^{i}_{\cris}(Y_n/R_n, J_{Y_n/R_n}^{[r]}) \rightarrow \rH^i_\cris (Y_n/R_n, J_{Y_n/R_n}^{[r-1]}) \otimes_R \wh \Omega_{R}$ satisfying the analogous properties as in (1) and (2). When $Y_n = X_n$, each $\rH^i_\cris (Y_n/R_n, J_{Y_n/R_n}^{[j]})$ is finite as an $R_n$-module, and 
 \[
 (\rH^{i}_{\cris}(Y_n/R_n,\cO_{Y_n/R_n}), \rH^i_\cris (Y_n/R_n, J_{Y_n/R_n}^{[r]}),\varphi_r, \nabla_r)
 \]
 gives an object of $\rM\rF_{\nabla, \mathrm{big}}(R)$. When $Y_n = \overline{X}_n$, the above tuple also gives an object of $\rM\rF_{\nabla, \mathrm{big}}(R)$ if we set in addition that $\rF^j = 0$ (so $\varphi_j = 0 = \nabla_j$) for $j \geq p$. 
\end{enumerate}
\end{corollary}

\begin{proof}
It remains to prove the last assertion. By passing to cohomology and using Lemma \ref{lem:divpowersheafcomp}, we see that the resulting data $(\rH^{i}_{\cris}(Y_n/R_n,\cO_{Y_n/R_n}), \rH^i_\cris (Y_n/R_n, J_{Y_n/R_n}^{[r]}),\varphi_r, \nabla_r)$ satisfies the analogous properties as in (1) and (2). When $Y_n = X_n$, each $\rH^i_\cris (X_n/S_n, J_{X_n/R_n}^{[j]})$ is a finite $R_n$-module by \cite[Thm. 14]{Faltings} since $X$ is proper and smooth over $R$. Recall from Section \ref{subsection:crys-cohom-base-change-general-case} that we have the injective map
\[
\bar{\alpha}\colon \rH^i_\cris (X_n/R_n, J_{X_n/R_n}^{[j]})/t_d \rightarrow \rH^i_\cris (X_{d, n}/R_{d, n}, J_{X_{d, n}/R_{d, n}}^{[j]}) 
\]
given by crystalline base change. Since $\rH^i_\cris (X_{s, n}/W_n(k), J_{X_{s, n}/W_n(k)}^{[j]}) = 0$ for $j \geq p$ by the classical Fontaine--Messing theory (see Theorem \ref{Thm-FM}), we deduce by induction on $d$ and Nakayama's Lemma that $\rH^i_\cris (X_n/R_n, J_{X_n/R_n}^{[j]}) = 0$ for $j \geq p$. This proves (3).
\end{proof}

\begin{rem}
When $Y= \Xbar$ and $S$ is either $\Spec(\wt R)$ or $\Spec(W(k))$, we construct $\varphi_r$ for $0 \leq r \leq p-1$ on $\cJ ^{[r]}_{Y_n /S_n}$ in an analogous manner. As before, these construction are functorial in $Y_n$. 
\end{rem}

\begin{rem} For $X_n / R_n$, one can also prove the existence of $\varphi_r$ and $\nabla$ as in the previous corollary via de Rham cohomology. In particular, one can identify crystalline cohomology with de Rham cohomology, and directly construct $\varphi_r$. In this case, $\nabla$ is the usual Gauss--Manin connection. One can show that these two constructions coincide.   
\end{rem}

\begin{proof}[C\textbf{ompletion of proof of Theorem} \ref{thm-cons-phi}]
The first part of the theorem is contained in the previous corollary. It remains to prove the second part. Suppose we have a cartesian diagram 
$$
\xymatrix{
X' \ar[r]^{f} \ar[d] & X \ar[d] \\
\Spec(R') \ar[r] & \Spec(R) }
$$
as in Theorem \ref{thm-cons-phi} (2). First note that it suffices to show that the following induced diagram of sheaves on the syntomic site commutes:
$$
\xymatrix{
\cJ ^{[r]}_{X_n /R_n} \ar[r] \ar[d]^{\varphi_r} & f_*\cJ ^{[r]}_{X'_n /R'_n} \ar[d]^{f_{*}(\varphi_r)}\\
\cO_{X_n/R_n} \ar[r]    & f_*\cO_{X'_n/R'_n} }
$$
Assuming that the previous diagram commutes, we obtain the desired result from the functoriality of the comparison isomorphisms between syntomic and crystalline cohomology (see Remark \ref{rem:topoi-functoriality}). We can check this locally on objects of the category $\cC_{X_{n}} \subset (X_n)_{\syn}$. Let $A_0 \in \cC_{X_{n}}$ be as above, and similarly for $A$. Then $A'_0 \coloneqq A_0 \otimes_{R_n} R'_n \in \cC_{X'_{n}}$ and we have $A' = A \otimes_{R_n} R'_n $. By construction of $\varphi_r$, it suffices to show that the following diagram commutes:
$$
\xymatrix{
\cJ ^{[r]}_{X_n  /R_n}(A) \ar[r] \ar[d]^{\varphi_r} &\cJ ^{[r]}_{X'_n /R'_n}(A') \ar[d]^{\varphi_r}\\
\cO_{X_n /R_n}(A) \ar[r]    &\cO_{X'_n /R'_n}(A') }
$$
Moreover, by definition of $\varphi_r$, it suffices to show that the corresponding diagram with $\varphi_r$ replaced by $\varphi$ commutes. We may further reduce to the setting where frobenius on $A/p$ is surjective. More precisely, we can consider $A$ as in the paragraph before Corollary \ref{cor:firstpartmainthm}. But then compatibility of Frobenius follows from the fact that $R \rightarrow R'$ is compatible with Frobenius and Witt vector functoriality. 
\end{proof}

We simplify notation by setting $\overline \cJ^{[r]}_n \coloneqq \cJ^{[r]}_{\Xbar _n/W_n(k)} ,\ \wt  \cJ^{[r]}_n \coloneqq \cJ^{[r]}_{\Xbar _n/\wt R_n},\ \cJ^{[r]}_n \coloneqq \cJ^{[r]} _{\Xbar _n/R_n} $, and as before $\cJ^{[0]}_{Y_n /S_n} = \cO_{Y_n /S_n}$ and $\wt \cJ^{[0]}_{Y_n /S_n} = \wt \cO_{Y_n /S_n}$. The natural inclusions $W(k) \subset R \subset \wt R$ and functoriality of crystalline cohomology induce the following morphisms on the syntomic site of $\Xbar_n$:
$$\overline \cJ^{[r]}_n \overset{i}{\to } \cJ^{[r]}_n \overset{j}{\to} \wt \cJ^{[r]}_n.$$ 

\begin{proposition}\label{prop-structureofJ} Let $0 \leq r \leq p-2$. 
\begin{enumerate}
\item[(1)] $i$ is injective and compatible with $\varphi_r$. Furthermore, it induces an isomorphism $\overline \cJ^{[r]}_n \simeq (\cJ^{[r]}_n )^ {\nabla =0}$.
\item[(2)] $j$ is surjective and compatible with $\varphi_r$. Moreover, $j \circ i$ is the identity map. 
\end{enumerate}
\end{proposition}

\begin{proof} It suffices to prove the assertions at the level of sections on $\Spec A$ where $A$ is assumed to be as in Proposition \ref{prop-calculateH0-1} and Proposition \ref{prop-calculateH0}. Using the isomorphism $f$ in Lemma \ref{lem-toOAcris} to identify these local sections, $i$ and $j$ can be identified with the morphisms $A_{\cris, n}(A) \to A_{\cris, n}(A) \langle X_1 , \dots , X_d \rangle \to A_{\cris, n}(A)$ where the first map is the natural inclusion and the second map is induced by $X_i \mapsto 0$. Both maps are compatible with $\varphi_r$, and 
\[
(A_{\cris, n}(A) \langle X_1 , \dots , X_d \rangle)^{\nabla = 0} = A_{\cris, n}(A).
\]
\end{proof}

The following lemma will be used later to build the comparison map between crystalline and \'etale cohomology.

\begin{lemma}\label{lem-varphir-surjective} For $0 \leq r \leq p-2$, the map $\varphi_r-1\colon \wt \cJ_n ^{[r]} \to \wt \cO_n$ is surjective. 
\end{lemma}

\begin{proof} 
By Lemma \ref{lem-info-fil}, we are reduced to showing the desired surjectivity for the case $n = 1$. Let $U = \Spec(A)$ with $A$ as in Lemma \ref{lem-info-fil} (2), i.e. $A= \wt R_n[x_1^{\frac{1}{p ^\infty}}, \ldots, x_m^{\frac{1}{p ^\infty}}]/ (f_1, \ldots, f_c)$ and Frobenius on $A_1$ is surjective. In this case, $\varphi_r-1\colon \wt \cJ_1^{[r]} (U) \to \wt \cO_1(U)$ can be identified with $\varphi_r-1\colon J^{[r]}_{W_1^\PD(A_1^\flat)} \to W_1^\PD( A_1^\flat)$. It suffices to show for any $y \in W_1^{\PD} (A_1^\flat)$, there exists a direct limit of syntomic $A$-algebras $B$ such that $\varphi$ on $B$ is surjective (we are assuming $n = 1$) and that there exists $x \in J^{[r]}_{W^{\PD}_1(B_1^\flat)}$ satisfying $\varphi_r(x) -x = y$. Moreover, it suffices to find such an algebra $B$ and $x \in J^{[r]}_{W^{\PD}_1(B_1^\flat)}$ so that $\varphi_r (x) - x \equiv y \mod J^{[p-1]}_{W_1^\PD(B_1^\flat)}$, since $\varphi_r$ on $J^{[p-1]}_{W_1^\PD(B_1^\flat)}$ is $0$.

We are reduced to showing that given $y = (y^{(i)})_{i \geq 0} \in A_1^\flat \cong W_1(A_1^\flat)$, we can find $B$ and $x$ as above such that $\varphi_r (x) -x \equiv y \mod J^{[p-1]}_{W_1^\PD(B_1^\flat)}$. Replacing $A$ by a suitable direct limit of syntomic $A$-algebras if necessary, we have a non-zero element $(q_i)_{i \geq 0} \in A_1^{\flat}$ such that $q_0 = -p$. Let $z = [(q_i)]+p \in  W(A_1^\flat)$, and write $\bar z = z \mod p \in W_1^\PD (A_1^\flat)$. Note that 
 $\varphi_r({\bar z}^r) -1 \in J^{[p-1]}_{W_1^\PD(A_1^\flat)}$. Let $B^{(0)} = A[X]/(X^p-q_1^rX-y^{(1)})$, which is syntomic over $A$, and consider $B = \varinjlim_{j} B^{(j)}$ where $B^{(j+1)} \coloneqq B^{(j)}[(T_b)_{b \in B^{(j)}}] / ((T_b^p-b)_{b \in B^{(j)}})$ for $j \geq 0$. Let $w = (w^{(i)})_{i \geq 0} \in B_1^\flat$ such that $w^{(1)} \in B$ is the image of $X$ in $B$, and let $x = \bar z^r w \in J^{[r]}_{W^{\PD}_1(B_1^\flat)}$. Since $w^{(1)}$ in $B$ satisfies the equation $X^p-q_1^rX-y^{(1)} = 0$, we conclude that $\varphi_r (x) - x \equiv y \mod J^{[p-1]}_{W_1^\PD(B_1^\flat)}$.
\end{proof}

The cohomology groups $\rH^i_\syn (\Xbar_n, \cJ^{[r]}_{n})$ have a natural $G_R$-action, and  there are several equivalent ways to define this action. Given an element $g \in G_R$, we have a natural induced morphism of the corresponding syntomic topoi, and moreover a natural morphism $g^*(\cJ^{[r]}_{n}) \rightarrow \cJ^{[r]}_{n}$ or equivalently a morphism $\cJ^{[r]}_{n} \rightarrow g_*\cJ^{[r]}_{n}$. The latter morphism is given by an application of crystalline functoriality. We recall an explicit construction of this action to see that this morphism is $\varphi_r$-compatible. Consider $U_0 = \Spec(A_0)$ with $A_0= \Rbar_n[x_1,\ldots,x_m]/I$, where $I$ is a Koszul-regular ideal and $A_0$ is flat over $\Rbar_n$. For each $g \in G_R$, $g(I) \subset \Rbar_n[x_1,\ldots,x_m]$ is a Koszul-regular ideal and $g(A_0) \coloneqq \Rbar_n[x_1,\ldots,x_m] / g(I)$ is $\Rbar_n$-flat. By functoriality of crystalline cohomology, the $g$-semi-linear map $A_0 \rightarrow g(A_0)$ of $\Rbar_n$-algebras induces a morphism $\cJ_n^{[r]}(U_0) \rightarrow \cJ_n^{[r]}(g(U_0))$ of $R_n$-algebras, where $g(U_0) \coloneqq \Spec(g(A_0))$. We have a similar statement for $U = \Spec(A)$ with $A =\Rbar[x_1^{\frac{1}{p^{\infty}}},\ldots,x_m^{\frac{1}{p^{\infty}}}]/I$. By Remark \ref{rem-GroupAction}, the induced map $\cJ_n^{[r]}(U) = \rF^r\cO A _{\cris,n}(A) \rightarrow \cJ_n^{[r]}(g(U)) = \rF^r\cO A _{\cris,n}(g(A))$ is compatible with $\varphi_r$ and connection. By a similar argument as above using exact sequences of the form 
\[
\xymatrix {\cJ^{[r]}_{n}(U_0)  \ar[r] &  \cJ^{[r]}_{n}(U)  \ar@<1ex>[r]\ar[r] & \cJ^{[r]}_{n}(U \times_{U_0} U)},
\] 
we conclude that the natural induced $G_R$-action on $\cJ^{[r]}_{n}$ is compatible with $\varphi_r$,  connection and natural map $\cJ^{[r]}_{n} \to \cJ^{[r-1]}_{n} $. Consequently, the $G_R$-action on $\rH^i_\syn (\Xbar_n, \cJ^{[r]}_{n})$ is compatible with $\varphi_r$ and connection.  

\begin{rem}
The $G_R$-action can also be defined directly over the crystalline-syntomic site or the crystalline site via functoriality. These are all compatible with the action defined above by the comparison isomorphisms between these cohomology groups.
\end{rem}

For $0 \leq r \leq p-2$, consider the following sheaf on $(\Xbar_n)_{\syn}$:
$$ S^{[r]}_n \coloneqq (\cJ_n^{[r]}) ^{\varphi_r =1, \nabla= 0}. $$ 
Since the $G_R$-action on $\cJ^{[r]}_{n}$ is compatible with $\varphi_r$ and connection, we have a natural induced $G_R$-action on $S^{[r]}_n$ and thus on $\rH^i_\syn (\Xbar_n, S^{[r]}_n)$. By Proposition \ref{prop-structureofJ}, we have $S_n ^{[r]} = (\wt \cJ ^{[r]}_n )^{\varphi_r=1}$, and by Lemma \ref{lem-varphir-surjective}, we have an exact sequence over $(\Xbar_n)_{\syn}$:
\begin{equation}\label{eqn-S-r}
\xymatrix{0 \ar[r] & S^{[r]}_n \ar[r] & \wt \cJ^{[r]}_n \ar[r]^{\varphi_r -1} & \wt \cO_n \ar[r] & 0. }
\end{equation}
We shall use this sequence in the following sections to construct the desired comparison map from crystalline to etale cohomology.

\subsection{Construction of the comparison map, Part 1}

We now construct a natural $G_R$-equivariant morphism for $i \leq r \leq p-2$:
$$\iota_R\colon T_\cris( \rH^i_\cris(X_n / R_n))(r) \rightarrow \rH^i_\syn (\Xbar_n, S^{[r]}_n).$$

First, consider the natural $G_R$-equivariant `Kunneth' morphisms:
$$\alpha'_{R}\colon \rH^{0}_{\cris}(\Rbar_n/R_n) \otimes_{R_n} \rH^i _\cris(X_n/R_n) \longrightarrow \rH^i_\cris (\Xbar_n / R_n)$$
and
$$\beta'_{R}\colon \bigoplus_{j =0}^r\left ( \rH^0_{\cris}(\Rbar_n/R_n, J_{\Rbar_n/R_n}^{[j]}) \otimes_{R_n} \rH^i_\cris (X_n/ R_n, J_{X_n/R_n}^{[r-j]} )\right ) \longrightarrow \rH^i _\cris (\Xbar_n/ R_n, J^{[r]}_{\Xbar _n /R_n} ).   $$
By Corollary \ref{cor-acris} and Lemma \ref{lem:divpowersheafcomp}, we obtain $G_R$-equivariant morphisms:
$$\alpha_{R}\colon \cO A_{\cris,n} (R) \otimes_{R_n} \rH^i_\syn(X_n, \cO_{X_n / R_n}) \longrightarrow \rH^i_\syn (\Xbar_n, \cO _{\Xbar_n /R_n}),$$
$$\beta_{R}\colon \bigoplus_{j =0}^r\left ( \rF^j \cO A_{\cris,n} (R) \otimes_{R_n} \rH^i_\syn (X_n, \cJ_{X_n/R_n}^{[r-j]} )\right ) \longrightarrow \rH^i _\syn (\Xbar_n, \cJ^{[r]}_{\Xbar_n /R_n} ).   $$
By construction, $\varphi_r$ and $\nabla$ on $\cJ ^{[r]}_{Y_n /S_n}$ are functorial in $Y_n$, and so the morphisms $\alpha_R$ and $\beta _R$ are compatible with these structures. 
When $i \leq p-2$, by Theorem \ref{thm-cons-phi},  $M \coloneqq (\rH_\cris ^i(X_n/ R_n), \rH^i_\cris (X_n/R_n, J^{[j]}_{X_n/ R_n}), \varphi_j, \nabla)$ is an object in $\MF_{\nabla}^{[0, i]}(R)$, and thus admits an adapted basis by Remark \ref{rem-adpated}. It follows that $\beta_R$ induces the following map (which we also denote by $\beta_R$):
$$ \beta_{R}\colon \rF^r (\cO A_{\cris,n} (R) \otimes_{R_n} \rH_\cris ^i (X_n/ R_n))\longrightarrow \rH^i _\syn (\Xbar_n, \cJ^{[r]}_{\Xbar_n /R_n} ).$$

A similar construction involving $\wt R_n$, $\wt X \coloneqq \wt R \times_R X$, and $\wt \cJ^{[r]}_n $ gives natural maps 
$$\alpha_{\wt R }\colon  A_{\cris,n} (R) \otimes_{R_n} \rH^i _\cris( X_n/R_n) \cong  A_{\cris,n} (R) \otimes_{\wt R_n} \rH^i _\cris(\wt X_n/\wt R_n)\longrightarrow \rH^i _\cris (\Xbar_n / \wt R_n)\cong \rH ^i_\syn ( \Xbar_n, \wt \cO_n),$$
$$ \beta_{\wt R}\colon \rF^r ( A_{\cris,n} (R) \otimes_{R_n} \rH_\cris ^i (X_n/ R_n))  \longrightarrow \rH^i _\cris (\Xbar_n/ \wt R_n, J^{[r]}_{\Xbar_n /\wt R_n} ) \cong \rH^i _\syn (\Xbar_n, \wt \cJ ^{[r]}_n ).$$
Here we use the isomorphism $\rH^i _\cris (\wt X_n / \wt R_n)\cong \wt R \otimes_R \rH^i _\cris (X_n /R_n)$, which follows from Theorem \ref{thm-cry-BC} and Lemma \ref{lem-faithflat}. 

Consider the natural map from $T_\cris (M)$ to the kernel of the map
$$\varphi_r -1\colon \rH^i_\syn (\Xbar_n, \cJ_n ^{[r]})^{\nabla = 0} \to \rH^i_\syn (\Xbar_n, \cO_n)^{\nabla = 0}. $$
We obtain the commutative diagram
\begin{equation}\label{diag-G_R-actions}
\xymatrix{ T_\cris (M)(r)\ar[r] \ar[d] &   \rH^i_\syn (\Xbar_n, \cJ_n ^{[r]})^{\nabla = 0} \ar[r]^{\varphi_r -1} \ar[d] &   \rH^i _\syn (\Xbar_n, \cO_n )^{\nabla = 0} \ar[d]\\ K \ar[r] &   \rH^i_\syn (\Xbar_n, \wt \cJ_n ^{[r]}) \ar[r]^{\varphi_r -1} &   \rH^i _\syn (\Xbar_n, \wt \cO_n )}
\end{equation}
where $K$ denotes the kernel of the bottom horizontal map $\varphi_r-1$. 

\begin{proposition}\label{prop-map-exists} For $i \leq  r\leq p -2$, we have a natural isomorphism $\rH^i_\syn (\Xbar, S_n ^{[r]}) \cong K$ such that the induced morphism $\iota_R\colon T_\cris (M)(r) \to \rH^i_\syn (\Xbar_n, S^{[r]}_n)$ is compatible with $G_R$-actions. 
\end{proposition}

\begin{proof}
The exact sequence \eqref{eqn-S-r} yields the long-exact sequence 
$$\xymatrix{\cdots \rightarrow \rH ^i _\syn(\Xbar_n, S^{[r]}_n)  \ar[r] & \rH^i_\syn (\Xbar_n, \wt \cJ^{[r]}_n)  \ar[r]^-{\varphi_r -1} & \rH^i _\syn(\Xbar_n, \wt \cO_n)\ar[r] &  \rH ^{i+1} _\syn(\Xbar_n, S^{[r]}_n)  \rightarrow \cdots  }.$$
It suffices to prove that $\varphi_r -1$ is surjective for $0\leq i\leq p-3$ and $i \leq r \leq p-2$. Consider the commutative diagram 
$$\xymatrix{ \rF^r(A_{\cris,n} (R) \otimes _{R_n} \rH^i_\cris (X_n / R_n))\ar[d]^{\beta_{\wt R}} \ar@{->>}[r]^-{\varphi_r -1} & A_{\cris, n}(R) \otimes_{R_n} \rH_\cris^i (X_n / R_n) \ar[d]^{\alpha_{\wt R} }_\wr  \\
 \rH^i_\syn (\Xbar_{n} , \wt \cJ^{[r]}_n)  \ar[r]^-{\varphi_r -1} & \rH^i _\syn(\Xbar_{n}, \wt \cO _n)} $$
By Theorem \ref{thm-Faltings2}, $\varphi_r -1$ in the top row is surjective. Since $\alpha_{\wt R}$ is an isomorphism by Lemma \ref{lem-fromRtoRbar} below, we conclude that the bottom row is surjective. 
 
Thus, we obtain a morphism $\iota_R\colon T_\cris (M)(r) \to \rH^i_\syn (\Xbar_n, S_n^{[r]})$. It remains to check that $\iota_R$ is compatible with $G_R$-actions. Note that this is not obvious a priori because $\wt \cJ^{[r]}_n$ does not have a natural $G_R$-action (although it has a natural $G_\infty$-action). By Proposition \ref{prop-structureofJ}, we can identify $\overline \cJ^{[r]}_n \cong \wt \cJ_n^{[r]}$ and get the following commutative diagram
$$\xymatrix {T_\cris (M)(r)\ar[dr]  \ar@/^1pc/[rr]& \rH^i_\syn (\Xbar_n, \overline \cJ_n ^{[r]}) \ar[dr]^\simeq  \ar@{^(->}[r] & \rH^i_\syn (\Xbar_n, \cJ_n^{[r]}) \ar[d] \\ &  \rH ^i_\syn (\Xbar_n, S_n^{[r]}) \ar@{^(->}[u]\ar@{^(->}[r] & \rH^i_\syn (\Xbar_n, \wt  \cJ_n^{[r]})}$$
In particular, the $G_R$-equivariant map $T_\cris (M)(r) \rightarrow \rH^i_\syn (\Xbar_n, \cJ_n^{[r]})$ factors through $\iota_R$ with $\rH^i_\syn (\Xbar_n, S_n^{[r]}) \rightarrow \rH^i_\syn (\Xbar_n, \cJ_n^{[r]})$ being injective, and thus $\iota_R$ is $G_R$-equivariant. 
\end{proof}

\begin{lemma}\label{lem-fromRtoRbar} 
For $i \leq p-3$, $\alpha_{\wt R}$ is an isomorphism.  
\end{lemma}

\begin{proof} 
The usual base change map in crystalline cohomology (see Appendix \ref{appendix-crystalline-site}) gives rise to a hypertor spectral sequence with
	$$\rE_{2}^{m,q} = \rH^{m}(A_{\cris,n}(R) \otimes_{\wt R_n}^{\bbL} \rH^q_\cris(\wt X_n/ \wt R_n, \cO_{\wt X_n/ \wt R_n}))$$
	which converges to $$\rH^{m+q}_{\cris}(\Xbar_n/ \wt R_n, \cO_{\Xbar_n/\wt R_n} ).$$
Note that the $\rE_{2}^{p,q}$ terms are concentrated in the second quadrant, and after reindexing cochain complexes to chain complexes, we have 
$${\mathrm Tor}^{m}_{{\wt R_n}}(A_{\cris,n}(R), ~\rH^{q}_\cris({\wt X}_n/\wt R_n,\cO_{\wt X_n/ \wt R_n} )) =  \rH^{-m}(A_{\cris,n}(R) \otimes_{\wt R_n}^{\bbL} \rH^q_\cris(\wt X_n/ \wt R_n, \cO_{\wt X_n/ \wt R_n})).$$
For $i \leq p-2$, $\rH^i_\cris (X_n / R_n, \cO_{X_n/R_n})$ is a Fontaine--Laffaille module, so by Theorem \ref{thm-Faltings1}, it is a direct sum of $R/p^{m_{\ell}} R$ for some integers $m_{\ell}$'s. By the flat base change $R_n \rightarrow \wt R_n$, $\rH^{q}_\cris({\wt X}_n/\wt R_n,\cO_{\wt X_n/ \wt R_n})$ is a direct sum of $\wt{R}/p^{m_{\ell}}\wt{R}$ for some integers $m_{\ell}$'s when $q \leq p-2$. Since $A_{\cris,n}(R)$ is $\Z/p^n\Z$-flat, we have $\rE_{2}^{m,q} = 0$ when $q \leq p-2$ and $m < 0$. Hence, $\alpha_{\wt R}$ is an isomorphism when $i \leq p-3$. 
\end{proof}

\begin{rem} 
We do not know whether $A_{\cris}(R)$ is flat over $\wt R$ in general, in which case $\alpha_{\wt R}$ would be an isomorphism for any $i$.
\end{rem}

\subsection{Construction of the comparison map, Part 2}

Following the approach of Fontaine--Messing, we explain how to obtain a natural $G_R$-equivariant morphism $$\rH^i_\syn(\Xbar_n, S^{[r]}_n) \to \rH ^i _\et (X_{\Rbar [\frac 1 p]}, \Z/ p^n \Z(r)).$$

For any $\Z$-module $A$, we denote by $\wh A$ its $p$-adic completion. For $\Xbar = X \times_R \Rbar$ as before, let $\wh {\Xbar}$ be the corresponding $p$-adic formal scheme. We shall construct a diagram of topoi
$$
\xymatrix{
(\Xbar_n)_{\syn} \ar[r]^{ \iota} &  \wh{\Xbar}_{\syn-\et} \ar[r]^{i}& \Xbar_{\syn-\et} &  (X_{\Rbar[1/p]})_{\et} \ar[l]_{j} \\
}
$$
and a sheaf $\cS_n^{[r]}$ on $\Xbar_{\syn-\et}$ whose restriction to $\wh{\Xbar}_{\syn-\et} $ is $\iota_*(S_n^{[r]})$ and whose pull-back to $(X_{\Rbar[1/p]})_{\et}$ is 
$\bbZ/p^n\bbZ(r)$.\footnote{In the following, we work with the small sites.} This will allow us to construct the desired comparison map. The constructions below follow the approach of Fontaine--Messing very closely. The main technical modification is that instead of passing to the rigid analytic generic fiber in defining the site $\wh{\Xbar}_{\syn-\et}$, we work directly with restricted power series algebras. 

We begin by defining the \emph{syntomic-\'etale} site of the $p$-adic formal scheme $\wh {\Xbar}$. For any morphism of $p$-adic affine formal schemes $g\colon \mathrm{Spf} \mathcal{B} \rightarrow \mathrm{Spf} \mathcal{A}$ over $\wh{\Xbar}$, we say $g$ is syntomic if for each integer $m \geq 1$, the associated morphism $\mathcal{A}_m \rightarrow \mathcal{B}_m$ is syntomic. We say $g$ has \'etale generic fiber if it is of finite presentation, and the associated ring map $\mathcal{A}[\frac{1}{p}] \rightarrow \mathcal{B}[\frac{1}{p}]$ is flat and the universal finite differential module $\Omega^f_{\mathcal{B}[\frac{1}{p}] / \mathcal{A}[\frac{1}{p}]}$ is $0$. Note that if $g$ has \'etale generic fiber, then $\mathcal{B}$ has a presentation $\mathcal{B} = \mathcal{A}\langle X_1, \ldots, X_m \rangle / (P_1, \ldots, P_m)$ such that the Jacobian $J = \det{(\frac{\partial P_i}{\partial X_j})}$ is invertible in $\mathcal{B}[\frac{1}{p}]$. The objects of the (small) syntomic-\'etale site of ${\wh{\Xbar}}$ are defined to be the morphisms of $p$-adic formal schemes $\mathcal{Y} \rightarrow {\wh{\Xbar}}$ which are locally quasi-finite, syntomic, and have \'etale generic fiber. The topology is generated by surjective families of open immersions and finite surjective families of locally quasi-finite syntomic morphisms having \'etale generic fibers with source and target both affine schemes. From the immersion $\Xbar_n \rightarrow \wh \Xbar$, we have a natural morphism of topoi 
$$\iota\colon (\Xbar_n)_{\syn} \rightarrow \wh{\Xbar}_{\syn-\et}$$ where $\iota_*$ is given by $\iota_*(F)(Y) = F(Y_n)$. 

\begin{lemma}\label{lem-iotastar}
	For any $n \geq 1$, $\iota_*$ is an exact functor. 
\end{lemma}

\begin{proof}
By \cite[\href{https://stacks.math.columbia.edu/tag/0570}{Tag 0570}]{stacks-project}, it suffices to prove that if $\mathcal{A}$ is a $p$-adically complete $p$-torsion free $R$-algebra and $f\colon \mathcal{A}_n \rightarrow B$ is a syntomic quasi-finite morphism, then locally $B$ can be lifted to a $p$-adically complete $R$-algebra $\mathcal{B}$ and $f$ can be lifted to $\tilde{f}\colon \mathcal{A} \rightarrow \mathcal{B}$ such that $\mathrm{Spf}\mathcal{B} \rightarrow \mathrm{Spf}\mathcal{A}$ is syntomic, quasi-finite with \'etale generic fiber. By \cite[\href{https://stacks.math.columbia.edu/tag/0DWJ}{Tag 0DWJ}]{stacks-project}, locally $B$ can be written as $B = \mathcal{A}_n[X_1, \ldots, X_m] / (g_1, \ldots, g_m)$ such that $\mathcal{A}_n \rightarrow B$ is a relative global complete intersection. Consider the syntomic map $B \rightarrow B' \coloneqq \mathcal{A}_n[Y_1, \ldots, Y_m] / (g_1(Y^{p^{n+1}}), \ldots, g_m(Y^{p^{n+1}}))$ given by $X_i \mapsto Y_i^{p^{n+1}}$, and let $\mathcal{B} = \mathcal{A}\langle Y_1, \ldots, Y_m\rangle / (\hat{g}_1(Y^{p^{n+1}})+p^nY_1, \ldots, \hat{g}_m(Y^{p^{n+1}})+p^nY_m)$ where $\mathcal{A}\langle Y_1, \ldots, Y_m\rangle$ denotes the $p$-adic completion of $\mathcal{A}[Y_1, \ldots, Y_m]$ and $\hat{g}_i$ is a lift of $g_i$. Note that $\mathcal{B} / p^n\mathcal{B} = B'$. 

Letting $\tilde{g}_i = \hat{g}_i(Y^{p^{n+1}})+p^nY_i$, we claim that $(\tilde{g}_1, \ldots, \tilde{g}_m)$ is Koszul-regular in $\mathcal{A}\langle Y_1, \ldots, Y_m\rangle$. For the Koszul complex $K_{\bullet}(\tilde{g}_1, \ldots, \tilde{g}_m)$, we have $\rH_i(K_{\bullet}(\tilde{g}_1, \ldots, \tilde{g}_m))/p^n \rH_i(K_{\bullet}(\tilde{g}_1, \ldots, \tilde{g}_m)) = 0$ for $i > 0$ since $\mathcal{A}_n \rightarrow B'$ is syntomic and $\mathcal{A}\langle Y_1, \ldots, Y_m\rangle$ is $p$-torsion free. Since the terms in $K_{\bullet}(\tilde{g}_1, \ldots, \tilde{g}_m)$ is $p$-adically complete, it follows that $\rH_i(K_{\bullet}(\tilde{g}_1, \ldots, \tilde{g}_m)) = 0$ for $i > 0$. Moreover, $\mathcal{B}$ is $p$-torsion free by \cite[\href{https://stacks.math.columbia.edu/tag/068M}{Tag 068M}]{stacks-project}.

Hence, the map $\tilde{f}\colon \mathcal{A} \rightarrow \mathcal{B}$ is quasi-finite, syntomic, and has \'etale generic fiber. 
\end{proof}

For $\iota\colon (\Xbar_{n})_{\syn} \rightarrow (\wh{\Xbar})_{\syn-\et}$, we write $S_n^{[r]}$ for $\iota_*S_n^{[r]}$ and $(\mathcal{J}_n^{[r]})^{\varphi_r = 1}$ for $\iota_*(\mathcal{J}_n^{[r]})^{\varphi_r = 1}$ to simplify notations. By the previous lemma, we have $$\rH^{i}((\Xbar_n)_{\syn}, S_n^{[r]}) \cong \rH^i({\wh{\Xbar}}_{\syn-\et}, S_n^{[r]}),$$
and similarly for $(\mathcal{J}_n^{[r]})^{\varphi_r = 1}$.

We define the syntomic-\'etale site of $\Xbar$ so that the objects are morphisms $f: U \rightarrow \Xbar$ which are syntomic, locally quasi-finite, and $f\otimes_\Rbar \Rbar[\frac{1}{p}]$ is \'etale. The topology is generated by surjective families of open immersions and finite surjective families of locally quasi-finite syntomic morphisms with source and target both affine. By definition, we have a natural morphism of topoi $j\colon (X_{\Rbar[\frac{1}{p}]})_{\et} \rightarrow {\Xbar}_{\syn-\et}$. Similarly, we have a natural morphism of topoi $i: \wh{\Xbar}_{\syn-\et} \rightarrow {\Xbar}_{\syn-\et}$ as follows. If $\mathrm{Spec} B$ is an affine object of the {\syn}tomic-\'etale site of ${\Xbar}$, then $i_*F(\mathrm{Spec}B) \coloneqq F(\mathrm{Spf} \wh{B})$. We obtain a description of $i^*$ via the following lemma. 

\begin{lemma}\label{lem-Zar-et}
	Let $\mathrm{Spec} A$ be Zariski open in ${\Xbar}$, and $\mathcal{B}$ a $p$-adically complete algebra equipped with a syntomic, quasi-finite morphism $\tilde{f}\colon \wh{A} \rightarrow \mathcal{B}$ with \'etale generic fiber. Then there exists a syntomic, quasi-finite map $f\colon A \rightarrow B$ with \'etale generic fiber such that $\wh{B} = \mathcal{B}$ and $\hat{f} = \tilde{f}$.	
\end{lemma}

\begin{proof}
	Choose a presentation $\mathcal{B} = \wh{A}\langle X_1, \ldots, X_m\rangle/(P_1, \ldots, P_m)$ such that $J = \det{(\frac{\partial P_i}{\partial X_j})} \in \wh{A}\langle X_1, \ldots, X_m\rangle$ is invertible in $\mathcal{B}[\frac{1}{p}]$. Then there exist $a, a_1, \ldots, a_m \in \wh{A}\langle X_1, \ldots, X_m\rangle$ such that 
	\begin{equation}\label{eqn-Ja}
	Ja = p^r+\sum_{i = 1}^m a_iP_i 
	\end{equation}
	for some integer $r \geq 0$. For each $i = 1,\ldots, m$, let $Q_i \in A[X_1, \ldots, X_m]$ such that $Q_i \equiv P_i$ mod $p^{3r+1}$ in $\wh{A}\langle X_1, \ldots, X_m\rangle$. Let $B = A[X_1, \ldots, X_m]/(Q_1, \ldots, Q_m)$. Note that by a similar argument as in the proof of Lemma \ref{lem-iotastar}, $(Q_1, \ldots, Q_m) \subset \wh{A}\langle X_1, \ldots, X_m\rangle$ is Koszul-regular and $\wh{B}$ is $p$-torsion free. We claim that there exists a unique $\wh{A}$-algebra map $h\colon \mathcal{B} \rightarrow \wh{B}$ which is the identity map mod $p^{3r+1}$.
	
	To prove the claim, we construct a compatible system of $\wh{A}$-algebra maps $h_n\colon \mathcal{B} \rightarrow \wh{B}/p^n\wh{B}$ for each integer $n \geq 3r+1$ inductively as follows. Let $h_{3r+1}$ be given by $X_i \mapsto X_i$. Suppose for $n \geq 3r+1$, $h_n$ is given by $X_i \mapsto X_i+u_i$ where $u_i \in p^{2r+1}\wh{B}$. Then $P_i(X_1+u_1, \ldots, X_m+u_m) \in p^n\wh{B}$ for each $i$. By Taylor series expansion, we have
	\[
	P_i(\mathbf{X}+\mathbf{u}+\mathbf{\mu}) = P_i(\mathbf{X}+\mathbf{u})+\sum_{\ell=1}^m \frac{\partial P_i}{\partial X_{\ell}}(\mathbf{X}+\mathbf{u})\cdot \mu_{\ell}+\sum_{|b| \geq 2}\frac{\partial^b P_i}{\partial X^b}(\mathbf{X}+\mathbf{u})\cdot \frac{\mu^b}{b!}
	\]  
	where $b = (b_1, \ldots, b_m) \in \mathbb{N}^m$ denotes a multi-index in the last term. From equation \eqref{eqn-Ja}, there exist $\mu_1, \ldots, \mu_m \in p^{n-r}\wh{B}$, uniquely determined mod $p^{n-r+1}\wh{B}$, such that $P_i(\mathbf{X}+\mathbf{u})+\sum_{\ell=1}^m \frac{\partial P_i}{\partial X_{\ell}}(\mathbf{X}+\mathbf{u})\cdot \mu_{\ell} \in p^{n+1}\wh{B}$, and thus $P_i(\mathbf{X}+\mathbf{u}+\mathbf{\mu}) \in p^{n+1}\wh{B}$ for each $i$. This proves the claim. Similarly, there exists a unique $\wh{A}$-algebra map $\wh{B} \rightarrow \mathcal{B}$ which is the identity mod $p^{3r+1}$, and therefore $h\colon \mathcal{B} \rightarrow \wh{B}$ is an isomorphism.
	
	Since $B$ is finitely presented over $A$ and $\wh{B}$ is $p$-torsion free, we have that $B$ is $p$-torsion free. Let $j \in B$ be the jacobian. Then for some $g, b \in B$, we have $jg = p^r+p^{r+1}b = p^r(1+pb)$. If we replace $B$ by $B[\frac{1}{1+pb}]$, then the map $f\colon A \rightarrow B$ has \'etale generic fiber. Furthermore, since $\mathcal{A} \rightarrow \mathcal{B}$ is syntomic, $A/pA \rightarrow B/pB$ is flat, and thus $f\colon A \rightarrow B$ is flat. It follows that $f\colon A \rightarrow B$ satisfies the desired properties.
\end{proof}

We may now define $i^*F$ as the sheafification of the presheaf $\mathcal{B} \mapsto \varinjlim F(B')$, where $\varinjlim F(B')$ is the filtered colimit of the commutative diagrams
$$\xymatrix{ B \ar[r] \ar[dr] & B'\ar[d] \\ & \cB}$$
such that $B$ is a fixed choice of a ring satisfying the conditions in Lemma \ref{lem-Zar-et}, $B \rightarrow B'$ is \'etale, and the map $B' \rightarrow \mathcal{B}$ induces an isomorphism $\wh{B}' \cong \mathcal{B}$. We claim that $\varinjlim F(B')$ does not depend on the choice of $B$, and $i^*F$ is well-defined. For this, suppose $B$ and $B'$ are $A$-algebras both satisfying the conditions in Lemma \ref{lem-Zar-et}, and let $B'' \subset \mathcal{B}$ be the $A$-algebra given by the image of the induced map $B\otimes_A B' \rightarrow \mathcal{B}$. Then we claim that the maps $B \rightarrow B''$ and $B' \rightarrow B''$ are \'etale, and $B''$ also satisfies the conditions in Lemma \ref{lem-Zar-et}. First, note that for each $n \geq 1$, the induced map $B/p^n B \rightarrow B''/p^n B''$ is an isomorphism, and $\wh{B}'' \cong \mathcal{B}$. In particular, $A/pA \rightarrow B''/pB''$ is flat, and $B''$ is $p$-torsion free as $B'' \subset \mathcal{B}$. Write $\overline{R} = \varinjlim_{\eta} R_{\eta}$ where the limit goes over finite normal $R$-algebras $R_{\eta} \subset \overline{R}$ such that $R_{\eta}[\frac{1}{p}] / R[\frac{1}{p}]$ is \'etale. Then there exists an index $\alpha$, finitely generated flat $R_{\alpha}$-algebra $A_{\alpha}$, and finitely generated $A_{\alpha}$-algebra $B_{\alpha}$ such that $\displaystyle A \cong \varinjlim_{\eta \geq \alpha} A_{\alpha}\otimes_{R_{\alpha}} R_{\eta}$ and $\displaystyle B'' \cong \varinjlim_{\eta \geq \alpha} B_{\alpha}\otimes_{R_{\alpha}} R_{\eta}$. For indices $\eta \geq \alpha$, denote $A_{\eta} = A_{\alpha}\otimes_{R_{\alpha}} R_{\eta}$ and $B_{\eta} = B_{\alpha}\otimes_{R_{\alpha}} R_{\eta}$. Then there exists an index $\beta \geq \alpha$ such that $A_{\beta}/pA_{\beta} \rightarrow B_{\beta}/pB_{\beta}$ is flat. Furthermore, the natural map $\mathrm{Tor}_1^{A_{\eta}}(B_{\eta}, A_{\eta}/pA_{\eta})\otimes_{R_{\eta}} R_{\lambda} \rightarrow \mathrm{Tor}_1^{A_{\lambda}}(B_{\lambda}, A_{\lambda}/pA_{\lambda})$ is surjective for any $\lambda \geq \eta \geq \beta$. Since $\mathrm{Tor}_1^A(B'', A/pA) = 0$, we have $\mathrm{Tor}_1^{A_{\gamma}}(B_{\gamma}, A_{\gamma}/pA_{\gamma}) = 0$ for some $\gamma \geq \beta$. Then by the local criterion of flatness, $B_{\gamma}$ is flat over $A_{\gamma}$, and so $B''$ is flat over $A$. By a similar argument, we deduce that $B''$ is flat over both $B$ and $B'$. Since $(B\otimes_A B')[\frac{1}{p}]$ is \'etale over $B[\frac{1}{p}]$ and $(B\otimes_A B')[\frac{1}{p}] \rightarrow B''[\frac{1}{p}]$ is a closed immersion, $B''[\frac{1}{p}]$ is unramified over $B[\frac{1}{p}]$. Since $B/pB \cong B''/pB''$,  $B''$ is unramified over $B$ and so it is \'etale over $B$. Similarly, $B''$ is \'etale over $B'$.

We denote $B^h \coloneqq \varinjlim B'$ where the colimit is taken over the diagrams as above. Note that $B^h$ is the henselization of $B$ with respect to $(p)$ as defined in  \cite{raynaud-henselien}.

\begin{proposition}\label{prop-i&j}
	The functor $G \mapsto (i^*G, j^*G, \alpha)$	 is an equivalence of categories between the category of sheaves on ${\Xbar}_{\syn-\et}$ to the category of triples $(F, H, \alpha)$ where $F$ is a sheaf on $({\wh{\Xbar}})_{\syn-\et}$, $H$ a sheaf on $(X_{\Rbar[1/p]})_{\et}$, and $\alpha\colon F \rightarrow i^*j_*H$ a morphism.
	
	Furthermore, $i_*$ is exact, and the functor $j_!$ from the category of sheaves on $(X_{\Rbar[\frac 1 p]})_{\et}$   to that on ${\Xbar}_{\syn-\et}$ given by $j_!(G) \coloneqq (0, j_*G, 0)$ is exact. For any sheaf $G$ on ${\Xbar}_{\syn-\et}$, the sequence
	\[
	0 \rightarrow j_! j^* G \rightarrow G \rightarrow i_* i^* G \rightarrow 0
	\]
	is exact. 
\end{proposition}

\begin{proof}
	Note that for any sheaf $F$ on $({\wh{\Xbar}})_{\syn-\et}$, we have $i^*i_* F \cong F$ and $j^*i_*F \cong 0$. And for any sheaf $H$ on $(X_{\Rbar[\frac 1 p]})_{\et}$, $j^*j_*H \cong H$. 
	
	For a sheaf $G$ on ${\Xbar}_{\syn-\et}$, consider the commutative diagram
	$$\xymatrix{G \ar[r] \ar[d] 
		& i_*i^*G \ar[d]\\
		j_*j^*G \ar[r]
		& i_*i^*j_*j^*G}$$	
	
	\noindent We claim that it is cartesian. Consider the natural morphism $G \rightarrow i_*i^*G\times_{(i_*i^*j_*j^*G)} j_*j^*G$. Since $i^*$ and $j^*$ are exact,
	\[
	i^*(i_*i^*G\times_{(i_*i^*j_*j^*G)} j_*j^*G) \cong i^*i_*i^*G\times_{(i^*i_*i^*j_*j^*G)} i^*j_*j^*G \cong i^*G\times_{i^*j_*j^*G}i^*j_*j^*G \cong i^*G,
	\]
	\[
	j^*(i_*i^*G\times_{(i_*i^*j_*j^*G)} j_*j^*G) \cong j^*i_*i^*G\times_{(j^*i_*i^*j_*j^*G)} j^*j_*j^*G \cong 0\times_0 j^*G \cong j^*G.
	\]
	Thus, the morphism $G \rightarrow i_*i^*G\times_{(i_*i^*j_*j^*G)} j_*j^*G$ is an isomorphism. This proves that the diagram is cartesian, which implies the first statement. 
	
	The exactness of $i_*$ follows from its definition. The remaining statements are direct consequences. 
\end{proof}

We now construct a sheaf $\mathcal{S}_n^{[r]}$ on ${\Xbar}_{\syn-\et}$ such that $i^*\mathcal{S}_n^{[r]} = S_n^{[r]}$ and $j^*\mathcal{S}_n^{[r]} = \Z/p^n\Z(r)$. By Proposition \ref{prop-i&j}, this is equivalent to defining a map $\alpha\colon {S_n^{[r]}} \rightarrow i^*j_* \Z/p^n\Z(r)$. For a map $\tilde{f}\colon \wh{A} \rightarrow \mathcal{B}$ with $f\colon A \rightarrow B$ as in Lemma \ref{lem-Zar-et}, we need to define a functorial map ${S_n^{[r]}}(\mathcal{B}) \rightarrow\Z/p^n\Z(r)(B^h[\frac{1}{p}])$. If we write $\overline{R} = \varinjlim_{\eta} R_{\eta}$ as above, then there exists an index $\lambda$, finitely generated smooth $R_{\lambda}$-algebra $A_{\lambda}$, and finitely generated $A_{\lambda}$-algebra $B_{\lambda}$ with $f_{\lambda}\colon A_{\lambda} \rightarrow B_{\lambda}$ such that $\displaystyle A \cong \varinjlim_{\eta \geq \lambda} A_{\lambda}\otimes_{R_{\lambda}} R_{\eta}$ and $\displaystyle B \cong \varinjlim_{\eta \geq \lambda} B_{\lambda}\otimes_{R_{\lambda}} R_{\eta}$, and $f_{\lambda}$ is syntomic and quasi-finite having \'etale generic fiber inducing $f$. Moreover, by taking large enough $R_{\lambda}$, we can assume that $R_{\lambda}$ contains a primitive $p^n$-th root of $1$.

Note that $S_n^{[r]}(\mathcal{B}) = S_n^{[r]}(\mathcal{B}/p^{n}) = S_n^{[r]}(B/p^{n}) = \varinjlim_{\eta} S_n^{[r]}(B_{\eta}/p^{n})$, and $\Z/p^n\Z(r)(B^h[\frac{1}{p}]) = \varinjlim_{\eta}\Z/p^n\Z(r)(B_{\eta}^h[\frac{1}{p}])$. Denote $A' = A_{\lambda}$ and $B' = B_{\lambda}$. Since $\pi_0(\mathrm{Spec}((B')^h[\frac{1}{p}])) = \pi_0(\mathrm{Spec}(\wh{B}'[\frac{1}{p}]))$ by \cite[Cor. 4.4]{Ferrand-Raynaud}, it suffices to construct a functorial morphism $S_n^{[r]}(B'/p^n) \rightarrow\Z/p^n\Z(r)(\wh{B}'[\frac{1}{p}])$.

Since $\wh{B}'[\frac{1}{p}]$ is \'etale over $\widehat{A'}[\frac{1}{p}]$ and Noetherian, we can write $\wh{B}'[\frac{1}{p}] \cong \prod_{a=1}^{\ell} B_a'$ where each $B_a'$ is an integral domain which is regular and thus normal. Each $B_a'$ is an affinoid algebra over $R[\frac{1}{p}]$, and we denote by $\mathring{B_a'}$ the subring of $B_a'$ consisting of elements with spectral norm $\leq 1$. Note that $\mathring{B_a'}$ is integrally closed in $B_a'$ and thus is a normal domain. We have a natural map $S_n^{[r]}(B'/p^{n}) \rightarrow \prod_{a=1}^{\ell} S_n^{[r]}(\mathring{B_a'}/p^{n})$. On the other hand, $\Z/p^n\Z(r)(\wh{B}'[\frac{1}{p}]) \cong \prod_{a=1}^{\ell} (\mathbf{\mu}_{p^n}(\overline{R}[\frac{1}{p}])^{\otimes r})_a$. 

Let $E = \mathrm{Frac}(A')$ equipped with the $p$-adic valuation, and let $C$ be the $p$-adic completion of the algebraic closure of $\widehat{E}$. For each $a = 1, \ldots, \ell$, choose an embedding $B_a' \rightarrow C$, which induces an embedding $\mathring{B_a'} \rightarrow \mathcal{O}_C$. This induces a map $S_n^{[r]}(\mathring{B_a'}/p^{n}) \rightarrow S_n^{[r]}(\mathcal{O}_C/p^{n})$. Note that $C$ contains $\overline{R}[\frac{1}{p}]$. We claim that $S_n^{[r]}(\mathcal{O}_C/p^{n}) = \mathbf{\mu}_{p^n}(\overline{R}[\frac{1}{p}])^{\otimes r}$. Indeed, we have $S_n^{[r]} (\cO_C)= \rH^0 _\cris (\cO_{C, n}/ R_n, J^{[r]}_{\cO_{C, n}/ R_n})^{\varphi_r = 1, \nabla= 0}$. Since $\varphi\colon \cO_C/ p  \to \cO_C/ p$ is surjective, we have by Proposition \ref{prop-calculateH0} and Lemma \ref{lem-toOAcris}
 $$S_n^{[r]} (\cO_C)= \rH^0 _\cris (\cO_{C, n}/ R_n, J^{[r]}_{\cO_{C, n}/ R_n} )^{\varphi_r = 1, \nabla= 0}= (\rF^r A_{\cris, n} (\cO_C))^{\varphi _r =1}=  \mathbf{\mu}_{p^n}(C)^{\otimes r}$$
where the last equality uses the well-known result in  \cite{Fontaine-BTrings}. We obtain the desired morphism $\alpha$ in this way. 

\begin{proposition}
	The natural map $\rH^*(\overline{X}_{\syn-\et}, \mathcal{S}_n^{[r]}) \rightarrow \rH^*(\widehat{\overline{X}}_{\syn-\et}, S_n^{[r]})$ induced by $i$ is an isomorphism.	
\end{proposition}

\begin{proof}
	By Proposition \ref{prop-i&j}, we need to show $\rH^*(\overline{X}_{\syn-\et}, j_! \mathbb{Z}/p^n\mathbb{Z}) = 0$. We have an exact sequence
	\[
	0 \rightarrow j_!(\mathbb{Z}/p^n\mathbb{Z}) \rightarrow \mathbb{Z}/p^n\mathbb{Z} \rightarrow i_*(\mathbb{Z}/p^n\mathbb{Z}) \rightarrow 0.
	\]	
	Since $i_*$ is exact, it suffices to show that the natural map
	\[
	\rH^*(\overline{X}_{\syn-\et}, \mathbb{Z}/p^n\mathbb{Z}) \rightarrow \rH^*({\widehat{\overline{X}}_{\syn-\et}, \mathbb{Z}/p^n\mathbb{Z}})
	\]
	is an isomorphism. 

By \cite[\href{https://stacks.math.columbia.edu/tag/0DDU}{Tag 0DDU}]{stacks-project}, the natural map 
	\[
	\rH^*(\overline{X}_{\et}, \mathbb{Z}/p^n\Z) \rightarrow \rH^*(\overline{X}_{\syn-\et}, \mathbb{Z}/p^n\Z)
	\] 
	induced by the Leray spectral sequence for $\overline{X}_{\syn-\et} \rightarrow \overline{X}_{\et}$ is an isomorphism since the syntomic-\'etale topology on $\overline{X}_{\syn-\et}$ is finer than the \'etale topology and coarser than the flat topology. Similarly, we have $\rH^*(\widehat{\overline{X}}_{\et}, \mathbb{Z}/p^n\Z) \cong \rH^*(\widehat{\overline{X}}_{\syn-\et}, \mathbb{Z}/p^n\Z)$. On the other hand, the natural map $\rH^*(\overline{X}_{\et}, \mathbb{Z}/p^n\mathbb{Z}) \rightarrow \rH^*({\widehat{\overline{X}}_{\et}, \mathbb{Z}/p^n\mathbb{Z}})$ is an isomorphism by the proper base change theorem.  
\end{proof}

By the previous proposition, we obtain a map 
\[
\rH^*_{\syn}(\overline{X}_n, S_n^{[r]}) \rightarrow \rH^*_{\et}(X_{\overline{R}[\frac{1}{p}]}, \mathbb{Z}/p^n\mathbb{Z}(r))
\]  
which is compatible with $G_R$-actions.

\subsection{Proof of the main theorem (Theorem \ref{thm-main})} 

In summary, for $i \leq r \leq p-2$, we have constructed a map 
$$\iota\colon T_\cris (\rH ^i_\cris (X_n/ R_n))(r) \to \rH^i_\syn (\Xbar_n, S^{[r]}_{n}) \to \rH ^i _\et (X_{\Rbar[\frac 1 p]}, \Z/ p^n \Z(r))$$ 
which is compatible with $G_R$-actions. On the other hand, consider the map $b_g\colon R \to R_g$ (which factors through $\wt R$) together with a choice of an extension $\overline{b_g}\colon \Rbar \to \overline{R_g}$. Setting $X_g \coloneqq X \times_R R_g$ (viewed over $R_g$), we similarly obtain a $G_{R_g}$-equivariant map 
$$\iota_g\colon  T_\cris (\rH^i_\cris (X_{g, n}/ R_{g, n}))(r) \to \rH^i_\syn (\Xbar_{g, n}, S^{[r]}_{g, n}) \to \rH^i _\et (X_{\overline{R_g}[\frac{1}{p}]}, \Z/ p^n \Z(r)). $$
The map $\Rbar \to \overline{R_g}$ induces a map of topoi $i\colon (\Xbar_{g, n})_\syn  \to (\Xbar_n)_\syn$, and we have a natural map $\cJ^{[r]}_{\Xbar_n/ \wt R_n} \to i _* \cJ^{[r]}_{\Xbar_{g, n}/ R_{g, n}}$. By a similar argument as in the proof of Theorem \ref{thm-cons-phi} (2), this map is compatible with $\varphi_r$. Therefore, we obtain a natural map $\rH^i_\syn (\Xbar_n, S_n^{[r]})\to \rH^i_\syn (\Xbar_{g, n}, S_{g, n}^{[r]})$. The construction of the comparison map is functorial, and we therefore have the following commutative diagram: 
$$\xymatrix{ T_\cris (\rH ^i_\cris (X_n/ R_n))(r) \ar[r]^{\iota}\ar[d] &  \rH ^i _\et (X_{\Rbar[\frac 1 p]}, \Z/ p^n \Z (r))\ar[d]\\ T_\cris (\rH ^i_\cris (X_{g, n}/ R_{g, n}))(r) \ar[r]^-{\iota_g} &  \rH ^i _\et (X_{\overline{R_g}[\frac{1}{p}]}, \Z/ p^n \Z (r)) }$$

Note that our construction for $\iota_g$ is the same as that of \cite{FontaineMessing}, and so $\iota_g$ is an isomorphism. By Theorem \ref{thm-cry-BC} and Corollary \ref{cor-abstract-comp}, the left vertical map is an isomorphism.  As the right vertical is also an isomorphism by the smooth and proper base change for \'etale cohomology, we conclude that $\iota$ is an isomorphism. This completes the proof of the main theorem.

\section{Applications} 

We study some applications of Fontaine-Messing theory over $R$ analogous to those in the original setting of Fontaine--Messing (see \cite[\S I 3, 4]{FontaineMessing}).

\subsection{Crystalline property of $\rH^i_\et(X_{\Rbar[\frac 1 p]}, \bbQ_p)$} 

As before, we consider a smooth proper scheme $X$ over $R$. Let  $\rH^i_\et(X_{\Rbar[\frac 1 p]}, \Z_p) \coloneqq \varprojlim_n \rH^i_\et(X_{\Rbar[\frac 1 p]}, \Z/ p ^n \Z) $, $T^i \coloneqq \rH^i_\et(X_{\Rbar[\frac 1 p]}, \Z_p)/ \text{(torsion)}$ and 
$\rH^i_\et(X_{\Rbar[\frac 1 p]}, \bbQ_p) \coloneqq \rH^i_\et(X_{\Rbar[\frac 1 p]}, \Z_p)\otimes_{\Z_p} \bbQ_p $.  Let $\rH ^i _\cris (X/ R) \coloneqq  \varprojlim_n \rH^i_\cris (X_n /R_n)$ and $M^i \coloneqq \rH ^i _\cris (X/ R) / \text{(torsion)}$. Recall that we say a finite continuous $\bbQ_p$-representation $V$ of $G_R$ is \textit{crystalline} if the natural injective map 
\[
\alpha_{\cris}\colon D_{\cris}(V)\otimes_{R[\frac{1}{p}]} \cO B_{\cris}(R) \rightarrow V \otimes_{\bbQ_p} \cO B_{\cris}(R)
\]
is an isomorphism, where $\cO B_\cris(R) \coloneqq \cO A _\cris (R) [\frac 1 p, \frac 1 \beta]$ with $\beta= \log([\tilde \epsilon])$ as defined before Remark \ref{rem:indep_r} and $D_{\cris}(V) \coloneqq (V\otimes_{\bbQ_p} \cO B_{\cris}(R))^{G_R}$ (see \cite[Sec. 8]{brinon-relative}). In this subsection, $\cO A_{\cris}(R)$ is an $R$-algebra via $\lambda_2$.

\begin{theorem}\label{Thm-cryofHet}
For $i \leq p -2$, the natural map
\[
\rH^i_\et(X_{\Rbar[\frac 1 p]}, \Z_p)\otimes_{\Z_p} \cO A_{\cris}(R)[\frac{1}{p}] \rightarrow \rH^i_{\cris}(X/R)\otimes_R \cO A_{\cris}(R)[\frac{1}{p}]
\] 
is injective, and its cokernel is killed by $\beta^i$. In particular, $\rH^i_\et(X_{\Rbar[\frac 1 p]}, \bbQ_p)$ is a crystalline representation of $G_R$. 
\end{theorem}

\begin{proof}
By our main theorem with $r = i$ (Theorem \ref{thm-main} and Theorem \ref{thm-Faltings2}), we have
$$T^i(i) = T_{\cris}(M^i)(i) \coloneqq \rF^i  (\cO A _\cris (R) \otimes_{R}M^i) ^{\varphi_i  =1, \nabla = 0}  \subset \cO A _\cris (R)\otimes_{R}M^i,$$
and $M^i$ is finite $R$-free with the same type as $T^i (i)$. Choose an $R$-basis $\{e_1 , \dots, e_m\}$ of $M^i$ and a $\Z_p$-basis $\{f_1 , \dots, f_m\}$ of $T^i$. Let $A = (a_{ij})$ be the $m \times m$ matrix with entries in $\cO A_\cris(R)$ such that $(f_1, \ldots, f_m) = (e_1, \ldots, e_m )A$. We claim that there exists an $m \times m$ matrix $B = (b_{ij})$ with entries in $\cO A_\cris(R)[\frac 1 p]$ such that $AB = \beta ^i I_m$.  

To prove the claim, we first project $\Mbar _R \coloneqq \cO A_\cris (R) \otimes_{\lambda_2, R} M^i $ to $\Mbar \coloneqq A_\cris (R) \otimes_{\lambda_1, R} M $ via $\pi_M$ as in Lemma \ref{lem-Tcris-compatible}. By the same argument as in the proof of Lemma \ref{lem-Tcris-compatible}, we obtain the following commutative diagram: 
$$\xymatrix{T ^i(i)\ar@{=}[d] \ar@{^{(}->}[r] &  \Mbar_R^{\nabla = 0}  \ar@{^{(}->}[r]\ar[rd]^{\cong} & \Mbar _R \ar@[->>][d]^{\pi_M} \\  T ^i(i) \ar@{^{(}->}[rr] &  & \Mbar } $$
Moreover, $\Mbar _R = \Mbar_R^{\nabla =0} \otimes_{A_\cris (R)} \cO A_\cris (R)$. So it suffices to prove the claim after replacing $e_j$'s and $f_j$'s by $\pi_M (e_j)$'s and $\pi_M(f_j)$'s respectively, and $A$ by $\pi(A)$ (the corresponding matrix with entries in $A_{\cris}(R)$). 

We now have
$$T^i (i) \cong \rF^i (A_\cris (R) \otimes _R M^i) ^{\varphi_i =1} \subset A_\cris (R) \otimes_R  M ^i$$ 
with $\Z_p$-basis $\{f_j\}$ of $T^i$ and an $R$-basis $\{e_j\}$ of $M^i$ with $(f_1, \ldots , f_m) = (e_1, \dots , e_m)A$. It suffices to show there exists an $m \times m$ matrix $B$ with entries in $A_\cris (R)[\frac 1 p]$ such that $AB = \beta ^i I _m$. By Remark \ref{rem-adpated}, we may assume $\{e_1, \ldots, e_m\}$ is an adapted basis, so that $e_j \in \rF^{r_j} (M^i) \setminus \rF^{r_j+1} (M^i)$ for $ r_j\leq r_{j +1}\leq i $. Then the tuple $(M^i, \rF^j (M^i), \varphi _j)$ corresponds to an invertible matrix $D \in {\mathrm{GL}}_m (R)$ via $(\varphi_{r_1} (e_1), \ldots, \varphi_{r_m} (e_m)) = (e_1, \ldots, e_m )D$. Since 
$$(f_1, \ldots, f_m) = (e_1, \ldots, e_m) A \subset \rF^i (A_\cris(R) \otimes_R M) ^{\varphi_i =1},$$
we conclude that $D\varphi_{i-\vec r}(A) = A$, where $\varphi_{ i-\vec r} (A)$ has $(k, \ell)$-entry $\varphi_{i - r_k} (a_{k \ell})$ with $a_{k\ell}\in \rF^{i-r_k}A_\cris (R)$ being $(k, \ell)$-entry of $A$. 

Consider $(N, \rF^j (N), \varphi_j)$ defined as follows: $N$ is a finite free $R$-module with basis $\{e_1^\vee , \dots, e_m^\vee\}$, and $\rF^j (N)$'s are direct summands of $N$ as $R$-modules such that $e_j^\vee \in \rF^{i - r_j} (N) \setminus \rF ^{i - r_j+1}(N)$. Define $\varphi_j $-structure via $(\varphi_{i -r_1} (e^\vee_1), \ldots, \varphi_{i - r_m} (e^\vee_m)) = (e_1^\vee, \ldots, e_m ^\vee) (D^T)^{-1}$. Then $N/p^n$ is an object of $\MF^{[0, i]}(R)$ for each $n$.  In particular, $\rF^i(N \otimes_R A_\cris (R))^{\varphi_i= 1}$ is finite $\Z_p$-free with the same type as $N$ by the proof of Theorem \ref{thm-Faltings2} (note that the proof does not need the existence of $\nabla$). Similarly as above, we then have a matrix $Z = (z_{k \ell})_{m \times m}$ with $z_{k\ell} \in \rF ^{r_k}  A_\cris (R)$ such that $(D^{-1})^ T\varphi_{\vec r} (Z)= Z$, where $(k, \ell)$-entry of $\varphi_{\vec r} (Z)$ is $\varphi_{r_k} (z_{k\ell})$. Now consider $Z^T A = (y_{k\ell})_{m \times m}$.  We have $y_{k\ell} \in \rF^ i A_\cris (R)$, and 
$$\varphi (Z^T A)= \varphi (Z^T)\varphi (A)=  Z^T D \Lambda' \Lambda D^{-1} A = p ^i Z^T A$$ 
where $\Lambda'$ and $\Lambda$ are diagonal matrices whose $(k, k)$-entries are $p^{r_k}$ and $p ^{i-r_k}$ respectively. Since $\varphi(\beta) = p\beta$ and $y_{k\ell} \in \rF^ i A_\cris (R)$, we have $y_{k\ell} = \beta^i w_{k\ell}$ for some $w_{k\ell} \in \mathbb{Q}_p$ by \cite[Cor. 6.2.20]{brinon-relative}. Let $W$ be the $m\times m$ matrix $(w_{k\ell})$ so that $Z^T A = \beta^i W$. 

For $M_g \coloneqq M^i\otimes_{R, b_g} R_g$ and $N_g \coloneqq N\otimes_{R, b_g} R_g$, $M_g/p^n$ and $N_g/p^n$ are objects in $\MF^{[0, i]}(R_g)$ for each $n \geq 1$. Moreover, $T_{\cris}(M^i) \cong T_{\cris}(M_g)$ as $G_{R_g}$-representations by Corollary \ref{cor-abstract-comp}. By classical Fontaine--Laffaille theory, we have $T^i(i)\otimes_{\Z_p} B_{\cris}(R_g) \cong M_g\otimes_{R_g} B_{\cris}(R_g)$. In particular, $\det{A} \neq 0$. Similarly, $\det{Z} \neq 0$, and thus $\det{W} \neq 0$. Then we can let $B = W^{-1}Z^T$, which proves the first part of the statement.

Now the natural map
$$T^i (i) \otimes_{\Z_p} \cO B_\cris (R) \rightarrow M^i\otimes_R \cO B_\cris (R)$$ 
is an isomorphism. Hence, for $V = \rH^i_\et(X_{\Rbar[\frac 1 p]}, \bbQ_p)$, we have $D_{\cris}(V) = M^i[\frac{1}{p}]$ by \cite[Prop. 6.2.9]{brinon-relative}, and $\alpha_{\cris}$ is an isomorphism. 
\end{proof}

\subsection{A local invariant cycle theorem} 

We assume $k = \bar k$ in this subsection. Let $X_k \coloneqq X \times_{R} k$. By Theorem \ref{thm-main}, $M \coloneqq \rH ^i _\cris (X_n /R_n )$ is an object of $\MF^{[0, i]}_\nabla (R)$ if $i \leq p-2$. We write $\varphi  = \varphi _0\colon M \to M$. 

\begin{proposition}\label{prop-etclosedfiber}Assume that $k = \bar k$ and $i \leq p-2$. Then 
	$$\rH^i _\et (X_k, \Z/p ^n \Z ) \simeq \rH ^i _\cris (X_n /R_n )^{\varphi =1}\simeq \rH ^i _\et ( X_{\Rbar[\frac{1}{p}]}, \Z/ p ^n \Z ) ^{G_R}. $$
\end{proposition}

\begin{proof} Let $X_0 \coloneqq X \times_{R} W_n (k)$ and $\Mbar \coloneqq \rH^i _\cris (X_0 / W_n(k))$. By Proposition \ref{prop-multi-key}, the natural base change map induces a $\varphi$-compatible isomorphism $M \otimes _R W(k)\simeq \Mbar$. Let $M^0= R \otimes_{\Z_p} M^{\varphi =1}$ and $\Mbar^0 = W(k)\otimes_{\Z_p} \Mbar^{\varphi =1}$. Since $\varphi$ is $\nabla$-horizontal, we have $\nabla (M^{\varphi=1})= \{0\}$. Consider the natural map $\iota\colon  M^0 \to M$. If we define $\rF ^1 (M^0) = \{0\}$ and $\nabla$ on $M^0$ via $\nabla (M ^{\varphi =1})= \{0\}$, then $M^0 $ is an object in $\MF_{\nabla,\text{big}} (R)$ and $\iota$ is a morphism in $\MF_{\nabla,\text{big}} (R)$. In a similar manner, we define $\Mbar^0$ and a morphism $\bar \iota\colon \Mbar ^0 \to \Mbar$ in $\MF_{\text{big}}(W(k))$. We claim that $\Mbar^0$ is a direct summand of $\Mbar$ and hence $\bar \iota$ is an injection. To see this, note that by Fitting's Lemma, we can write $\Mbar \simeq \Mbar^* \oplus \Mbar'$ where $\varphi$ acts on $\Mbar^* \coloneqq \bigcap_n \varphi^n (\Mbar )$ bijectively and $\Mbar^{\varphi =1} \subset \Mbar ^* $. By a standard Frobenius descent, we deduce that $\Mbar^* \simeq W(k) \otimes _{\Z_p}\Mbar ^{\varphi = 1} = \Mbar ^{0}$.  
	
Let $\fm \coloneqq (t_1, \dots , t_d)\subset R$ be the maximal ideal, and let $\pi\colon M \onto \Mbar = M / \fm M$ denote the natural projection. We claim that the induced map $\pi\colon M^{\varphi=1} \to \Mbar ^{\varphi =1}$ is an isomorphism; note this implies that $\bar \iota \simeq \iota \otimes_{R} W(k)$ and that $\iota$ is a morphism in $\MF_{\nabla}(R)$. To prove the claim, first note that for $x,y \in M^{\varphi=1}$, if $\pi (x) = \pi (y)$ then $x -y \in \fm M $ and $\varphi (x-y) = x-y  $. So $x-y \in \varphi^j  (\fm)M$ for any $j > 0$, which implies $x-y= 0$. On the other hand, let $\bar x \in \Mbar^{\varphi =1}$ and choose $x \in M$ such that $\pi (x)= \bar x$. Consider $\wt x \coloneqq x+ \sum\limits_{j = 0}^\infty (\varphi^j (\varphi(x)-x))$. $\wt x $ converges in $M$ since $\varphi (x)-x \in \fm M$, and $\pi(\wt x) = \bar x$ with $\varphi (\wt x) = \wt x$.

We now show $\ker(\iota) = 0$ as follows. Consider the image $\im(\iota)$ in the abelian category $\MF_{\nabla}(R)$. Note that ${\mathrm Tor}_1^{R}(W(k), N) = 0$ for any object $N$ in $\MF_{\nabla}(R)$ since such an object is a direct sum $\oplus R/p^{n_i}R$. In particular, tensoring the exact sequence 
$$0 \rightarrow \ker(\iota) \rightarrow M^0 \rightarrow \im(\iota) \rightarrow 0$$ 
with $W(k)$ over $R$ gives the analogous exact sequence for $\bar \iota$. So $\im(\iota) \otimes_R W(k) \cong \im(\bar \iota)$, and $\ker(\iota) = 0$ by Nakayama's Lemma.

Thus, $M^0$ is a subobject of $M$. By Theorem \ref{thm-Faltings2}, $T_\cris(M^0)$ embeds into $T_\cris (M)$. Since the map $\varphi$ on $M^{\varphi=1}$ is trivial and $\nabla (M^{\varphi =1})= 0$, we conclude that $T_\cris (M^0)$ is a representation on which $G_R$ acts trivially, i.e.,   
	$M^{\varphi =1}= T_\cris (M^0) \inj  T_\cris (M) ^{G_R}$. Similarly, $T_\cris (\Mbar ^0) \inj T_\cris (\Mbar) ^{G_K}$. 
	
On the other hand, since the classical Fontaine--Laffaille theory is stable under subobjects, there exists a subobject $N \subset \Mbar$ such that 
	$T_\cris (N) \simeq T_\cris(\Mbar) ^{G_K}$. $T_\cris (\Mbar) ^{G_K}$ is the maximal sub-representation of $T_\cris (\Mbar)$ on which $G_K$-acts trivially, and $\Mbar^0$ is the maximal $W(k)$-submodule of $\Mbar$ on which $\varphi$ acts bijectively. Thus, $N = \Mbar^0$ and $T_\cris (\Mbar) ^{G_K} \simeq T_\cris (\Mbar^0) = \Mbar ^{\varphi =1}$. 
	
From above observations, we obtain a commutative diagram:
	$$
\xymatrix{
	M ^{\varphi =1} = T_\cris (M ^0) \ar@{^{(}->}[r]  \ar[d]^{\cong} & T_\cris (M)^{G_R} \ar@{^{(}->}[d]   \\
	\Mbar^{\varphi =1}=T_\cris (\Mbar ^0) \ar[r]^-{\cong}  &  T_\cris (\Mbar)^{G_K} }
	$$
	Thus, all maps in this diagram are isomorphisms. In particular,
	$$\rH ^i _\cris (X_0 /W_n (k))^{\varphi =1} \simeq \rH ^i _\cris (X_n /R_n )^{\varphi =1}\simeq \rH ^i _\et ( X_{\Rbar[\frac{1}{p}]}, \Z/ p ^n \Z ) ^{G_R}\simeq  \rH ^i _\et ( X_{\Kbar}, \Z/ p ^n \Z )^{G_K}. $$
The result now follows from the well-known fact $\rH ^i _\cris (X_0 /W_n (k)) ^{\varphi =1} \simeq \rH^i _\et (X_k, \Z/ p ^n \Z)$ (see the proof of \cite[Thm. 5.2]{IllusiedR}). 
\end{proof}

\subsection{Cohomological non-defectness}

Let $k$ be a perfect field of characteristic $p > 2$ as before. In this subsection, we consider a finite totally ramified extension $K$ over $W(k)[\frac{1}{p}]$ with ramification index $e$. Let $\cO_K$ denote the ring of integers of $K$. Let $X$ be a proper smooth scheme over $\cO_K$, and $X_k \coloneqq X \times_{\cO_K} k$. It is interesting to study the relationship between $\rH^i _\cris(X_k/W_n(k))$ and $\rH^i_\et (X_{\overline K}, \Z/p^n \Z)$. By \cite[Thm. 1.1 (ii)]{SBM}, we have
\[
\ell_{W(k)} (\rH^ i _\cris(X_k/W_n(k))) \geq \ell_{\Z_p} (\rH^i_\et (X_{\overline K}, \Z/p^n \Z)).
\]
It is known that the inequality above can be strict in general when $ei\geq p-1$ (see \S 2 \emph{loc. cit.}). Moreover, \emph{loc. cit.} also gives an example where $\rH^ i _\cris(X_k/W(k))_{\mathrm{tor}} \cong k\oplus k$ and $\rH^i_\et (X_{\overline K}, \Z_p)_{\mathrm{tor}} \cong \Z/p^2\Z$. This motivates us to give the following definition.

\begin{defn} We say $X$ has \emph{cohomological non-defect in degree $i$} if 
\[
\rH^i_\et (X_{\overline K}, \Z/p^n \Z)\otimes_{\Z_p} W(k) \cong \rH^ i _\cris(X_k/W_n(k))
\]
as $W(k)$-modules for each $n$. Otherwise, we say $X$ has \emph{cohomological defect in degree $i$}. 
\end{defn}

We say that $X$ admits a \emph{model} $\cX$ over $R=W(k)[\![t]\!]$ if there exists a smooth proper scheme $\cX$ over $R$ and a map $b\colon R \to \cO_K$ such that $X \simeq \cX \times_{R, b} \cO_K$. 

\begin{proposition} \label{prop-cohom-non-defect} If $X$ admits a model $\cX$ over $R$, then $X$ has cohomological non-defect in degree $i$ for all $i \leq p-2$.
\end{proposition}

\begin{proof}
By Theorems \ref{thm-main} and \ref{thm-Faltings2}, $\rH^i_\cris (\cX_n / R_n)$ has the same type as $\rH^i_\et( \cX_{\Rbar [\frac 1 p]}, \Z/ p^n \Z)$. By considering base change maps $b\colon R \rightarrow \cO_K$ and $\bar{b}\colon R \rightarrow W(k)$ (via crystalline base change  theorem,  cf. Theorem \ref{thm-cry-BC} and Example \ref{Ex-BC}) and from Proposition \ref{prop-alphabar-suj}, we conclude that $\rH^i_\et (X_{\overline K}, \Z/p^n \Z)$ and $\rH^i _\cris(X_k/W_n(k))$ have the same type. 
\end{proof}

\begin{corollary}If $X$ has cohomological defect in degree $i$ for $i\leq p-2$, then $X$ does not admit a model over $R$. 
\end{corollary}

\appendix

\section{Background on the Crystalline site}\label{appendix-crystalline-site}

In this appendix, we recall some basic facts about the crystalline site. The results are standard and recalled here only for the convenience of the reader.

We begin by recalling the definition of the big and small crystalline site. Let $(\Sigma,\cI,\gamma)$ denote the data of a scheme $\Sigma$ with a quasi-coherent divided power ideal $(\cI,\gamma)$. Let $X$ be a $\Sigma$-scheme such that $p$ is locally nilpotent on $X$, and that the divided powers $\gamma$ extend to $X$. We shall moreover assume that the divided powers $\gamma$ extend to all $X$-schemes. The latter condition always holds when the ideal $\cI$ is principal, and this will be the case in all our applications.

\begin{defn}\label{defn:cris-site}
The big (resp. small) crystalline site $\CRIS(X/\Sigma)$ (resp. $\cris(X/\Sigma)$) is defined as follows (cf. \cite{BBM}, 1.1.1):
\begin{enumerate}
	\item[(1)] The objects of $\CRIS(X/\Sigma)$ are quadruples $(U,T,i,\delta)$ where $U$ is a $X$-scheme, $T$ is a $\Sigma$-scheme with $p$ locally nilpotent on $T$, $i\colon U \hookrightarrow T$ is a closed $\Sigma$-immersion, and $\delta$ is a divided power structure on the ideal of definition of $i$ compatible with $\gamma$. In particular, $i$ is a nil-immersion. Below, we shall denote such data simply by $(U,T)$. Morphisms are given by morphisms of the underlying schemes compatible with all of the above structures.
	\item[(2)] A morphism $(U',T') \rightarrow (U,T)$ is cartesian if the natural map $U' \rightarrow U \times_{T} T'$ is an isomorphism.
	\item[(3)] A morphism in $\CRIS(X/\Sigma)$ is an open immersion if the morphism is cartesian and the induced map $T' \rightarrow T$ is an open immersion.
	\item[(4)] The topology on $\CRIS(X/\Sigma)$ is generated by surjective families of open immersions.
	\item[(5)] The small site is defined to be the full subcategory of the big site where $U \rightarrow X$ is an open immersion, with the induced topology.
\end{enumerate}
\end{defn}

We now recall the standard base change theorem for crystalline cohomology. We shall denote by $(X/\Sigma)_{\cris}$ the topos of abelian sheaves corresponding to the small crystalline site. We will work with the small site in this section, but the usual comparison theorem between the big and small sites allows one to conclude analogous results for the big site. Let $\cO_{X/\Sigma}$ denote the usual crystalline structure sheaf whose section over $(U,T)$ is given by $\cO_T(T)$. We have natural morphisms of topoi:
$$(X/\Sigma)_{\cris} \xrightarrow{f_{X/\Sigma}} X_{\zar} \xrightarrow{f} \Sigma_{\zar}.$$
Here $X_{\zar}$ and $\Sigma_{\zar}$ denote corresponding topoi of abelian sheaves on the (small) zariski site.
Let $f_{X}$ denote the composite. The morphism $f_{X/\Sigma}$ comes equipped with a natural splitting
$i_{X/\Sigma}\colon X_{\zar} \rightarrow (X/\Sigma)_{\cris}$. By abuse of notation, let $\cO_{X} \coloneqq i_{X/\Sigma,*}\cO_{X}$. We have a natural exact sequence of sheaves of $\cO_{X/\Sigma}$-modules:
$$0 \rightarrow J_{X/\Sigma} \rightarrow \cO_{X/\Sigma} \rightarrow \cO_{X} \rightarrow 0.$$
We note that $J_{X/\Sigma}$ is a divided power ideal, and denote by $J^{[m]}_{X/\Sigma}$ the corresponding ideal of $m$-th divided powers.

Let $u\colon (\Sigma', I', \gamma') \rightarrow (\Sigma, I, \gamma)$ be a PD-morphism, and set $X' \coloneqq X \times_{\Sigma} \Sigma'$. The cartesian diagram of schemes
$$
\xymatrix{
X' \ar[r]^{u'} \ar[d]^{f'} & X \ar[d]^{f} \\
\Sigma' \ar[r]^{u} & \Sigma }
$$
gives rise to a canonical base change morphism
$$ B_{X/X'}\colon \bbL u^{*}\bbR f_{*}E \rightarrow  \bbR f'_{*}\bbL u^{'*}(E)$$
in $\rD(\cO_{\Sigma'})$ for any $E \in \rD(f^{-1}\cO_{\Sigma})$. Here $\rD(\cO_{\Sigma'})$ and $\rD(f^{-1}\cO_{\Sigma})$ denote the corresponding derived categories of sheaves in the Zariski topology.

By \cite[V Cor. 2.3.4, 2.3.7]{Berthelot}, there is a natural commutative diagram
\begin{equation}\label{compdiff}
\xymatrix{
\bbR f_{X/\Sigma,*}(J_{X/\Sigma}^{[m]}) \ar[r] \ar[d] & \tau_{\leq m}\Omega^{\bullet}_{X/\Sigma} \ar[d] \\
\bbR u'_{*} (\bbR f_{X'/\Sigma',*}(J_{X'/\Sigma'}^{[m]})) \ar[r] & \bbR u'_{*} (\tau_{\leq m}\Omega^{\bullet}_{X'/\Sigma'} )}
\end{equation}
where the horizontal maps are isomorphisms if $X/\Sigma$ is smooth. Applying the previous base change map to the left vertical map (and counit of adjunction $\bbL u^{'*} \bbR u'_{*}(E) \rightarrow E$) gives rise to a morphism
$$B_{X/X'}\colon \bbL u^{*} \bbR f_{*}(\bbR f_{X/\Sigma,*}(J_{X/\Sigma}^{[m]})) \rightarrow \bbR f'_{*} (\bbR f_{X'/\Sigma',*}(J_{X'/\Sigma'}^{[m]}))$$
in $\rD(\cO_{\Sigma'})$.

\begin{theorem}\label{thm-cry-BC}
With notations as above, suppose that $X \rightarrow \Sigma$ is a smooth morphism (quasi-separated and quasi-compact). Then the base change map $B_{X/X'}$ is an isomorphism.
\end{theorem}

\begin{proof}
We can proceed exactly as in \cite[Thm. 7.8]{crystal}. Suppose first that $X$ is a smooth affine scheme over $\Sigma$. The horizontal maps in \eqref{compdiff} are isomorphisms. Moreover, in this case, the base change map $B_{X/X'}$ identifies with the usual base change map for coherent sheaves in the Zariski topology:
$$ u^{*} f_{*} (\sigma_{\geq m}\Omega^{\bullet}_{X/\Sigma}) \rightarrow  f'_{*} (\sigma_{\geq m}\Omega^{\bullet}_{X'/\Sigma'}).$$
Here we used the compatibilty of Kh\"aler differential with base change to identify $u'^{*}\Omega^{\bullet}_{X/\Sigma}$ with
$\Omega^{\bullet}_{X'/\Sigma'}$ (and similarly for its truncations). Moreover, since $f$ is smooth, these complexes are flat and therefore derived pull-back maybe replaced by ordinary pull back. Similarly, since the morphisms are affine, we may replace derived push-forward by ordinary push-forward. On the other hand, the above base change map in the Zariski topology is known to be an isomorphism for affine morphisms. We may now proceed exactly as in the proof of \cite[Thm. 7.8]{crystal}. Briefly, we use a Cech cover and cohomological descent to reduce to the case of an affine morphism.
\end{proof}

\begin{example}\label{Ex-BC}
Let $k$ be a perfect field of characteristic $p$, $R = W(k)[\![t_{1},\ldots,t_{d}]\!]$, and $X$ be a proper smooth scheme over $R$. Let $R_n \coloneqq R/p^{n}R$ and $R'_n \coloneqq W_{n}(k)$. Then we have a natural isomorphism:
$$\bbR\Gamma((X_n/S_{n})_{\cris}, J^{[m]}_{X_n/S_n}) \otimes^{\bbL}_{R_{n}} R'_{n} \rightarrow \bbR\Gamma((X_n'/S_{n}')_{\cris}, J^{[m]}_{X_n'/S_n'}),$$
where $S_{n} = \Spec(R_{n})$ and $S'_{n} = \Spec(R'_{n})$.
\end{example}

We will also use a slightly more general form of the base change theorem, which we recall below for the convenience of the reader. Let $(\Sigma,J,\lambda)$ and $(\Sigma',J',\lambda') $ be divided power schemes with $p$ nilpotent. We shall assume as before that the divided power structure extends to all schemes over $\Sigma$ (resp. $\Sigma'$). 
Consider a commutative diagram:
\begin{equation*}
\xymatrix{
X' \ar[r]^{f'} \ar[d]^{g'} & Y' \ar[r] \ar[d]^g &     \Sigma ' \ar[d]^u \\
X \ar[r]^f & Y  \ar[r] & \Sigma .}
\end{equation*}
This gives rise to a commutative diagram of topoi:
$$
\xymatrix{
(X'/ \Sigma')_{\cris} \ar[r]^{f'_{\cris}} \ar[d]^{g'_{\cris}} & (Y'/\Sigma')_{\cris} \ar[d]^{g_{\cris}} \\
(X/\Sigma)_{\cris} \ar[r]^{f_{\cris}} & (Y/\Sigma)_{\cris}. }
$$
Let $E$ be a flat quasi-coherent crystal of $\cO_{X/\Sigma}$-modules. Suppose $f$ is a smooth morphism and the left square involving $g$ and $f$ is cartesian. Then by \cite[Thm. V 3.5.1]{Berthelot}, there is a natural base change isomorphism:
$$ \bbL g_{\cris}^*(\bbR f_{\cris,*}(E)) \rightarrow \bbR f'_{\cris,*}(\bbL g^{'*}_{\cris}(E)) $$
in the derived category of $\cO_{Y'/\Sigma'}$-modules. We will apply this base change morphism to the following diagram for a proper smooth scheme $X$ over $R$:
\begin{equation*}
\xymatrix{
	{\Xbar}_n \ar[r]^-{f'} \ar[d]^{g'} &  \Spec(\Rbar_n) \ar[r] \ar[d]^g &     \Spec(\wt R_n) \ar[d]^u \\
	\wt X_n \ar[r]^-f &  \Spec(\wt R_n) \ar[r]^-u  &  \Spec(\wt R_n)}
\end{equation*}
Here $\wt X_n$ is a proper smooth $\wt R_n$-scheme via $f$ and $\Xbar_n = \wt X_n \times_{\wt R} \Rbar$, and $u$ is the identity map. We then have an isomorphism
$$\bbL g_{\cris}^*(\bbR f_{\cris,*}(\cO_{{\wt X}_n/{\wt R}_n})) \rightarrow \bbR f'_{\cris,*}( \bbL{g'}^*_{\cris}(\cO_{{\wt X}_n/{\wt R_n} })) .$$
In particular, we have an isomorphism
$$\bbR f_{\cris,*}(\cO_{{\wt X}_n/{\wt R}_n}) \otimes^{\bbL}_{\wt R_n} \bbR g _{\cris, * }(\cO_{\Rbar_n/\wt R_n}) \rightarrow \bbR f'_{\cris,*}( \cO_{{\Xbar}_n/{\wt R_n} }) .$$
By  Remark \ref{rem-calculateH0-1}, we have $\bbR\Gamma((\Rbar_n/\wt R_n)_{\cris},\cO_{\Rbar_n/\wt R_n}) = A_{\cris,n}(R)$. Applying global section to the previous isomorphism gives an isomorphism
$$ A_{\cris,n}(R) \otimes^{\bbL}_{{\wt R}_n} \bbR\Gamma(({\wt X}_n/{\wt R}_n)_{\cris}, \cO_{\widetilde{X}_n/\wt R_n}) \simeq \bbR\Gamma((\Xbar_n/\wt R_n)_{\cris}, \cO_{{\Xbar}_n/{\wt R}_n}). $$ This gives rise to the following ${\mathrm Tor}$-spectral sequence (see \cite[\href{https://stacks.math.columbia.edu/tag/0662}{Tag 0662}]{stacks-project}):
$$
\rE_2^{m,q} \coloneqq {\mathrm Tor}^m_{\wt R_n}(A_{\cris,n}(R), \rH^q_{\cris}(\wt X_n/\wt R_n)) \Rightarrow
\rH^{m+q}_{\cris}(({\Xbar}_n/\wt R_n)_{\cris}, \cO_{{\Xbar}_n/{\wt R}_n}).
$$

\section{Background on the Syntomic and Crystalline-Syntomic sites}\label{appendix:syn-cris}

In this appendix, we recall some basic facts about the crystalline-syntomic site, and its comparison with the syntomic site. The results presented here are standard and well-known to the experts. In particular, we do not claim any originality. We recall them for readers' convenience due to the lack of an appropriate reference.

\begin{defn}\label{defn:syn-site}
Given a scheme $X$, let $\SYN(X)$ (resp. $\syn(X)$) denote the big (resp. small) syntomic site defined as follows (see \cite[Sec. 1]{Bauer}):
\begin{enumerate}
	\item[(1)] The objects of $\SYN(X)$ are schemes over $X$, and the topology is generated by surjective families of open immersions and finite surjective families of syntomic morphisms with source and target both affine schemes.
	\item[(2)] $\syn(X)$ is the full subcategory of $\SYN(X)$ with syntomic structure map and induced topology.
\end{enumerate}
\end{defn}

We now recall the (big) crystalline-syntomic site. Let $(\Sigma,\cI,\gamma)$ denote the data of a scheme $\Sigma$ with a quasi-coherent divided power ideal $(\cI,\gamma)$. Let $X$ be a $\Sigma$-scheme such that $p$ is locally nilpotent on $X$, and the divided powers $\gamma$ extend to $X$. We shall moreover assume that the divided powers $\gamma$ extend to all $X$-schemes. The latter condition always holds when $\cI$ is principal, which will be the case in all our applications.

\begin{defn}\label{defn:cris-syn-site}
Let $\CRIS(X/\Sigma)_{\SYN}$ denote the (big) crystalline syntomic site defined as follows (see \cite[Def. 1.8]{Bauer}):
\begin{enumerate}
	\item[(1)] The underlying category is $\CRIS(X/\Sigma)$. 
	\item[(2)] A morphism $(U',T') \rightarrow (U,T)$ in $\CRIS(X/\Sigma)$ is syntomic if it is cartesian and the map $T' \rightarrow T$ is syntomic.
	\item[(3)] The topology on $\CRIS(X/\Sigma)_{\SYN}$ is generated by surjective families of open immersions and finite surjective families of syntomic morphisms with base and target affine (i.e. objects $(U,T)$ with both $U$ and $T$ affine).
\end{enumerate}
\end{defn}

Let $X_{\SYN}$ (resp. $X_{\syn}$) denote the topos of sheaves associated to $\SYN(X)$ (resp. $\syn(X))$. Similarly, let $(X/\Sigma)_{\CRIS-\SYN}$ denote the topos associated to $\CRIS(X/\Sigma)_{\SYN}$.

\begin{rem}\label{rem:mor-crissyn-syn}
We have natural morphisms of topoi:
\begin{enumerate}
	\item[(1)] $u_{X/\Sigma}\colon (X/\Sigma)_{\CRIS-\SYN} \rightarrow X_{\SYN}$ where $u_{X/\Sigma}^{*}(\cG)(U,T) = \cG(U)$ for $\cG \in X_{\SYN}$ and $u_{X/\Sigma_{*}}(\cF)(U) = \Gamma(U, \cF|_{\CRIS(U/\Sigma)})$ where $U \in \SYN(X)$ and $\cF \in (X/\Sigma)_{\CRIS-\SYN}$ (see \cite[Prop. 1.10]{Bauer}). 
	\item[(2)] $\alpha\colon (X/\Sigma)_{\CRIS-\SYN} \rightarrow (X/\Sigma)_{\CRIS}$ where $\alpha_{*}$ just forgets the syntomic topology and the pull-back of a sheaf is given by taking the associated sheaf (see \cite[Rem. 1.12]{Bauer}).
\end{enumerate}
\end{rem}

{
\begin{rem}\label{rem:topoi-functoriality}
The topoi $X_{\SYN}$ is functorial in $X$, and similalry $(X/\Sigma)_{\CRIS-\SYN}$ is functorial in $X/\Sigma$. In particular, given a commutative diagram
$$
\xymatrix{
X' \ar[r]^{f} \ar[d] & X \ar[d]\\
\Sigma' \ar[r] & \Sigma }
$$
with $\Sigma$ and $\Sigma'$ as before (i.e. divided power schemes with $p$ nilpotent), we have morphisms of topoi 
$f_{\cris,\syn}\colon (X'/\Sigma')_{\CRIS-\SYN} \rightarrow (X/\Sigma)_{\CRIS-\SYN}$ and
$f_{\syn}\colon X'_{\SYN} \rightarrow X_{\SYN}$. The construction of $f_{\cris,\syn}$ is exactly analogous to the case of usual crystalline topoi. We also note that the resulting diagrams of topoi given by $u_{X/\Sigma}$ (resp. $\alpha$) commute.  
\end{rem}
}

\begin{rem}
\begin{enumerate}
\item[(1)] The site $\CRIS(X/\Sigma)_{\SYN}$ is naturally a ringed site where the structure sheaf $\cO_{S/\Sigma}$ is defined by
$$\Gamma((U,T), \cO_{X/\Sigma}) \coloneqq \Gamma(T,\cO_T).$$
\item[(2)] Similarly, we define the divided powers ideal sheaves $J^{[k]}_{X/\Sigma}$.
\end{enumerate}
\end{rem}

We now recall a result on the higher direct images of $\alpha$ due essentially to \cite[Prop. 1.1.19]{BBM}. Before stating the result, we make some preliminary remarks.

\begin{rem}(see \cite[1.1.3]{BBM})
Given a morphism of schemes $v\colon T' \rightarrow T$, there is a natural inverse image functor $v^{-1}\colon T_{\syn} \rightarrow T'_{\syn}$ (with a natural right adjoint denoted by $v_{*}$). Giving a sheaf on $\CRIS(X/\Sigma)_{\SYN}$ is equivalent to giving the the data of a sheaf $F_{(U,T)}$ on $T_{\syn}$ for each object $(U,T) \in \CRIS(X/\Sigma)_{\SYN}$ and for all morphisms $(u,v)\colon (U,T) \rightarrow (U',T')$, a morphism of sheaves $\rho_{(u,v)}\colon v^{-1}(F_{(U,T)}) \rightarrow F_{(U',T')}$ such that:
\begin{enumerate}
\item[(1)] $\rho_{(Id,Id)} = Id$,
\item[(2)] $ \rho_{(u \circ u',v \circ v')} = \rho_{(u',v')} \circ v'^{-1}(\rho_{(u,v)})$,
\item[(3)] If $(u,v)$ is syntomic, then $\rho_{(u,v)}$ is an isomorphism.
\end{enumerate}
\end{rem}

\begin{rem}(see \cite[1.1.18]{BBM})\label{rem:pullbacksheaf}
\begin{enumerate}
\item[(1)] Let $\alpha $ be as above. If $E$ is a sheaf on $\CRIS(X/\Sigma)$, then the pull back $\alpha^{*}(E)$ is by definition the sheaf associated to the presheaf $E$ on $\CRIS(X/\Sigma)_{\SYN}$.
\item[(2)] Suppose that $E$ is a sheaf of $\cO_{X/\Sigma}$-modules on $\CRIS(X/\Sigma)$ such that for all $(U,T)$, $E_{(U,T)}$ is a quasi-coherent $\cO_T$-module and for each syntomic morphism 
$(u,v)\colon (U',T') \rightarrow (U,T)$, the natural pull back map $v^{*}(E_{(U,T)}) \rightarrow E_{(U',T')}$ is an isomorphism of $\cO_{T'}$-modules. Then by part (1) above and flat descent, the presheaf $E$ is a sheaf in the syntomic topology, and the presheaves underlying $\alpha^{*}(E)$ and $E$ are the same, i.e. $\alpha^{*}(E) = E$.
\end{enumerate}
\end{rem}

\begin{proposition} (see \cite[Prop. 1.1.19]{BBM})\label{Prop-B.8}
Let $E$ be an $\cO_{X/\Sigma}$-module satisfying the conditions of Remark \ref{rem:pullbacksheaf} (2). Then for all $i > 0$, 
$\bbR^{i}\alpha_{*}(E) = 0$.
\end{proposition}

\begin{proof}
This is proved in \textit{loc. cit.} with the fppf and etale topologies on $\CRIS(X/\Sigma)$	instead of the syntomic topology. However, the same proof goes through verbatim in our setting. We briefly outline the steps:\\
(1) As in \textit{loc. cit.}, we are reduced to showing that the sheaf associated to the presheaf 
$$ (U,T) \mapsto \rH^{i}((X/\Sigma)_{\CRIS-\SYN}/\tilde{T}, E)$$
	vanishes for $i > 0$. Here $\tilde{T}$ is the sheaf represented by $(U,T)$, and  $(X/\Sigma)_{\CRIS-\SYN}/\tilde{T}$ is the corresponding topos of objects over $\tilde{T}$.\\
(2) Since this is a local question, we may reduce to the case where $T$ is affine.\\
(3) It suffices to show that the correspoding Cech cohomology groups vanish. \\
(4) The Cech complex of cochains along a syntomic cover $(U,T ) \rightarrow (U',T')$ can be computed by considering the analogs Cech complex of cochains for the syntomic cover $T' \rightarrow T$ with coefficients in $E_{(U,T)}$. Since syntomic morphisms are flat, the latter complex is exact by faithfully flat descent.
\end{proof}

\begin{corollary}\label{cor-alpha-vanish}
With notations as above, $\bbR^{i}\alpha_{*}(J^{[k]}_{X/\Sigma}) = 0$ for all $i>0$.
\end{corollary}

\begin{proof}
We must verify the hypotheses of Remark \ref{rem:pullbacksheaf} (2). For each $(U,T)$, $(J^{[k]}_{X/\Sigma})_{(U,T)}$ is by definition the divided powers of the ideal defining the closed immersion $U \hookrightarrow T$. Let $(u,v)\colon (U',T') \rightarrow (U,T)$ be a syntomic morphism. Since both $u$ and $v$ are flat, the ideal of definition of $U'$ in $T'$ is the pull back $v^{*}((J_{X/\Sigma})_{(U,T)})$, and $(J_{X/\Sigma})_{(U',T')} = (J_{X/\Sigma})_{(U,T)} \otimes _{\cO_T} \cO_{T'}$. Since divided power envelopes are compatible with flat base change by \cite[Prop. 3.21]{crystal}, we have
$ v^{*}((J^{[k]}_{X/\Sigma})_{(U,T)}) = (J^{[k]}_{X/\Sigma})_{(U',T')}$.
\end{proof}

\subsection*{Funding} This work was partially supported by National Science Foundation [DMS-1406926 to T.L., DMS-1502296 to D.P.]

\subsection*{Acknowledgements}

The first author would like to thank Fucheng Tan and Shizang Li for useful discussion and communications during the preparation of this paper. We also would like to thank anonymous referees for providing valuable comments and suggestions to improve the paper.

\bibliographystyle{agsm}

\begin{thebibliography}{0}

\bibitem[Bauer(1992)]{Bauer}
Bauer, W.
``On the conjecture of Birch and Swinnerton-Dyer for abelian varieties over function fields in characteristic $p > 0$.''
\textit{Invent. Math.} 108 (1992): 263--287.

\bibitem[Berthelot, Breen, and Messing(1982)]{BBM}
Berthelot, P., Breen, L., and Messing, W.
``Th\'eorie de Dieudonn\'e cristalline II.''
\textit{Lecture Notes in Math.} 930 (1982).
Berlin: Springer-Verlag.

\bibitem[Berthelot and Messing(1990)]{BM}
Berthelot, P. and Messing, W.
``Th\'eorie de Dieudonn\'e cristalline III.''
\textit{Progr. Math.} 86 (1990): 171--247.

\bibitem[Berthelot(1974)]{Berthelot}
Berthelot, P.
``Cohomologie cristalline des sch\'emas de caract\'eristique $p > 0$.''
\textit{Lecture Notes in Math.} 407 (1974).
Berlin-New York: Springer-Verlag.

\bibitem[Bhatt, Morrow, and Scholze(2018)]{SBM}
Bhatt, B., Morrow, M., and Scholze, P.
``Integral $p$-adic Hodge theory.''
\textit{Publ. Math. Inst. Hautes \'{E}tudes Sci.} 128 (2018): 219--397.

\bibitem[Berthelot and Ogus(1978)]{crystal}
Berthelot, P. and Ogus, A.
``Notes on crystalline cohomology.''
Princeton University Press, Princeton, NJ. (1978).

\bibitem[Breuil(2000)]{Breuil}
Breuil, C.
``Groupes $p$-divisibles, groupes finis et modules filtr\'es.''
\textit{Ann. of Math.} (2). 152 (2000): 489--549.

\bibitem[Brinon(2008)]{brinon-relative}
Brinon, O.
``Repr\'esentations $p$-adiques cristallines et de de rham dans le cas relatif.''
\textit{M\'em. Soc. Math. Fr.} 112 (2008).

\bibitem[Caruso(2008)]{CarusoInvent}
Caruso, X.
``Conjecture de l'inertie mod\'er\'ee de Serre.''
\textit{Invent. Math.} 171 (2008): 629--699.

\bibitem[Faltings(1989)]{FaltingsCcris}
Faltings, G.
``Crystalline cohomology and $p$-adic Galois representations.'' 
\textit{Algebraic analysis, geometry, and number theory.} (1989): 25--80.
Johns Hopkins University Press, Baltimore, MD. 

\bibitem[Faltings(1999)]{Faltings}
Faltings, G.
``Integral crystalline cohomology over very ramified valuation rings.''
\textit{J. Amer. Math. Soc.} 12 (1999): 117--144.

\bibitem[Fontaine and Laffaille(1982)]{Fonatine-Laffaille}
Fontaine, J.-M. and Laffaille, G.
``Construction de repr\'esentations $p$-adiques.''
\textit{Ann. Sci. \'Ec. Norm. Sup\'{e}r.} (4). 15 (1982): 547--608.

\bibitem[Fontaine and Messing(1987)]{FontaineMessing}
Fontaine, J.-M. and Messing, W.
``$p$-adic periods and $p$-adic \'etale cohomology,'' in Current trends in arithmetical algebraic geometry.
\textit{Contemp. Math.} 67 (1987): 179--207.

\bibitem[Fontaine(1982)]{Fontaine-BTrings}  
Fontaine, J.-M.
``Sur certains types de repr\'esentations $p$-adiques du groupe de Galois d’un corps local; construction d’un anneau de Barsotti-Tate.''
\textit{Ann. of Math.} (2). 115 (1982): 529--577.

\bibitem[Ferrand and Raynaud(1970)]{Ferrand-Raynaud}
Ferrand, D. and Raynaud, M.
``Fibres formelles d’un anneau local noeth\'erien.''
\textit{Ann. Sci. \'Ec. Norm. Sup\'{e}r.} (4). 3 (1970): 295--311.

\bibitem[Grothendieck(1967)]{EGAIV.4}
Grothendieck, A.
``\'El\'ements de g\'eom\'etrie alg\'ebrique IV. \'Etude locale des sch\'emas et des morphismes de sch\'emas IV.''
\textit{Publ. Math. Inst. Hautes \'{E}tudes Sci.} 32 (1967).

\bibitem[Illusie(1979)]{IllusiedR}
Illusie, L.
``Complexe de deRham-Witt et cohomologie cristalline.''
\textit{Ann. Sci. \'Ec. Norm. Sup\'{e}r.} (4). 12 (1979): 501--661.

\bibitem[Li and Liu(2023)]{li-liu-prism-deRham-cohom}
Li, S. and Liu, T.
``Comparison of prismatic cohomology and derived de Rham cohomology.''
\textit{J. Eur. Math. Soc. (JEMS)} (2023).

\bibitem[Moon(2021)]{moon-BTdeform}
Moon, Y.S.
``Extending $p$-divisible groups and Barsotti-Tate deformation ring in the relative case.''
\textit{Int. Math. Res. Not. IMRN.} 2021 (2021): 13182--13201.

\bibitem[Raynaud(1970)]{raynaud-henselien}
Raynaud, M.
``Anneaux locaux hens\'eliens.''
\textit{Lecture Notes in Math.} 169 (1970).
Berlin: Springer-Verlag.

\bibitem[Stacks Project(2023)]{stacks-project}
The Stacks Project Authors.
``Stacks Project.''
(http://stacks.math.columbia.edu).
 
\bibitem[Xu(2019)]{DXU}
Xu, D.
``Lifting the Cartier transform of Ogus-Vologodsky modulo $p^n$.''
\textit{M\'em. Soc. Math. Fr.} (N.S.) 163 (2019). 

\end{thebibliography}

\end{document}